%% file: main.tex
\documentclass{article}

\usepackage[final,nonatbib]{neurips_2021}

\usepackage{microtype}
\usepackage{graphicx}
\usepackage{subfigure}
\usepackage{booktabs}
\usepackage{amsmath}
\usepackage{mathtools}
\usepackage{amssymb}
\usepackage{amsthm}
\usepackage{tcolorbox}
\usepackage{thm-restate}
\usepackage{tcolorbox}
\usepackage{hyperref}
\tcbset{boxsep=0mm,boxrule=0pt,colframe=white,arc=0mm,left=0.5mm,right=0.5mm}
\usepackage{soul}
\usepackage{enumitem}
\usepackage{algorithm}
\usepackage{algorithmic}
\usepackage{wrapfig}

\input{defs}

\definecolor{bulgarianrose}{rgb}{0.28, 0.02, 0.03}

\title{On the Second-order Convergence Properties of Random Search Methods}

\author{%
  Aurelien Lucchi\thanks{Alphabetical ordering, all authors contributed equally.}
   \qquad
 Antonio Orvieto
 \qquad
  Adamos Solomou
 \vspace{3.5mm}\\
 Department of Computer Science\\
 ETH Zurich
}

\begin{document}

\maketitle

\begin{abstract}
We study the theoretical convergence properties of random-search methods when optimizing non-convex objective functions without having access to derivatives. We prove that standard random-search methods that do not rely on second-order information converge to a second-order stationary point. However, they suffer from an exponential complexity in terms of the input dimension of the problem. In order to address this issue, we propose a novel variant of random search that exploits negative curvature by only relying on function evaluations. We prove that this approach converges to a second-order stationary point at a much faster rate than vanilla methods: namely, the complexity in terms of the number of function evaluations is only linear in the problem dimension. We test our algorithm empirically and find good agreements with our theoretical results.
\end{abstract}

\input{01_introduction}
\input{02_related_work}

\input{03_analysis}

\input{04_experiments}

\input{05_conclusion}

\bibliography{main}
\bibliographystyle{plain}

\input{appendix}

\end{document}

%% file: defs.tex
\newcommand{\bigO}{\mathcal{O}}

\DeclareMathOperator*{\argmin}{arg\,min}

\newcommand{\e}{{\bf e}}
\newcommand{\g}{{\bf g}}

\newcommand{\s}{{\bf s}}

\renewcommand{\u}{{\bf u}}
\renewcommand{\v}{{\bf v}}
\newcommand{\vDelta}{{\bf \Delta}}

\newcommand{\x}{{\bf x}}
\newcommand{\y}{{\bf y}}

\newcommand{\vxi}{{\boldsymbol{\xi}}}

\newcommand{\DFPI}{\text{\tiny DFPI}}

\makeatletter
\newenvironment{proof*}[1][\proofname]{\par
  \pushQED{\qed}%
  \normalfont \partopsep=\z@skip \topsep=\z@skip
  \trivlist
  \item[\hskip\labelsep
        \itshape
    #1\@addpunct{.}]\ignorespaces
}{%
  \popQED\endtrivlist\@endpefalse
}
\makeatother

\def\Am{{\bf A}}

\def\Hm{{\bf H}}
\def\Im{{\bf I}}
\def\Mm{{\bf M}}

\def\Vm{{\bf V}}

\newcommand{\A}{{\mathcal A}}

\newcommand{\D}{{\mathcal D}}
\newcommand{\E}{{\mathbf E}}

\newcommand{\R}{{\mathbb{R}}}
\renewcommand{\S}{{\mathcal S}}

\newtheorem{theorem}{Theorem}
\newtheorem{lemma}[theorem]{Lemma}
\newtheorem{proposition}[theorem]{Proposition}

\newtheorem{assumption}{Assumption}

\newtheorem{remark}{Remark}

%% file: 01_introduction.tex

\setcounter{footnote}{0} 

\section{Introduction}
\vspace{-2mm}
We consider solving the non-convex optimization problem $\min_{\x\in\R^d} f(\x)$, where $f(\cdot)$ is differentiable but its derivatives are not directly accessible, or can only be approximated at a high computational cost. This setting recently gained attention in machine learning, in areas such as black-box adversarial attacks~\cite{chen2017zoo}, reinforcement learning~\cite{salimans2017evolution}, meta-learning~\cite{bergstra2012random}, online learning~\cite{bubeck2012regret}, and conditional optimization~\cite{yang2020devil}.

We focus our attention on a popular class of derivative-free methods known as random direct-search methods of directional type~\footnote{In this manuscript, we will use the terms ``random direct-search" and ``random search" interchangeably.}. These methods optimize $f(\cdot)$ by evaluating the objective function over a number
of (fixed or randomized) directions, to ensure descent using a sufficiently small stepsize.
Direct-search algorithms date to the 1960's, including e.g.~\cite{matyas1965random,nelder1965simplex}. More recent variants include deterministic direct search~\cite{conn2009introduction}, random direct search~(e.g.~\cite{vicente2013worst}, or the Stochastic Three Points (STP) method~\cite{bergou2019stochastic}), which randomly sample a direction and accept a step in this direction if it decreases the function $f(\cdot)$. As discussed in~\cite{lewis2000direct}, direct-search methods have remained popular over the years for a number of reasons, including their good performance and known global convergence guarantees~\cite{vicente2013worst}, as well as their straightforward implementation that makes them suitable for many problems. We refer the reader to~\cite{lewis2000direct,conn2009introduction} for a survey. 

In machine learning, objective functions of interest are often non-convex, which poses additional challenges due to the presence of saddle points and potentially suboptimal local minima~\cite{jin2019nonconvex}. Instead of aiming for a global minimizer, one often seeks a second-order stationary point (SOSP): i.e. a solution where the gradient vanishes and the Hessian is positive definite. Indeed, as shown by~\cite{choromanska2015loss,kawaguchi2016deep,ge2016matrix,ge2017no}, many machine learning problems have no spurious local minimizers, yet have many saddle points which yield suboptimal solutions and are often hard to escape from~\cite{du2017gradient}. While convergence to SOSPs and saddle escape times have been extensively studied in the context of gradient-based methods~\cite{jin2017escape,daneshmand2018escaping,carmon2018accelerated,xu2018first}, \textit{prior analyses for (random) direct search have mainly focused on convergence to first-order stationary points}~\cite{bergou2019stochastic, vicente2013worst}. One exception is the Approximate Hessian Direct Search (AHDS) method~\cite{gratton2016second}, that explicitly computes the Hessian of the objective to obtain second-order worst-case guarantees. However, computing or storing a full Hessian is prohibitively expensive in high dimensions.\looseness=-1

Towards a better understanding of the complexity of finding second-order stationary points with random search methods, we make the following contributions:

\begin{itemize}[leftmargin=*]
\setlength\itemsep{0.2em}
    \item We study the complexity of a simple random search similar to STP (Algorithm~\ref{algo:random_search}) to reach SOSPs. We find that the (worst-case) complexity requires a number of function evaluations that scales \textit{exponentially} in terms of the problem dimension $d$. As we will see, the exponential scaling is not an artefact of the analysis. This is intuitive, indeed, if we are at a saddle point where we just have one direction of negative curvature, finding good alignment with a \textit{random} direction becomes exponentially difficult as the dimension of the space increases.
    \item To solve this issue, we design a variant of random search (RSPI, Algorithm~\ref{algo:random_search_PI}) that, instead of randomly sampling directions from the sphere, relies on an approximate derivative-free power iteration routine to extract negative curvature direction candidates (unlike~\cite{gratton2016second} that requires an approximation to the full Hessian). This approach is inspired from recent work on gradient-based methods for finding SOSPs~\cite{carmon2018accelerated,xu2018first} and effectively decouples \textit{negative curvature estimation} from progress in the large gradient setting. As a result, RSPI does not suffer from the exponential scaling of the vanilla random search approach : we show that the overall complexity of finding a SOSP in terms of function evaluations becomes \textit{linear} in the problem dimension.\looseness=-1
    \item Finally, we verify our results empirically and compare our novel algorithm to standard random-search methods. Our results show improvements of RSPI both in terms of algorithm iterations and (crucially) wall-clock time.
\end{itemize}

%% file: 02_related_work.tex

\vspace{-2mm}
\section{Related work}
\label{sec:related_work}

\paragraph{Direct-search vs methods that approximate gradients}
A wide variety of algorithms enter the class of DFO methods. One common distinction is made according to whether or not the algorithm explicitly computes a gradient approximation. Direct-search (DS) and pattern-search (PS) methods only rely on function evaluations to validate a step along a direction, sampled according to some scheme. In contrast, a second group of methods explicitly compute a gradient approximation~\cite{nesterov2017random,ghadimi2013stochastic,chen2020accelerated}. Most methods in the latter group rely only on first-order information, except~\cite{ye2018hessian} that incorporates some second-order information to compute an estimate of the gradient. However, their approach focuses on convex optimization and therefore does not discuss second-order stationarity. Instead,~\cite{vlatakis2019efficiently} showed that $\mathcal{O}(\epsilon^{-2})$ approximate gradient computations are enough to reach an $\epsilon-$SOSP. Unlike DS, the gradient-free method analyzed in~\cite{vlatakis2019efficiently} computes an approximation of the gradient and, as a result, it matches the convergence rate guarantees of their exact gradient-based counterparts (up to constants). However, existing lower bounds clearly show that DS methods have a worst rate of convergence. Of special interest to us is the complexity w.r.t. the dimension of the problem, which is exponential, see paragraph ``Lower bounds" below.

\vspace{-2mm}
\paragraph{First-order guarantees of DS.} A recently proposed variant of DS is the stochastic-three-points~(STP) approach proposed in~\cite{bergou2019stochastic} that simply samples a direction at random on the sphere and accepts the step if it decreases the function. This is, in some sense, a simpler variant of direct search compared to adaptive methods such as~\cite{vicente2013worst} that also include an adaptive mechanism of the step-size selection to ensure convergence without having to rely on knowing certain quantities such as the smoothness constant of the objective. As shown in~\cite{bergou2019stochastic}, finding a point $\x$ such that $\|\nabla f(\x)\|\le\epsilon$ with STP requires $\mathcal{O}(d\epsilon^{-2})$ function evaluations. If instead the perturbation directions are sampled from a fixed basis and not from the unit sphere, the best known complexity increases to $\mathcal{O}(d^2\epsilon^{-2})$~\cite{vicente2013worst,bergou2019stochastic}. Ignoring the dependency on the dimension, these results match the iteration complexity of steepest descent~\cite{nesterov2018lectures}.

\vspace{-2mm}
\paragraph{Second-order guarantees for PS and DS.}
Second-order convergence guarantees have been developed in the context of (generalized) pattern search (GPS) methods, which share similarities with direct search methods. These methods sample update directions from a positive spanning set $\D$. The work of~\cite{abramson2005second} proved a type of ``pseudo-second-order" stationarity condition for GPS, where the Hessian is positive semidefinite in the directions of the basis $\D$ (but not with respect to all vectors in the space). A variant of GPS that explicitly constructs an approximate Hessian was shown by~\cite{abramson2014subclass} to converge to second-order stationary points under some assumptions on the quality of the approximation. The results discussed so far typically consider the properties of the limit points of the sequence of iterates and do not provide worst-case complexity rates of convergence. In contrast,~\cite{gratton2016second} proved convergence to a second-order stationary point as well as derived worst-case guarantees (upper bounds only) for a variant of a direct search approach that constructs an approximate Hessian matrix using a finite difference approach\looseness=-1.

\vspace{-2mm}
\paragraph{Model-based approaches.}
Model-based methods construct a function approximation which is used to compute the next step and updated at every iteration. The literature on such methods is broad and we refer the reader to~\cite{larson2019derivative} for a survey of the relevant work. To the best of our knowledge, these methods also require an approximation of the Hessian to obtain second-order guarantees.

\vspace{-2mm}
\paragraph{Lower bounds.} \cite{zabinsky2013stochastic, golovin2019gradientless} show that random search suffers from an exponential dependency to the dimension of the problem. Similar results are discussed in~\cite{duchi2015optimal} for derivative-free methods that use two function evaluations to approximate derivatives. These results serve as a motivation for this work to improve the complexity of random search methods in terms of the input dimension. We mostly focus on designing a new type of random search that achieves better worst-case guarantees than existing lower bounds for vanilla random search. 

\vspace{-2mm}
\paragraph{Inexact Power Iteration.}
Computing the largest eigenvector of a matrix has many applications in statistics and data analysis. Iterative methods such as the power iteration and the Lanczos algorithm are commonly used to solve this problem~\cite{golub2013matrix}. Perhaps the most complete and up-to-date analysis of convergence of inexact power methods for eigenspace estimation can be found in~\cite{balcan2016improved,hardt2014noisy}. Crucially, these convergence rates depend on the eigenvalue distribution. If one instead only seeks a direction aligned with a suboptimal but large eigenvalue, the rate becomes independent of the eigenvalue distribution~\cite{kuczynski1992estimating}. Alternatively to the power method, stochastic methods such as Oja's algorithm~\cite{oja1985stochastic} benefit from a cheaper iteration cost. Improvements have been proposed in~\cite{shamir2015stochastic} that analyzes a variance reduced method that achieves a linear rate of convergence. An accelerated variant of stochastic PCA is also proposed in~\cite{de2018accelerated}.\looseness=-1

%% file: 03_analysis.tex
\vspace{-2mm}
\section{Analysis}
\label{sec:analysis}
\vspace{-2mm}
We work in $\R^d$ with the standard Euclidean norm $\|\cdot\|$. Our goal is to optimize a twice continuously differentiable non-convex function $f(\cdot):\R^d\to \R$ without having access to gradients. We need the following assumption, standard in the literature on SOSPs~\cite{jin2019nonconvex}, also in the DFO setting~\cite{vlatakis2019efficiently}.
\begin{assumption}
The function $f(\cdot)$ is lower bounded, $L_1$-smooth and $L_2$-Hessian Lipschitz.
\label{ass:smoothness}
\end{assumption}
\vspace{-2mm}
In line with many recent works on non-convex optimization~\cite{jin2017escape,jin2018accelerated,Jin2019leadingeigval,daneshmand2018escaping,carmon2018accelerated,vlatakis2019efficiently,arjevani2020second}, our goal is to find an $(\epsilon, \gamma)$-second-order stationary point~(SOSP), i.e. a point $\x$ such that:
\vspace{-1mm}
\begin{equation}
\| \nabla f(\x) \| \leq \epsilon \text{ and } \lambda_{\min}(\nabla^2 f(\x)) \geq -\gamma,
\label{eq:second_order_stat}
\end{equation}
with $\epsilon,\gamma>0$. We analyze two different DFO algorithms: (1) a two-step random search method~(RS, Algorithm~\ref{algo:random_search}) -- similar to STP~\cite{bergou2019stochastic} and (2) a novel random search method~(RSPI, Algorithm~\ref{algo:random_search_PI}) that extracts negative curvature via~Algorithm~\ref{alg:inexact_PI}. We show that, while the approximation of negative eigendirections requires extra computation, the overall complexity of RSPI to find a SOSP is much lower than RS \textit{in terms of number of function evaluations.}

\begin{tcolorbox}
\begin{theorem}[Main result, informal]
The complexity --- in terms of number of function evaluations --- of vanilla Random Search~(similar to STP, Algorithm~\ref{algo:random_search}) for finding second-order stationary points depends \textbf{exponentially} on the problem dimension~(see Lemma~\ref{lemma:probability} and Theorem~\ref{thm:rand_search}).\\This dependency can be reduced to \textbf{linear} by computing an approximate negative curvature direction at each step~(i.e. Algorithm~\ref{algo:random_search_PI}, see Theorem~\ref{thm:conv_random_PI}).
\label{thm:main}
\end{theorem}
\end{tcolorbox}
We will use $\mathcal{O}(\cdot)$, $\Theta(\cdot)$, $\Omega(\cdot)$ to hide constants which do not depend on any problem parameter.

\vspace{-1mm}
\subsection{Two-step Random Search}
\vspace{-2mm}
We analyze a variant of random search that uses a strategy consisting of two steps designed to maximally exploit gradient and curvature. We will use this variant as a surrogate for vanilla random searches such as STP~\cite{bergou2019stochastic} that only use one step. We note that the convergence rate of the latter approach does not theoretically outperform the two-step approach since the gradient step is the same in both approaches. Algorithm~\ref{algo:random_search} shows that the two-step random search samples two symmetric perturbations from a sphere in $\R^d$ with some radius, and updates the current solution approximation $\x_k$ if one of the two perturbations decreases the objective function. The sampling procedure is repeated twice at each iteration, using different sampling radii $\sigma_1$ and $\sigma_2$. One sampling radius ($\sigma_1$) is tuned to exploit the gradient of the objective function while the other ($\sigma_2$) is tuned to exploit negative curvature. One could also interpret the two steps as a one step algorithm with an adaptive step-size, although we explicitly write down the two steps for pedagogic reasons. Note that, in our analysis, the values $\sigma_1$ and $\sigma_2$ are constant, although we will see that in practice, decreasing values can be used, as is typically the case with step-sizes in gradient-based methods (see experiments).
We divide the analysis into two cases: one for the case when the gradient is large and the other when we are close to a strict saddle~(i.e. we still have negative curvature to exploit).

\begin{algorithm}[ht]
\caption{{\sc Two-step Random search (RS)}. {\small Similar to the STP method~\cite{bergou2019stochastic}, but we alternate between two perturbation magnitudes: $\sigma_1$ is set to be optimal for the large gradient case, while $\sigma_2$ optimal to escape saddles.}}
\begin{algorithmic}[1]
\STATE {\bf{Parameters}} $\sigma_1,\sigma_2>0$ (see Theorem~\ref{thm:rand_search})
\STATE Initialize $\x_0$ at random
\FOR{$k=0, 2, 4,\cdots, 2K$}
\STATE $\s_1 \sim \S^{d-1}$ (uniformly)
\STATE \small$\x_{k+1} = \arg \min \{ f(\x_k), f(\x_k + \sigma_1 \s_1), f(\x_k - \sigma_1 \s_1) \} \qquad \qquad \quad$ \textit{\# Gradient step} 
\STATE \normalsize$\s_2 \sim \S^{d-1}$ (uniformly)
\STATE \small$\x_{k+2} = \arg \min \{ f(\x_{k+1}), f(\x_{k+1} + \sigma_2 \s_2), f(\x_{k+1} - \sigma_2 \s_2) \} \qquad$ \textit{\# Curvature step}
\STATE \normalsize Optional: Update $\sigma_1$ and $\sigma_2$ (see experiments)
\ENDFOR
\end{algorithmic}
\label{algo:random_search}
\end{algorithm}

\vspace{-5mm}
\paragraph{Case 1: Large gradients.} First we consider the case $\| \nabla f(\x) \| \geq \epsilon$. Under this assumption, the stepsize can be tuned~\footnote{Such tuning is optimal, please refer to the proof for details.} to yield a decrease proportional to $\epsilon^2/d$.
\begin{tcolorbox}
\begin{restatable}{lemma}{RSLargeGradient}
\label{lemma:large_gradient}
Let $f(\cdot)$ be $L_1$-smooth and $\| \nabla f(\x_k) \| \geq \epsilon$. Algorithm~\ref{algo:random_search} with $\sigma_1 =\epsilon/(L_1 \sqrt{2 \pi d})$ yields $\E [f(\x_{k+1}) - f(\x_{k})|\x_k] \leq - \Omega \left (\frac{\epsilon^2}{L_1 d} \right)$, where $ \E [\cdot|\x_k]$ denotes the conditional expectation w.r.t. $\x_k$.
\end{restatable}
\end{tcolorbox}
\vspace{-2mm}
The proof follows directly from Lemma 3.4 in~\cite{bergou2019stochastic} and is presented in the appendix.

\begin{wrapfigure}[9]{r}{0.35\textwidth}
  \centering
\includegraphics[width=\linewidth]{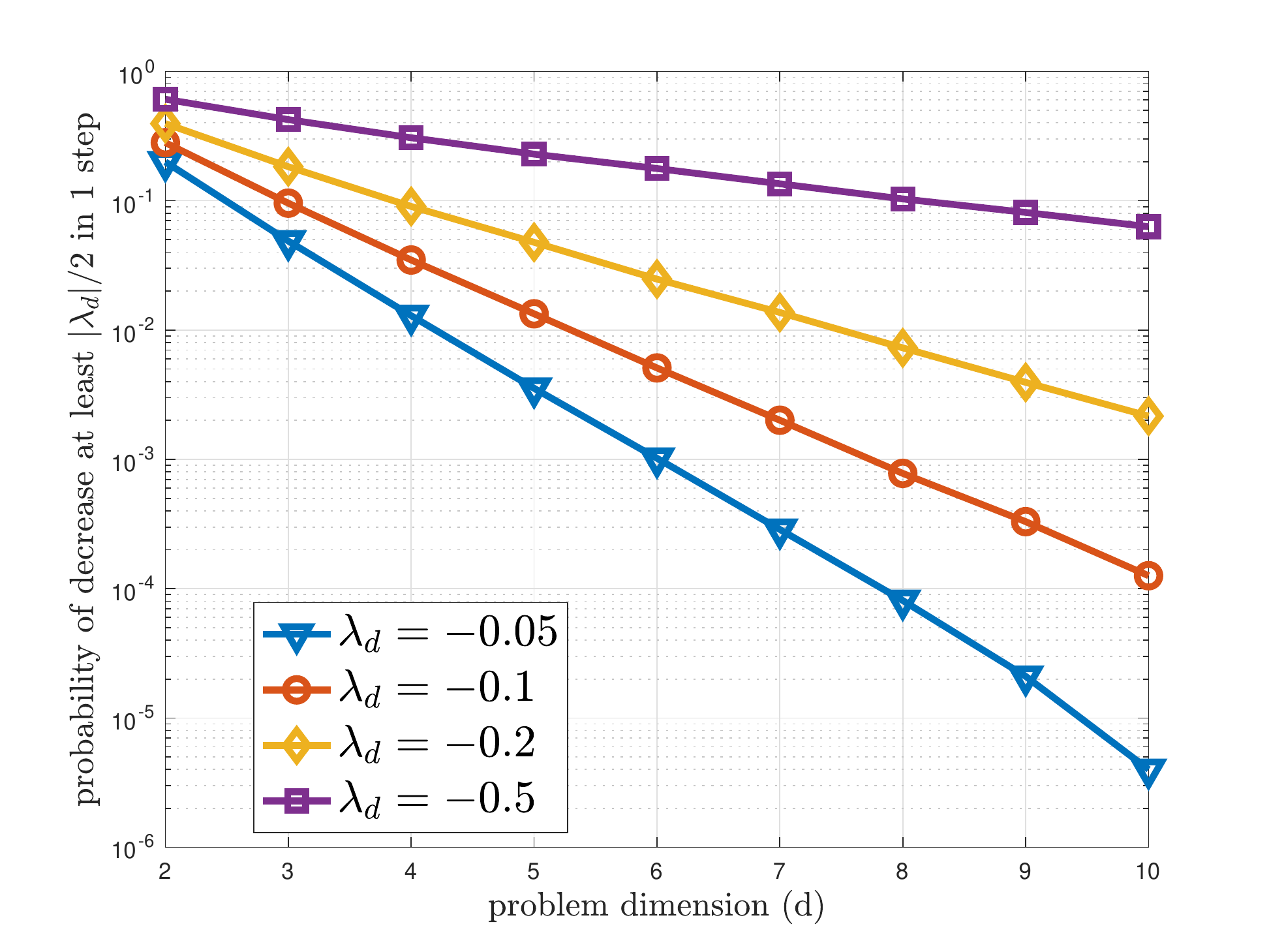}
\vspace{-4mm}
\caption{\small Behavior of a step of vanilla random search on a quadratic saddle centered at the origin. The Hessian has $d-1$ positive eigenvalues equal to $1$ and one negative eigenvalue equal to $\lambda_d$. Plotted is the probability of a decrease~($1e6$ runs) of at least $|\lambda_d|/2$ starting from the origin ($\sigma_2=1$). Performance degrades exponentially with the problem dimension, as predicted by Lemma~\ref{lemma:probability}.}
\label{fig:exponential}
\end{wrapfigure}
\vspace{-2mm}
\paragraph{Case 2: Close to a strict saddle.} We now address the case where $\| \nabla f(\x) \| \leq \epsilon$ but $\lambda_{\min}(\nabla^2 f(\x) )\leq -\gamma$, with $\gamma,\epsilon>0$. Similarly to the analysis of gradient-based methods~\cite{jin2019nonconvex,levy2016power}, our approach consists in first approximating $f(\x)$ around $\x$ with a quadratic $\tilde{f}(\x)$, and then bounding the error on the dynamics. Hence, our first task is to estimate the probability of having a decrease in function value with a single step of random search around a \textit{quadratic saddle}.

\begin{tcolorbox}
\begin{restatable}[Curse of dimensionality of RS around a saddle]{lemma}{RSSquadratic}
\label{lemma:probability}
Consider a $d$-dimensional ($d\ge 4$) quadratic saddle $\tilde f(\cdot)$ centered at the origin with eigenvalues $\lambda_1\ge\dots\ge\lambda_d$, with $\lambda_d<0$. Set $\gamma:=|\lambda_d|$ and $L_1:=\max\{\lambda_1,|\lambda_d|\}$~(cf. definition SOSP in Equation~\ref{eq:second_order_stat}). Starting from the origin, a random step $\s_2 \sim \S^{d-1}(\sigma_2)$ is such that
\vspace{-2mm}
\begin{equation}
    \Pr\left[\tilde f(\s_2) - \tilde f(\mathbf{0})\le-\frac{\gamma}{2}\sigma^2_2\right]\ge\left(\frac{\gamma}{4L_1}\right)^{\frac{d}{2}}.
\end{equation}
Moreover, if $\lambda_1=\lambda_2=\dots =\lambda_{d-1}>0$ (worst case scenario), we also have
\begin{equation}
    \Pr\left[\tilde f(\s_2) - \tilde f(\mathbf{0})\le-\frac{\gamma}{2}\sigma^2_2\right]\le\mathcal{O}\left(2^{-d} \sqrt{d} \right).
\end{equation}
\end{restatable}
\end{tcolorbox}
\vspace{-2mm}
The proof is presented in the appendix and is based on arguments on the geometry of high-dimensional spheres: finding an escape direction at random becomes increasingly hard as the dimension increases.

Lemma~\ref{lemma:probability} provides both an upper and a lower bound on the expected\footnote{Let $V$ be an event that happens in a trial with probability $p$. Then the expected number of trials to first occurrence of $V$ is $1/p$.} number of function evaluations needed by vanilla random search to escape a saddle. As it is the case for vanilla gradient descent~\cite{du2017gradient}, the complexity grows exponentially with the dimension.
We emphasize that the exponential dependency itself is not a new result (see Section~\ref{sec:related_work}) and that it can be extended to hold globally (not only for one step)~\cite{zabinsky2013stochastic}. We will use the upper bound of Lemma~\ref{lemma:probability} to derive a new result in Theorem~\ref{thm:rand_search} that proves convergence to a \emph{second-order} stationary point.

Further, we note that, for the bounds in Lemma~\ref{lemma:probability}, the decrease in function value is proportional to $\sigma_2^2$ --- yet the decrease probability is independent from $\sigma_2^2$. This might seem unintuitive at first, but it is simply due to the radially isotropic structure of quadratic landscapes. However, the landscape around $\x_k$ is only approximately quadratic : under Assumption~\ref{ass:smoothness}, we have~(for a proof please refer to~\cite{levy2016power}) that $|f(\y) - \tilde{f}(\x)| \leq \frac{L_2}{2} \| \y - \x \|^3$ for all $\x,\y \in \R^d$. Hence, in order to lift the bound to the non-quadratic case we need to consider ``small'' perturbations $\sigma_2=\mathcal{O}(1/L_2)$, so that the quadratic approximation remains valid. The proof of the next result can be found in the appendix.
\begin{tcolorbox}
\begin{restatable}{lemma}{RSSaddle}
\label{lemma:saddle_case}
Let $f(\cdot)$ be $L_1$-smooth and $L_2$-Hessian-Lipschitz. Assume $\| \nabla f(\x_k) \| \leq \epsilon$ and $\lambda_{\min}(\nabla^2 f(\x_k) )\leq -\gamma = - \epsilon^{2/3}$.
Then Algorithm~\ref{algo:random_search} with $\sigma_2 = \frac{\epsilon^{2/3}}{2 L_2}$ is s.t.
\vspace{-2mm}
\begin{equation}
\E[f(\x_{k+1}) - f(\x_{k})|\x_k] \leq - \Omega\left(\left(\frac{\gamma}{4L_1}\right)^{\frac{d}{2}}\epsilon^2\right).
\label{eq:exp_decrease}
\end{equation}
\end{restatable}
\end{tcolorbox}
With no surprise, we see that the exponential complexity already found for the quadratic case in Lemma~\ref{lemma:probability} directly influences the magnitude of the per-iteration decrease. We remark that the choice $\gamma = \epsilon^{2/3}$ is optimal for our proof, based on Hessian smoothness~(see proof of Lemma~\ref{lemma:saddle_case} \& \ref{lemma:saddle_case_PI}).
\vspace{-2mm}
\paragraph{Joint analysis.} We define the following sets that capture the different scenarios:
\begin{align*}
\A_1 = \left \{ \x : \| \nabla f(\x) \| \geq \epsilon \right \};\
\A_2 = \left \{ \x : \| \nabla f(\x) \| < \epsilon \text{ and } \lambda_{\min}(\nabla^2 f(\x)) \leq -\gamma \right \}.
\vspace{-2mm}
\end{align*}
By Lemma~\ref{lemma:large_gradient} and ~\ref{lemma:saddle_case}, we have that for any $\x_k \in \A_1 \cup \A_2$, Equation~\ref{eq:exp_decrease} holds --- the bottleneck scenario being $\mathcal{A}_2$. Since under Assumption~\ref{ass:smoothness} the function is lower bounded and Algorithm~\ref{algo:random_search} is a descent method, we get convergence in a finite number of steps to an SOSP. 
\begin{tcolorbox}
\begin{theorem}
\label{thm:rand_search}
Fix the desired level of accuracy $\epsilon>0$. Assume $f(\cdot)$ satisfies Assumption~\ref{ass:smoothness}. Two-step RS~(Alg.~\ref{algo:random_search}) with parameters $\sigma_1 = \frac{\epsilon}{L_1 \sqrt{2\pi d}}$ and $\sigma_2 = \frac{\epsilon^{2/3}}{2 L_2}$ returns in expectation an $(\epsilon, \epsilon^{2/3})$-second-order stationary point in $\bigO\left(\kappa^{d} \epsilon^{-2}\right)$ function evaluations, where $\kappa=\Theta(L_1/\gamma)= \Theta(L_1\epsilon^{-2/3})$.
\end{theorem}
\end{tcolorbox}

\vspace{-2mm}
\subsection{Power Iteration Random Search}
\label{sec:power-iter-random-search}
\vspace{-2mm}

\begin{wrapfigure}[14]{r}{0.35\textwidth}
\vspace{-5mm}
  \centering
\includegraphics[width=\linewidth]{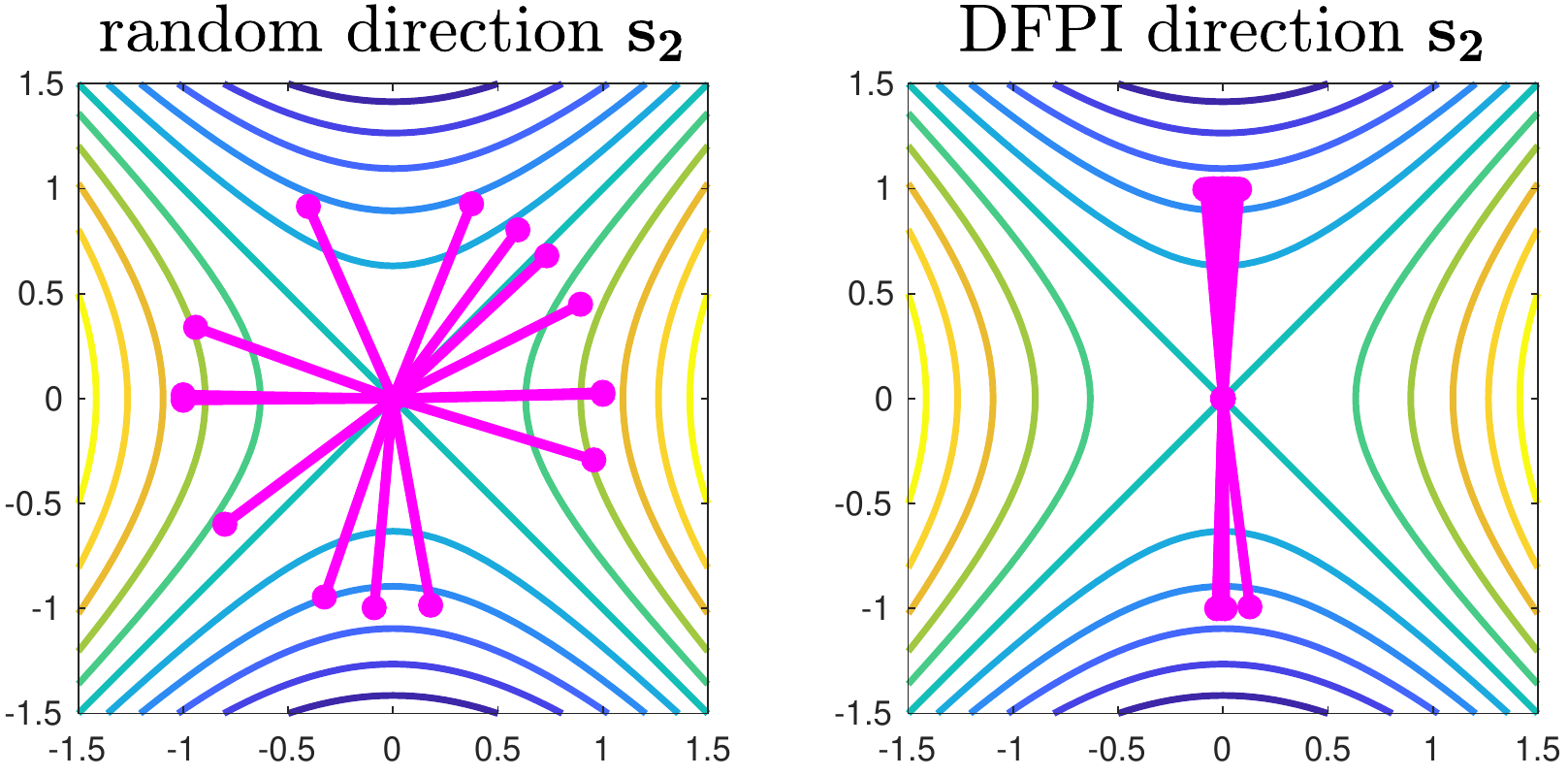}
\vspace{-5mm}
\caption{\small Difference between the two procedures to generate $\s_2$ when initialized at a quadratic saddle (blue $=$ small values). Each magenta point is a random update direction $\s_2\sim\S^{d-1}$ (left) or the result of the DFPI procedure with 3 iterations (right).}
\label{fig:DFPI}
\end{wrapfigure}
The purpose of this section is to modify vanilla random search to reduce the exponential dependency of the convergence rate on the problem dimension, illustrated by the lower bound in Lemma~\ref{lemma:probability}. We propose a new type of random search that is coupled with a derivative-free routine (named DFPI) in order to overcome the curse of dimensionality. We note that this type of ``hybrid" method is not completely new in the literature, see e.g. ~\cite{gratton2016second}, but in contrast to prior work, the computational cost of the DFPI routine in terms of the input dimension is low. This is achieved by computing an approximation of the eigenvector corresponding to the most negative eigenvalue of the Hessian based on a noisy power method~\cite{hardt2014noisy}, which is inspired from recent works on Hessian-free negative curvature exploitation techniques in gradient-based optimization~\cite{liu2017noisy,carmon2018accelerated}. The resulting method is shown as Algorithm~\ref{algo:random_search_PI} and the DFPI routine is presented as Algorithm~\ref{alg:inexact_PI}. We highlight that RSPI does not require computation and storage of the full Hessian~(in contrast to AHDS that performs~$d^2$ function calls). We instead approximate only the leading eigenvector using approximate zero-order Hessian-vector products~($\approx d\log d$ function calls:~Lemma 9). This provides a large speed-up in high dimensions, see wall-clock times reported later.

\begin{algorithm}[ht]
\caption{{\sc Random search PI (RSPI)}.
{\small Same structure as Alg.~\ref{algo:random_search}, only difference is in line 6: the curvature exploitation step is based on a perturbation sampled from a non-isotropic distribution aligned with negative eigendirections, computed by Alg.~\ref{alg:inexact_PI}~(see Fig.\ref{fig:DFPI}).}}
\begin{algorithmic}[1]
\STATE {\bf{Parameters}} $\sigma_1,\sigma_2>0$ (see Thm~\ref{thm:conv_random_PI})
\STATE Initialize $\x_0$ at random
\FOR{$k=0, 2, 4, \cdots 2K$}
\STATE $\s_1 \sim \S^{d-1}$ (uniformly)
\STATE \small $\x_{k+1} = \arg \min \{ f(\x_k), f(\x_k + \sigma_1 \s_1), f(\x_k - \sigma_1 \s_1) \}$
\STATE \normalsize \textcolor{bulgarianrose}{$\s_2 = \text{DFPI}(\x_k) \leftarrow$ Algorithm 3}
\STATE \small$\x_{k+2} = \arg \min \{ f(\x_{k+1}), f(\x_{k+1} + \sigma_2 \s_2), f(\x_{k+1} - \sigma_2 \s_2) \}$
\STATE Optional: Update $\sigma_1$ and $\sigma_2$ (see experiments)
\ENDFOR
\end{algorithmic}
\label{algo:random_search_PI}
\end{algorithm}

\vspace{-1mm}


\begin{algorithm}[ht] 
\begin{algorithmic}[1]
\STATE {{\bf{Parameters}} $c,r,\eta>0 \text{ and } T_\DFPI \in \mathbb{Z}^+$}~(see Theorem\ref{thm:conv_random_PI})
\STATE {\bf{INPUTS} : $\x \in \R^d$, }
\STATE $\s_2^{(0)} \sim \S^{d-1}$ (uniformly) \\
\FOR{$t = 0 \dots T_\DFPI-1$}
\STATE Set $\g_+ = \sum\limits_{i=1}^d \frac{f(\x + r\cdot\s^{(t)}_2 + c \cdot  \e_i) - f(\x + r\cdot\s^{(t)}_2 - c \cdot  \e_i)}{2c} \e_i$ \\
\STATE Set $\g_- = \sum\limits_{i=1}^d \frac{f(\x - r\cdot\s^{(t)}_2 + c \cdot \e_i) - f(\x - r\cdot\s^{(t)}_2 - c \cdot  \e_i)}{2c} \e_i$ \\
\STATE \textbf{Update}: $\s^{(t+1)}_2 = \s^{(t)}_2 - \eta \frac{\g_+ - \g_-}{2r}$ \\
\STATE Normalize $\s^{(t+1)}_2 = \s^{(t+1)}_2/\|\s^{(t+1)}_2\|$
\ENDFOR
\STATE {\bf RETURN : $\s^{(T_{\textmd{DFPI}})}_2$ }
\end{algorithmic}
\caption{{\sc{\small Derivative-Free Power Iteration (DFPI)}}
{\small A noisy derivative-free power method to approximate negative curvature directions, as shown in Lemma~\ref{prop:DFPI_error}. Every iteration of DFPI requires $4d$ function evals, hence the routine requires $K_\DFPI = 4d T_\DFPI$ function evals overall.}}
\label{alg:inexact_PI}
\end{algorithm}

As in the last subsection, we split the analysis into two cases. First, in the case of large gradients, i.e. $\| \nabla f(\x) \| \geq \epsilon$, we can simply re-use the result of Lemma~\ref{lemma:large_gradient}. The second case (sufficient negative curvature) requires to show that the vector $\s_2$ returned by the DFPI procedure (Alg.~\ref{alg:inexact_PI}) yields a decrease of the objective.
\vspace{-2mm}
\paragraph{Exploiting negative curvature close to a saddle.}
We consider again the case where $\| \nabla {f}(\x_k) \| \leq \epsilon$ and $\lambda_{\min}(\nabla^2 {f}(\x_k) )\leq -\gamma$. We saw in Lemma~\ref{lemma:probability}, that isotropic sampling in the neighborhood of $\x_k$ provides a function decrease only after an exponential number of iterations, in the worst case. We show that if the perturbation is instead sampled from a distribution which satisfies the next assumption, the number of required iterations drastically decreases. 

\begin{assumption}
Consider $\x$ s.t. $\lambda_{\min}(\nabla^2 {f}(\x) )\leq -\gamma$. The direction $\s_2$, output of DFPI after $K_\DFPI$ function evaluations, returns in expectation a good approximation to the most negative eigenvalue of $\nabla^2 f(\x)$. Specifically, $\E[\s_2^\top \nabla^2 f(\x) \s_2]\le \lambda_{\min}(\nabla^2 {f}(\x) )+\frac{\gamma}{2}$.
\label{ass:error_eigenvector}
\end{assumption}

This assumption --- \textit{which we will soon formally validate} --- is motivated by the fact that DFPI is an approximate (noisy) power method on $\Am(\x) :=\Im-\eta\nabla^2 f(\x)$, hence can be used to estimate the maximum eigenvalue of $\Am(\x)$ --- which is minimum eigenvalue of $\nabla^2 f(\x)$ if $\eta\le 1/L_1$.

\begin{tcolorbox}
\begin{restatable}{lemma}{ERRORDFPI}
\label{prop:DFPI_error}
Let $f(\cdot)$ be $L_1$-smooth and $L_2$-Hessian-Lipschitz. The iteration of DFPI can be seen as a step of a noisy power method: $\s_2^{(t+1)} = (\Im-\eta\nabla^2 f(\x))\s_2^{(t)} + \vxi_\DFPI^{(t)}$, with $\|\vxi_\DFPI^{(t)}\| = \mathcal{O}(r L_2 + \frac{c}{r} L_1\sqrt{d})$. In particular, $\|\vxi_\DFPI^{(t)}\|\to 0$ as $r,\tfrac{c}{r}\to 0$; hence the error can be made as small as needed within the limits of numerical stability. In addition, if $f(\cdot)$ is quadratic, we have $\|\vxi_\DFPI^{(t)}\|=0$.
\end{restatable}
\end{tcolorbox}
\vspace{-2mm}
The proof is presented in the appendix. The bound on $\vxi_\DFPI^{(t)}$ is enough for us to apply the well-known convergence rates for the noisy power method~\cite{hardt2014noisy} and to motivate a $\mathcal{O}(\log(d))$ bound on the DFPI iterations needed to satisfy Assumption~\ref{ass:error_eigenvector}. Before diving into this, we show that directly using Assumption~\ref{ass:error_eigenvector} actually makes the iteration\footnote{The dependency on the dimension will however show up in the final number of function evaluations~(Theorem~\ref{thm:conv_random_PI}). However, such dependency is not exponential even in the worst case.} complexity dimension independent.

\begin{tcolorbox}
\begin{restatable}{lemma}{RSPISaddle}
\label{lemma:saddle_case_PI}
Let $f(\cdot)$ be $L_1$-smooth and $L_2$-Hessian-Lipschitz, and assume $\| \nabla f(\x_k) \| \leq \epsilon$ and $\lambda_{\min}(\nabla^2 f(\x_k) )\leq -\gamma = - \epsilon^{2/3}$. Under Assumption~\ref{ass:error_eigenvector},~RSPI~(Algorithm~\ref{algo:random_search_PI}) with $\sigma_2 = \frac{\gamma}{2 L_2}$~(choice as Theorem~\ref{thm:rand_search}) yields $\E[f(\x_{k+1}) - f(\x_{k})|\x_k]\leq - \Omega(\epsilon^2)$, independent of the problem dimension.
\end{restatable}
\end{tcolorbox}

\vspace{-2mm}
\paragraph{Iteration complexity.} We now combine Lemma~\ref{lemma:large_gradient} and~\ref{lemma:saddle_case_PI}.
\begin{tcolorbox}
\begin{proposition}
Fix the desired level of accuracy $\epsilon>0$. Assume $f(\cdot)$ satisfies Assumptions~\ref{ass:smoothness} and \ref{ass:error_eigenvector}. RSPI~(i.e. Alg.~\ref{algo:random_search_PI}) with parameters $\sigma_1 = \frac{\epsilon}{L_1 \sqrt{2\pi d}}$ and $\sigma_2 = \frac{\epsilon^{2/3}}{2 L_2}$~(same choice as~Theorem~\ref{thm:rand_search}) returns in expectation an $(\epsilon, \epsilon^{2/3})$-second-order stationary point in $\bigO\left(\epsilon^{-2}\right)$ iterations.
\label{prop:conv_random_PI_it}
\end{proposition}
\end{tcolorbox}
\vspace{-2mm}
\paragraph{Overall number of function evaluations.}
While Proposition~\ref{prop:conv_random_PI_it} shows that the number of RSPI iterations is, conditioned on Assumption~\ref{ass:error_eigenvector}, independent of the problem dimension, it hides the number of function evaluations needed for the assumption to hold. To include this into the final complexity~(Theorem~\ref{thm:conv_random_PI}) -- \textit{i.e. to drop Assumption~2} -- we need a bound for convergence of noisy power methods~\cite{hardt2014noisy}.

\begin{tcolorbox}
\begin{restatable}[Consequence of Corollary 1.1 in~\cite{hardt2014noisy}] {lemma}{complexity_power}
\label{lemma:complexity_power} Let the parameters of DFPI be such that the error $\vxi_\DFPI$ is small enough~(always possible, as shown in Lemma~\ref{prop:DFPI_error}). Let $\eta\le 1/L_1$. Let $\gamma=\epsilon^{2/3}$; for a fixed RSPI iteration, $T_\DFPI =\mathcal{O}\left(\epsilon^{-2/3}L_1\log\left(\frac{d}{\delta^2}\right)\right)$ DFPI iterations are enough to ensure validity of Assumption~\ref{ass:error_eigenvector} at $\x_k$ with probability $1-\delta^{\Omega(1)}-e^{\Omega(d)}$.
\label{lemma:Hardt}
\end{restatable}
\end{tcolorbox}
A derivation is included for completeness in the appendix. We are ready to state the main result.
\begin{tcolorbox}
\begin{theorem}
Fix the accuracy $\epsilon>0$. Assume $f(\cdot)$ satisfies Assumption~\ref{ass:smoothness}. Algorithm~\ref{algo:random_search_PI} with parameters $\sigma_1 = \frac{\epsilon}{L_1 \sqrt{2\pi d}}$ and $\sigma_2 = \frac{\epsilon^{2/3}}{2 L_2}$~(same choice as~Theorem~\ref{thm:rand_search}) combined with DFPI~(Alg.~\ref{alg:inexact_PI}) with parameters $\eta\le1/L_1$, $T_\DFPI=\mathcal{O}\left( \epsilon^{-2/3}L_1\log(d)\right)$ and $c,r$ sufficiently small, returns in expectation an $(\epsilon, \epsilon^{2/3})$-second-order stationary point in $\bigO\left(\epsilon^{-8/3} d\log(L_1 d)\right)$ function evaluations.
\label{thm:conv_random_PI}
\end{theorem}
\end{tcolorbox}
\vspace{-2mm}
\begin{proof}
Fix $\delta$ s.t. the probability of success $p := 1-\delta^{\Omega(1)} -e^{-\Omega(d)} = \mathcal{O}(1)$ in Lemma~\ref{lemma:Hardt} is positive. The number of function evaluations needed for Lemma~\ref{lemma:Hardt} to hold in expectation at each step~(i.e. the burden of Assumption~\ref{ass:error_eigenvector}) is then $K_\DFPI/p = \mathcal{O}\left( \epsilon^{-2/3}L_1 d\log(d)\right)$. We conclude by Proposition~\ref{prop:conv_random_PI_it}.
\end{proof}

\vspace{-2mm}
\paragraph{Significance and Novelty of the analysis.} To the best of our knowledge, the result in Theorem~\ref{thm:conv_random_PI} is the first to prove convergence of a type of random search to a \emph{second-order} stationary point with a \emph{linear} dependency to the input dimension. Our experimental results (section~\ref{sec:experiments}) confirm this significant speed-up is observed in practice. The analysis presented in appendix relies on some new geometric arguments for high-dimensional spaces coupled with more classical optimization bounds. We again emphasize that the analysis technique is different from gradient-free methods that approximate the gradients. For instance~\cite{vlatakis2019efficiently} define the error between the approximate gradient $q(x_k, h_k)$ and the exact gradient $\nabla f(x_k)$ as $\epsilon_k = q(x_k, h_k) - \nabla f(x_k)$. They then require this error to be bounded, i.e. $|\epsilon_k| \leq c_h |h|$ where $c_h$ is a constant and $h$ is the size of the step used in the finite difference formula. The constant $c_h$ is directly controlled by $h$ and can be made as small as possible to better approximate the result of the gradient-based method. In our case, the update direction is sampled at random, and we therefore have to rely on a different proof technique that involves  probabilistic geometric arguments in high-dimensional spaces.

\vspace{-2mm}
\paragraph{How to further speed up DFPI with SPSA.} Each DFPI iteration requires $4d$ function evaluations. While this complexity is necessary to build an arbitrarily good finite difference (FD) approximation of second order information~(needed by Lemma~\ref{lemma:Hardt}), in practice a more rough estimate of Hessian-vector products can be obtained using cheap randomized techniques. In particular, in the experiments in the next section, we show that an SPSA estimator~\cite{spall1992multivariate} of $\g_+$ and $\g_-$ is in fact sufficient to achieve acceleration in the performance with respect to the two-step random search. In particular, SPSA computes $\g_+,\g_-$ as $\g_{\pm} = \sum\limits_{i=1}^d \frac{f(\x \pm r\s^{(t)}_2 + c \vDelta) - f(\x \pm r\s^{(t)}_2 - c \vDelta)}{2c \Delta_i} \e_i,$
where $r,c>0$ and $\vDelta$ is a vector of $d$ random variables, often picked to be symmetrically Bernoulli distributed. SPSA is asymptotically unbiased and only requires $4$ function calls, as opposed to the $4d$ needed by FD. However, the variance of SPSA does not arbitrarily decrease to zero as $c,r$ vanish, as it is instead the case for FD: it saturates~\cite{spall1992multivariate}. As a result, Lemma~\ref{lemma:Hardt} would not always hold. In the appendix, we provide an extensive comparison between FD and SPSA for DFPI: we show that, for the sake of Assumption~\ref{ass:error_eigenvector}, the error in the SPSA estimator is  acceptable for our purposes, for small enough $\eta$~(see use in Alg.~\ref{alg:inexact_PI}).

%% file: 04_experiments.tex
\section{Experiments}
\vspace{-2mm}
\label{sec:experiments}

In this section we verify our theoretical findings empirically. Specifically, we set two objectives: (1) to evaluate the performance RS and  RSPI, to verify the validity of Theorem~\ref{thm:main}, and (2) to compare the performance of these two methods against existing random search methods. For the latter, we run the two-step Random Search (RS) and Random Search Power Iteration (RSPI) against the Stochastic Three Points (STP) method, the Basic Direct Search (BDS)~\cite{vicente2013worst} and the Approximate Hessian Direct Search (AHDS). We recall that AHDS explicitly constructs an approximation to the Hessian matrix in order to extract negative curvature. Descriptions of each algorithm are provided in the~appendix.

\vspace{-2mm}
\paragraph{Setup.}All experiments follow a similar procedure. In each task, all algorithms are initialized at a strict saddle point and executed for the same number of iterations. We report the optimality gap as a function of iterations and wall-clock time. Additionally, we report the norm of the gradient vector in the appendix. Since all algorithms are stochastic, the experimental process is repeated multiple times (wherever possible using a different saddle point as initialization) and the results are averaged across all runs. For each task, the hyperparameters of every method are selected based on a coarse grid search refined by trial and error. For RS and RSPI the parameters $\sigma_1$ and $\sigma_2$ are initialized and updated in the same manner, hence the only difference between the two is that RSPI extracts negative curvature explicitly whereas the two-step RS samples a direction at random (see Figure~\ref{fig:DFPI}). We choose to run DFPI for $20$ iterations for all the results shown in the paper. Empirically, we observed that performing more iterations does not further improve the overall performance of RSPI. The hyperparameters used for each method are provided in the appendix and the code for reproducing the experiments is available online\footnote{\url{https://github.com/adamsolomou/second-order-random-search}}. 

\begin{figure*}[!t]
    \includegraphics[width=\textwidth]{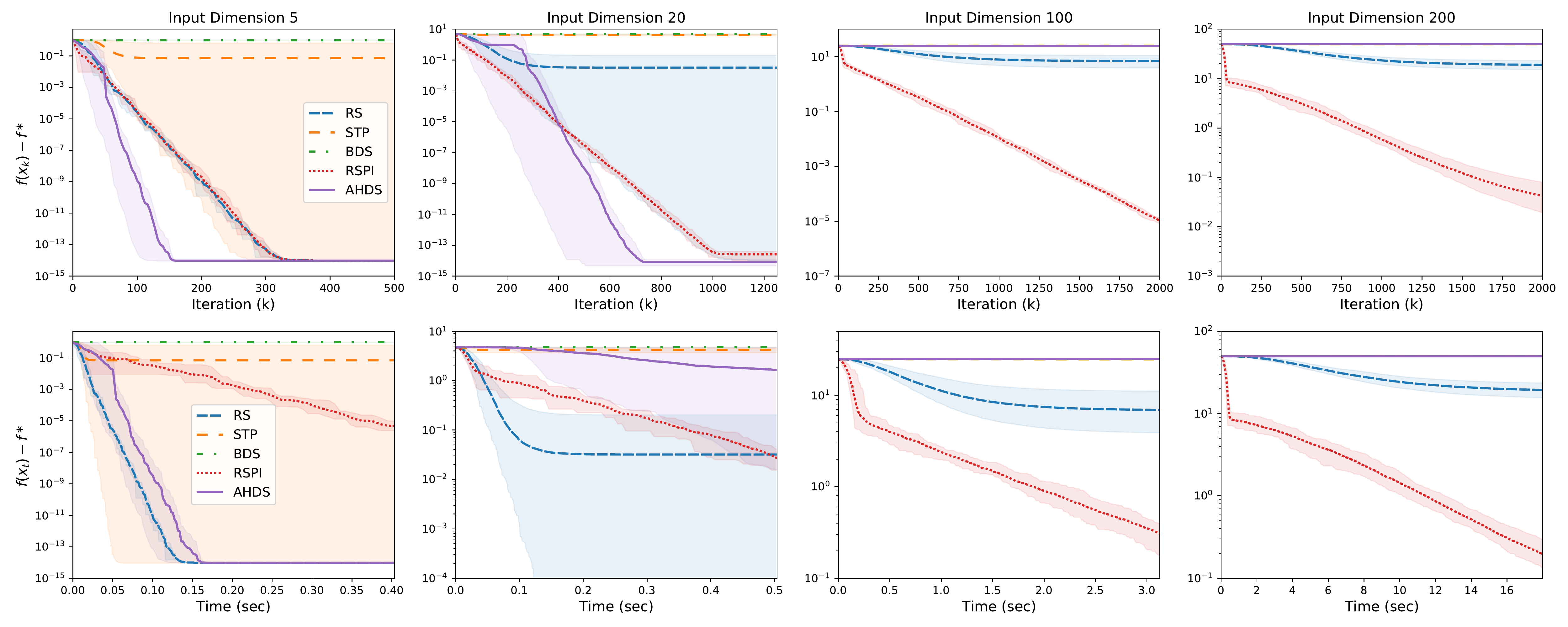}
    \vspace{-7mm}
    \caption{\small{Performance while minimizing the objective in Eq. (\ref{eq:sb-objective-def}) for different $d$. Confidence intervals show min-max intervals over ten runs. All algorithms are initialized at the strict saddle point across all runs. For $d=100,200$, the lines for STP, BDS and AHDS overlap each other as none of the methods achieve progress in terms of function value.} }
    \label{fig:sb-performance}
    \vspace{-2mm}
\end{figure*}

\vspace{-2mm}
\paragraph{Function with growing dimension.}
We start by considering the following benchmarking function (see App.~\ref{app:exp} for an illustration)
\vspace{-5mm}
\begin{equation}\label{eq:sb-objective-def}
    f(x_1, \cdots, x_d, y) = \frac{1}{4}\sum_{i=1}^d x_i^4 - y\sum_{i=1}^d x_i + \frac{d}{2}y^2,  
\end{equation}
which has a unique strict saddle point (i.e. with negative curvature) at $P1 = (0,\cdots,0,0)$ and two global minima at $P2 = (1,\cdots,1,1)$ and $P3 = (-1,\cdots,-1,-1)$. The results in Fig.~\ref{fig:sb-performance} illustrate that both the two-step RS method as well as the RSPI algorithm are able to consistently escape the saddle across all dimensions. While in low-dimensional settings ($d=5$) the RSPI algorithm is outperformed by the two-step RS and the AHDS in terms of their behavior as a function of run-time, the situation is clearly reversed as the dimensionality grows. For~$d=100,200$ the two-step RS achieves progress at a very slow rate and a high number of iterations is needed in order to achieve convergence to a second-order stationary point. In contrast, RSPI successfully approximates the negative curvature of the objective to efficiently escape the saddle point, allowing the algorithm to achieve progress at a faster rate. The fact that for $d=5$ the AHDS algorithm requires less time than RSPI to converge to a second-order stationary point, indicates that in low-dimensional settings the cost incurred by the power iterations within RSPI is higher than the cost of approximating the entire Hessian matrix. However, for higher values of $d$, AHDS quickly becomes inefficient and expensive. Further, STP performs worse than RS, simply because it employs only one sampling radius~(RS uses two).

\begin{figure*}[!b]
    \centering
    \includegraphics[width=\textwidth]{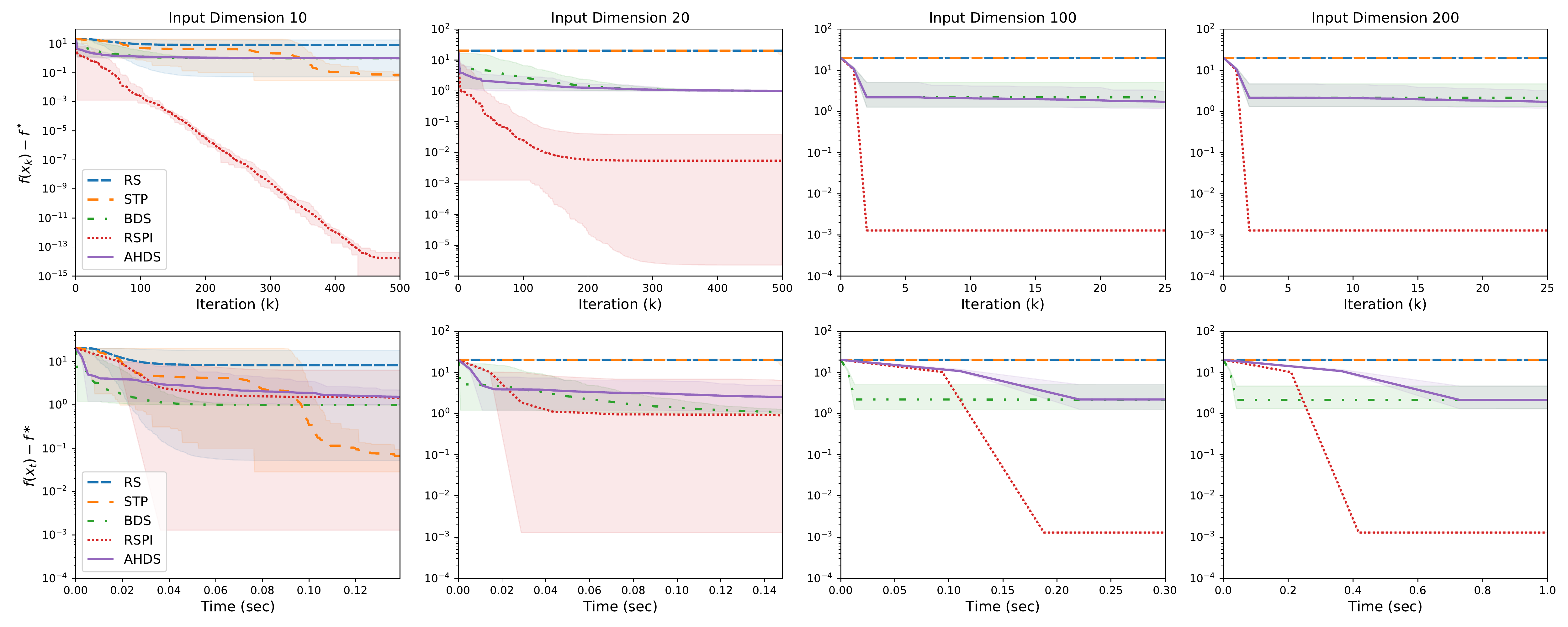}
    \vspace{-6mm}
    \caption{\small{Optimality gap on the Rastrigin function as a function of iterations (top) and running time (bottom). Confidence intervals show min-max intervals over ten runs. All algorithms are initialized at a strict saddle point and executed for a total of $500$ iterations across all runs, however for $d=100,200$ no further improvement is achieved after $25$ iterations.}}
    \label{fig:rastrigin-performance}
    \vspace{-2mm}
\end{figure*}

\paragraph{Rastrigin function.} Next, we conduct experiments on the Rastrigin function, a popular benchmark in the literature~\cite{hansen2009rastrigin}. For any $\x \in \R^d$, the $d$-dimensional Rastrigin function is defined as 
\vspace{-3mm}
\begin{equation}
    f(\x) = 10d + \sum_{i=1}^d (x_i^2 - 10\cos(2\pi x_i)).
\end{equation}
The function has a unique global minimizer at $\x^* = \mathbf{0}$, whereas the number of stationary points (including strict saddles) grows exponentially with $d$. Based on Lemma~\ref{lemma:probability}, we expect that having a single direction of negative curvature will challenge the core mechanism of each algorithm while trying to escape the saddle. To that end, we ensure that at each initialization point there exist a single direction of negative curvature across all settings of~$d$. More details about the initialization process and the implications on the results are given in the appendix. \newline
The results in Figure~\ref{fig:rastrigin-performance} illustrate the strong empirical performance of RSPI, not only in comparison to the two-step RS algorithm but also against the rest of the algorithms. RSPI successfully approximates the single direction of negative curvature and escapes the saddle point after one iteration. On the contrary, the two-step RS achieves minimal progress even for low dimensional settings, whereas for $d=20$ it requires more than $300$ iterations to escape the saddle (see the gradient norm plot in the appendix). For higher dimensional settings, two-step RS does not escape at all, supporting our theoretical argument that as the problem dimension grows the probability of sampling a direction that is aligned with the direction of negative curvature decreases exponentially. Lastly, both BDS and AHDS consistently escape the saddle point across all values of $d$. However, their performance remains suboptimal compared to RSPI.

\begin{figure}
    \centering
    \vspace{-2mm}
    \includegraphics[width=0.8\linewidth]{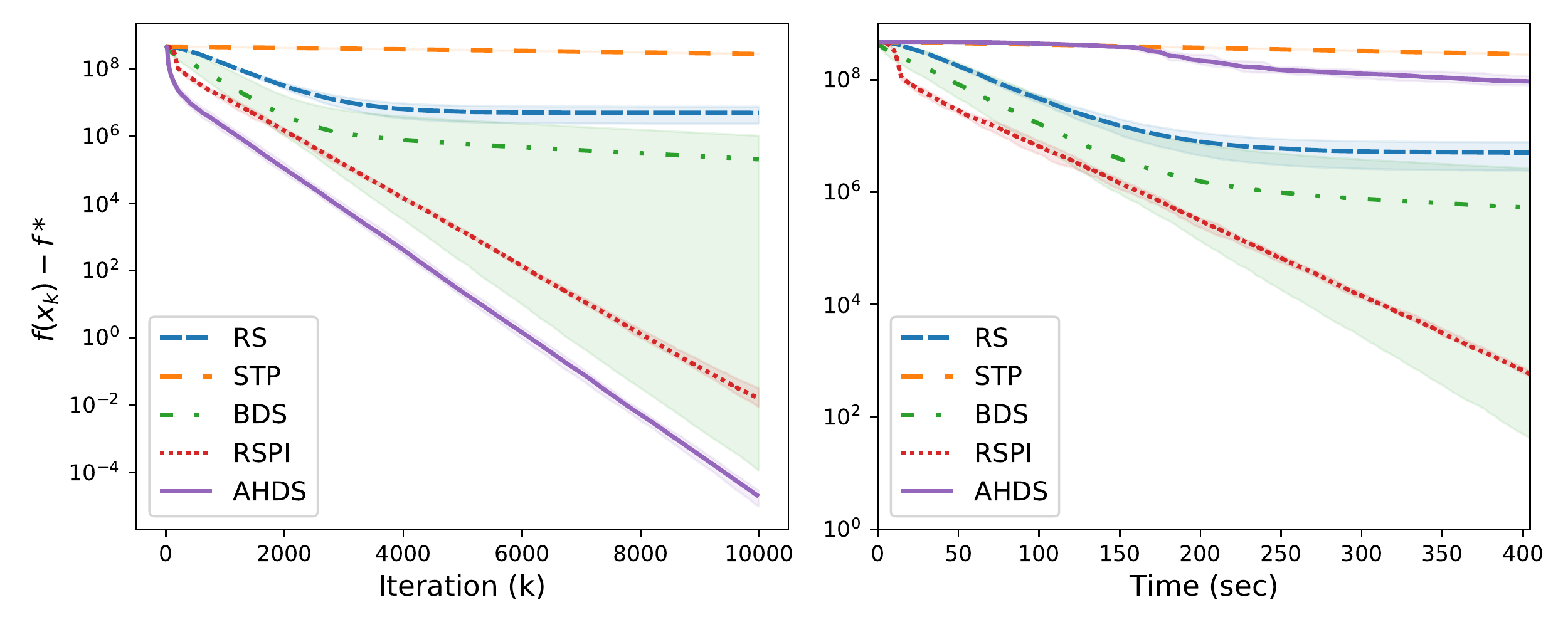}
    \vspace{-4mm}
    \caption{\small{Empirical performance in finding the leading eigenvector of a $350$-dimensional random matrix. Confidence intervals show min-max intervals over five runs. All algorithms are initialized at a strict saddle point.}}
    \label{fig:le-performance}
\end{figure}

\paragraph{Leading eigenvector problem.} Finally, we consider the task of finding the leading eigenvector of a positive semidefinite matrix $\Mm \in \R^{d\times d}$. The problem is equivalent to minimizing $f(\x) = \|\x \x^T - \Mm \|_F^2$~\cite{Jin2019leadingeigval}. Figure~\ref{fig:le-performance} shows the empirical performance in finding the leading eigenvector of a $350$-dimensional random matrix. While at an iteration level AHDS appears to be very effective, when the iteration complexity of each method is taken into account it is the slowest to escape. Notably, a single iteration of RSPI takes (on average) $0.06$ seconds, whereas a single iteration of AHDS takes approximately $9.01$ seconds. This experiment clearly illustrates the computational advantages that RSPI provides while provably ensuring convergence to second order stationary points.  

%% file: 05_conclusion.tex

\vspace{-2mm}
\section{Conclusion}
\label{sec:conclusion}
\vspace{-1mm}
We analyzed the convergence of two types of random search methods. The first approach is a slight variation of standard random search that converges to a second-order stationary point but whose worst-case analysis demonstrates exponential complexity in terms of the function input dimension. The second random search approach we propose extracts negative curvature using function evaluations. Importantly, the dependency in terms of the function input dimension becomes linear, a result which we clearly observed in our experimental results, especially in terms of run-time.\\
There are a number of avenues to pursue as future work. To start off, (1)~a simple extension would be to allow DFPI to store multiple candidates for negative eigendirections. As discussed in~\cite{musco2015randomized} and formally shown in~\cite{gu2015subspace}, this can directly boost performance. Similarly~(2) one could study the application of a zero-th order version of Neon~\cite{xu2018first}. (3)~It would be then interesting to also understand if injecting noise in the process~(additional exploration) can help in escaping saddles~\cite{du2017gradient}. (4)~Further, instead of using constant values or a predefined schedule for $\sigma_1$ and $\sigma_2$, one could analyze the commonly used adaptive strategy where these values are adapted according to whether the function is being decreased~\cite{vicente2013worst}. (5) We note that one could in principle relax Assumption~\ref{ass:smoothness} and instead work with a smoothed version of $f(\cdot)$, without requiring differentiability. (6) Finally, it would be interesting to benchmark DPFI on other machine learning problems, including for instance reinforcement learning tasks where random search methods are becoming more prevalent~\cite{mania2018simple, maheswaranathan2019guided}. Another potentially interesting direction would be to extend our analysis to random search methods for min-max optimization problems~\cite{anagnostidis2021direct}.

%% file: appendix.tex

\newpage
\onecolumn
\appendix


{\huge \begin{center}
    Appendix
\end{center}}

\section{Analysis of random search~(Algorithm~\ref{algo:random_search})}
\label{app:random_search}
We start by studying some properties of high dimensional spheres. We then apply these properties to show how the rate of the two-step random search~(Algorithm~\ref{algo:random_search}) depends exponentially on the problem dimension.

\subsection{High-dimensional spheres and curse of dimensionality}
We denote by $V_R(d)$ and $A_R(d)$ the volume and the surface area of the $(d-1)$ sphere with radius $R$: $\mathcal{S}^{d-1}(R)=\{\x\in\R^d \ | \ \|\x\|=R\}$. It is well known \cite{hopcroft2012computer} that the following formulas hold:
\begin{equation}
    A_R(d) = \frac{2\pi^{d/2}}{\Gamma\left(\frac{d}{2}\right)}R^{d-1},\quad V_R(d) = \frac{A(d) R}{d} = \frac{2\pi^{d/2}}{d \ \Gamma\left(\frac{d}{2}\right)}R^d.
    \label{eq:area_volume}
\end{equation}

Moreover, we have the following important lemma, which can also be found in Section~1.2.4 of~\cite{hopcroft2012computer}.
\begin{tcolorbox}
\begin{lemma}\label{lemma:initial_volume}
 Let $\varsigma\ge 0$ and define $A^\varsigma_1(d)$ to be the surface area of the cap $\{\x\in\R^d \ | \ \|\x\|=1, x_1\ge\varsigma\}$, $d\ge 2$. We have:
 \begin{equation}
     A^\varsigma_1(d) = A_1(d-1)\int_{\varsigma}^1(1-x_1^2)^{\frac{d-2}{2}}dx_1.
 \end{equation}
\end{lemma}
\end{tcolorbox}
\begin{proof}
The radius of the spherical cap at height $x_1$ is $\sqrt{1-x_1^2}$, and we have that $A_{\sqrt{1-x_1^2}}(d-1) = A_1(d-1) \left(\sqrt{1-x_1^2}\right)^{d-2}$ by the surface area formula in Equation~\ref{eq:area_volume}. Since $A^\varsigma_1(d) = \int_{\varsigma}^1 A_{\sqrt{1-x_1^2}}(d-1) dx_1$, we conclude.
\end{proof}

We will need both an upper and a lower bound on the integral above. The next result shows that both these bounds are exponential.
\begin{tcolorbox}
\begin{lemma}\label{lemma:integral_bounds} For any $\alpha>1$,
$$\left[\frac{1}{2}(1-\varsigma^2)\right]^{\alpha+1}\le\int_{\varsigma}^1(1-x^2)^{\alpha}dx\le (1-\varsigma^2)^\alpha.$$
\label{lemma:bound_integral}
\end{lemma}
\end{tcolorbox}
\begin{proof}
The upper bound is straightforward. The lower bound in an application of H\"older's inequality~(see e.g. Equation 1.1 in~\cite{cheung2001generalizations}):  for a real number $p>1$ and functions $f$ and $g$ regular enough,
\begin{equation}
    \int_\varsigma^1 |f(x)g(x)|dx\le\left[\int_\varsigma^1 |f(x)|^p dx\right]^{1/p}\left[\int_\varsigma^1 |g(x)|^{\frac{p}{p-1}} dx\right]^{\frac{p-1}{p}}.
\end{equation}
Take $g$ to be constant equal to one. Then, taking everything to power $p$
\begin{equation}
    \left[\int_\varsigma^1 |f(x)|dx\right]^p\le (1-\varsigma)^{p-1} \int_\varsigma^1 |f(x)|^p dx .
\end{equation}
By applying this formula and after performing a few algebraic manipulations, we get
\begin{multline}
 \int_{\varsigma}^1(1-x^2)^{\alpha}dx\ge\frac{\left(\int_{\varsigma}^1(1-x^2)dx\right)^{\alpha}}{(1-\varsigma)^{\alpha-1}} = (1-\varsigma)\left(\frac{\varsigma^3-3\varsigma+2}{3(1-\varsigma)}\right)^\alpha\\ = (1-\varsigma)\left(\frac{(1-\varsigma)^2 (\varsigma+2)}{3(1-\varsigma)}\right)^\alpha = (1-\varsigma)\left(\frac{1}{3}(1-\varsigma)(\varsigma+2)\right)^\alpha\ge \left[\frac{1}{2}(1-\varsigma^2)\right]^{\alpha+1},
\end{multline}
where in the last inequality we used the fact that for $\varsigma\in[0,1]$, $\frac{1}{3}(1-\varsigma)(\varsigma+2)\ge \frac{1}{2}(1-\varsigma^2)$.
\end{proof}

\begin{figure}[ht!]
    \centering
    \includegraphics[width=\textwidth]{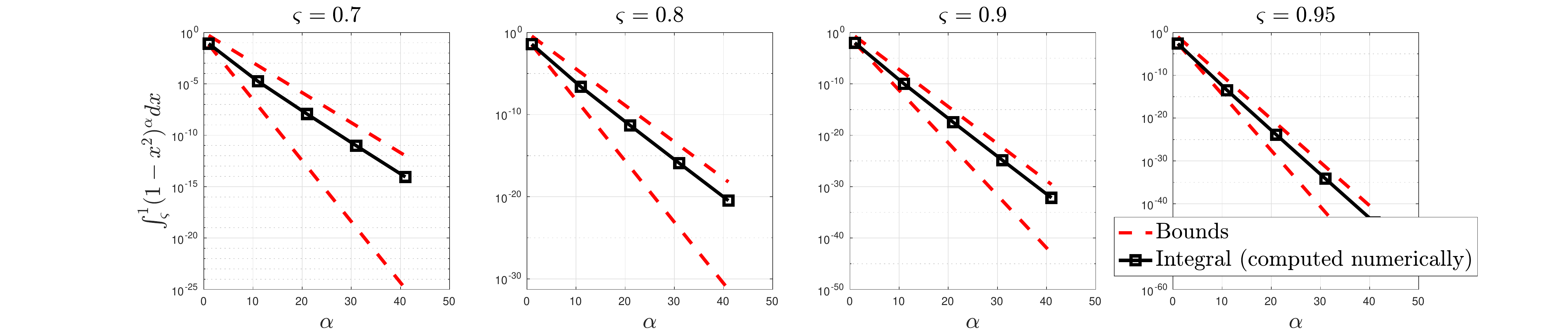}
    \caption{Numerical verification of Lemma~\ref{lemma:bound_integral}. Integral computed numerically using the MATLAB \texttt{integral} function.}
    \label{fig:bound_integral}
\end{figure}

A verification of the bound above can be found in Figure~\ref{fig:bound_integral}. We note that the upper bound becomes tight as $\alpha\to\infty$, and that the lower bound becomes less pessimistic as $\varsigma\to1$.

\begin{figure}[ht!]
    \centering
    \includegraphics[width=\textwidth]{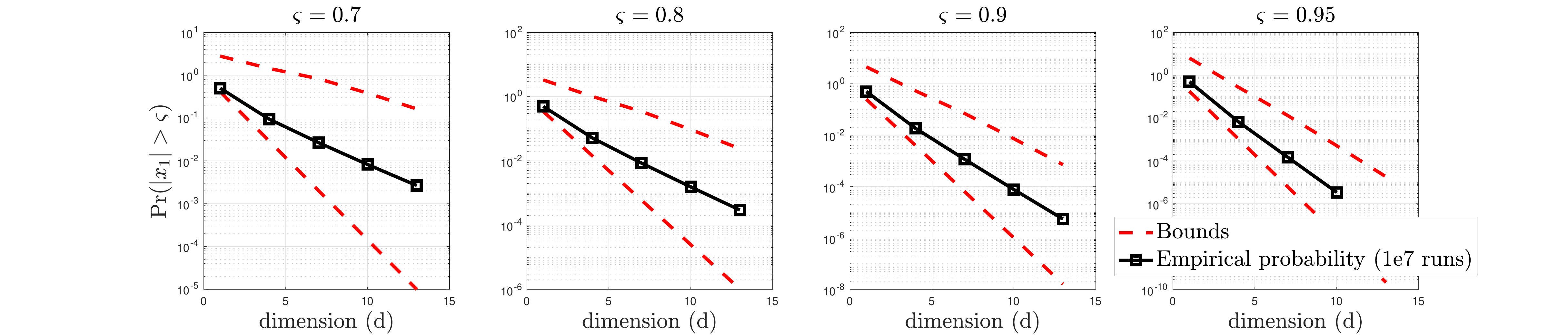}
    \caption{Numerical verification of Lemma~\ref{lemma:spheres}. Bounds can be found in Lemma~\ref{lemma:spheres}.}
    \label{fig:bound_spheres}
\end{figure}

Putting the previous two lemmas together, we get
\begin{equation}
   A_1(d-1) \left[\frac{1}{2}(1-\varsigma)\right]^{\alpha+1} \ \ \le \ \  A^\varsigma_1(d) \ \ \le \ \ A_1(d-1) (1-\varsigma^2)^\alpha,
\end{equation}
where $\alpha = (d-2)/2$. Now we are ready to state the final lemma for high-dimensional spheres, which is verified empirically in Figure~\ref{fig:bound_spheres}.

\begin{tcolorbox}
\begin{lemma}[Curse of dimensionality]
Let $\x$ be a random point on the surface of the unit $d$-ball in Euclidean space, with $d\ge 4$. For $\varsigma\in(0,1)$, we have
\begin{equation}
     \left[\frac{1}{2}(1-\varsigma^2)\right]^{\frac{d}{2}}\le \Pr[|x_1|>\varsigma]\le 2\sqrt{d-2}\left[1-\varsigma^2\right]^{\frac{d}{2}-1} .
\end{equation}
\label{lemma:spheres}
\end{lemma}
In particular, the probability of being $\varsigma$-away from the equator \textbf{decays exponentially} with the number of dimensions.
\end{tcolorbox}
\begin{proof}
The proof is just a matter of finding good upper and lower bounds on $A(d)$ as a function of $A(d-1)$, to combine with the result of Lemma~\ref{lemma:integral_bounds}. We are going to use the lower bound on the surface area by \cite{hopcroft2012computer}~(Equation 1.3): $A_1(d)\ge\frac{1}{\sqrt{d-2}}A_1(d-1)$. For an easy upper bound, we can instead pick $A_1(d)\le 2 A_1(d-1)$ (surface of the enclosing cylinder). Combining Lemma~\ref{lemma:initial_volume} with Lemma~\ref{lemma:integral_bounds} and the bounds we just found, we get
\begin{align}
    &\Pr[|x_1|\ge\varsigma] = \frac{A_1^\varsigma(d)}{\frac{1}{2}A_1(d)}\le\frac{(1-\varsigma^2)^{\frac{d-2}{2}}A_1(d-1)}{\frac{1}{2\sqrt{d-2}}A_1(d-1)},\\
    &\Pr[|x_1|\ge\varsigma] = \frac{A_1^\varsigma(d)}{\frac{1}{2}A_1(d)}\ge\frac{\left[\frac{1}{2}(1-\varsigma^2)\right]^{\frac{d}{2}}A_1(d-1)}{A_1(d-1)}.
\end{align}
\end{proof}

\subsection{Non-convex dynamics --- the quadratic case}
Here we seek to understand the behaviour of random search around a point $\x_k$ with negative curvature, we consider the quadratic approximation $\tilde{f}(\x) = f(\x_k) + \nabla f(\x_k)^\top(\x-\x_k) + \frac{1}{2}(\x-\x_k)^\top \nabla^2 f(\x_k) (\x-\x_k)$ where $\nabla^2 f(\x_k) \in \R^{d \times d}$. By the spectral theorem, we have $\nabla^2 f(\x_k) = \Vm^\top \Lambda \Vm$, where $\Vm = [\v_i]_{i=1}^d$, $\v_i \in \R^d$ contains an orthonormal basis of eigenvectors of $\nabla^2 f(\x_k)$ and $\Lambda$ is a diagonal matrix containing the eigenvalues $\lambda_1\ge\lambda_2\ge\cdots\ge\lambda_d$ of $\nabla^2 f(\x_k)$ (counted together with their multiplicity). For the setting considered in this paragraph, we have $\lambda_d<0$.

In our first result, we consider the case $f(\x_k)=0$ and $\nabla f(\x_k) = \textbf{0}$.

\begin{tcolorbox}
\RSSquadratic*
\end{tcolorbox}

\begin{proof}
Let $x_0$ be any initial point. We seek the probability of the event
\begin{equation}
    E_{\text{decr}}:= \{\s_2^T \nabla^2 \tilde{f}(\x_0) \s_2\le -\zeta\},
    \label{eq:e}
\end{equation}
for some positive $\zeta$. First, we divide everything by $\|\s_2\|^2 = \sigma^2_2$, to effectively reduce the problem to the special case $\sigma_2^2=1$. We get
\begin{equation}
    E_{\text{decr}}= \left\{\left(\frac{\s_2}{\|\s_2\|}\right)^T \nabla^2 \tilde{f}(\x_0) \frac{\s_2}{\|\s_2\|}\le -\tilde\zeta\right\},
\end{equation}
where $\tilde\zeta = \zeta/\sigma^2_2$. Let us now write $\s_2/\|\s_2\|$ in the eigenbasis $\{\v_i\}_{i=1}^d$ of the Hessian $\nabla^2 \tilde{f}(\x_0)$. We have that 
\begin{equation}
    \frac{\s_2}{\|\s_2\|} = \sum_{i=1}^d a_i\v_i, \quad \sum_{i=1}^d a_i^2 = 1.
\end{equation}
Hence, we can write
\begin{equation}
    E_{\text{decr}}= \{\lambda_1 a_1^2 +\lambda_2 a_2^2+\dots + \lambda_d a_d^2\le -\tilde\zeta\}.
\end{equation}
To bound the probability of this event, we construct the smaller event $E^*_{\text{decr}}\subseteq E_{\text{decr}}$:
\begin{equation}
    E^*_{\text{decr}}:= \{\lambda_1 a_1^2 +\lambda_1 a_2^2+\dots +\lambda_1 a_{d-1}^2 \le |\lambda_d| a_d^2 -\tilde\zeta\}.
    \label{eq:e_star}
\end{equation}

This event can be written in a reduced form, using the fact that $\sum_{i=1}^d a_i^2 = 1$; indeed

\begin{align}
    &\lambda_1 a_1^2 +\lambda_1 a_2^2+\dots +\lambda_1 a_{d-1}^2 \le |\lambda_d| a_d^2 -\tilde\zeta\\
    \iff & \lambda_1 a_1^2 +\lambda_1 a_2^2+\dots +\lambda_1a_{d-1}^2 + \lambda_1 a_d^2 \le (|\lambda_d| + \lambda_1) a_d^2 -\tilde\zeta\\
    \iff & \lambda_1\le(|\lambda_d| + \lambda_1) a_d^2 -\tilde\zeta\\
    \iff & a_d^2\ge\frac{\lambda_1+\tilde\zeta}{\lambda_1+|\lambda_d|}.
\end{align}

In conclusion, we find
\begin{equation}
    \Pr[E^*_{\text{decr}}] = \Pr\left[|a_d| \ge \varsigma \right],\quad \varsigma:= \sqrt{\frac{\lambda_1+\tilde\zeta}{\lambda_1+|\lambda_d|}}.
\end{equation}

\begin{figure}[ht!]
    \centering
    \includegraphics[width=0.35\textwidth]{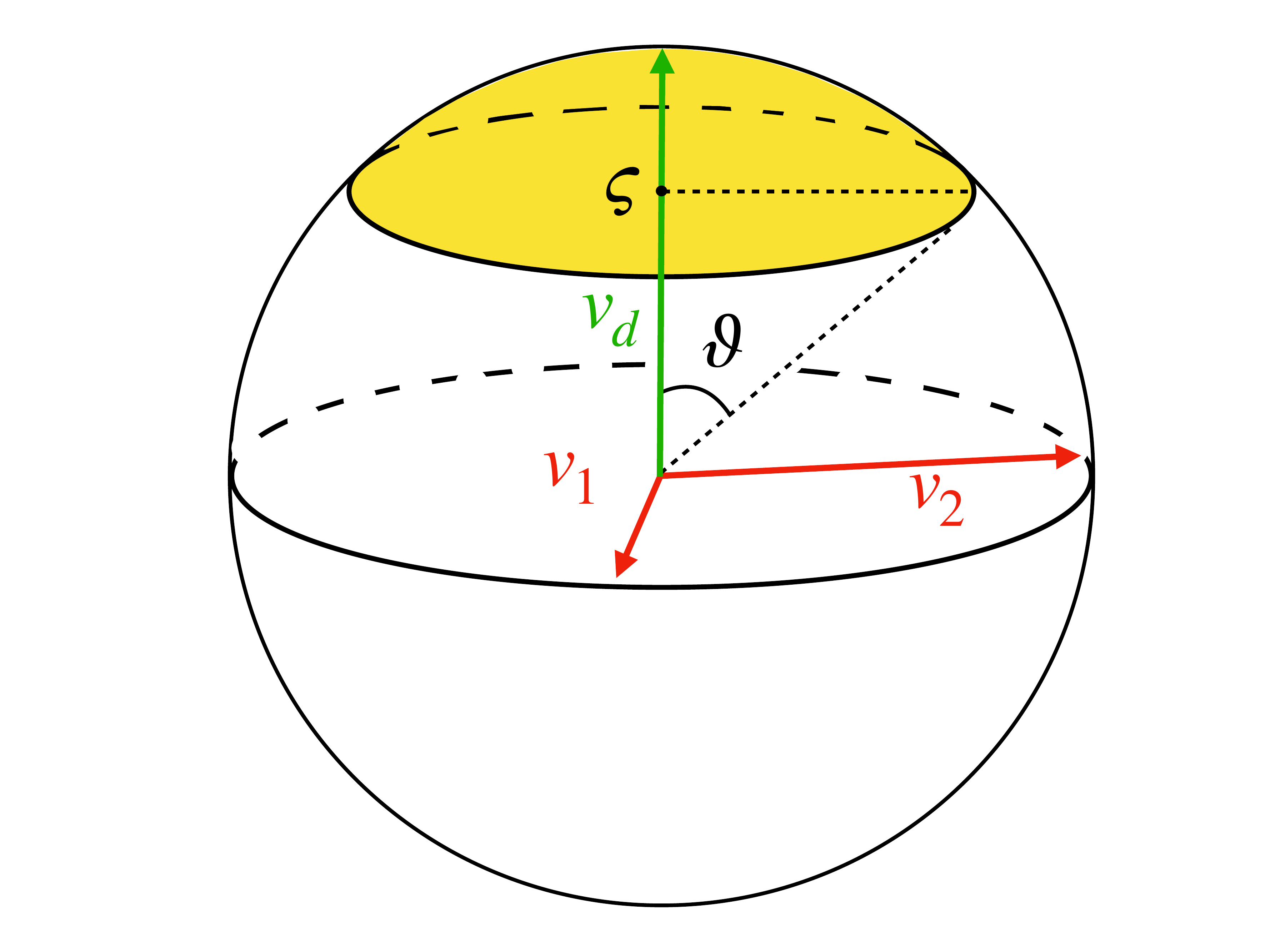}
    \caption{Illustration for the proof of Lemma~\ref{lemma:probability}. Any vector on the unit sphere whose angle that is less than $\vartheta = \cos^{-1}(\varsigma)$ away from $\v_d$ belongs to cap colored in yellow. Our goal is to bound the surface area of this spherical cap.}
    \label{fig:eig_vectors}
\end{figure}

Therefore, since $\mathbf a = (a_1,a_2,\dots, a_d)$ is uniformly distributed on the surface of the unit sphere in $\R^d$, we have reduced the problem to finding the surface of a spherical cap~(see Figure~\ref{fig:eig_vectors}). From \eqref{lemma:spheres}, we directly get
\begin{equation}
     \left[\frac{1}{2}(1-\varsigma^2)\right]^{\frac{d}{2}}\le \Pr[E^*_{\text{decr}}] \le 2\sqrt{d-2}\left[1-\varsigma^2\right]^{\frac{d}{2}-1} .
\end{equation}
Plugging in $\zeta = \frac{1}{2}|\lambda_d|\sigma_2^2$, we get $1-\zeta^2 = \frac{|\lambda_d|}{2(\lambda_1 +|\lambda_d|)}$, so by setting $\gamma:=|\lambda_d|$ and $L_1:=\max\{\lambda_1,|\lambda_d|\}$~(cf. definition SOSP in Equation~\ref{eq:second_order_stat}):

\begin{equation}
    \frac{\gamma}{4 L_1}\le \frac{1}{2}(1-\zeta^2) \le  \frac{1}{4},
\end{equation}
this completes the proof.
\end{proof}

\subsection{Analysis for general function}

\begin{tcolorbox}
\RSLargeGradient*
\end{tcolorbox}
\begin{proof}
One can show (see e.g. Lemma 3.4 in~\cite{bergou2019stochastic}) that $\E_{\s_1 \sim\S^{d-1}}[\nabla f(\x_k)^\top \s_1|\x_k] = \frac{1}{\sqrt{\mu_d}} \| \nabla f(\x_k) \|$,
with $\mu_d := 2\pi d$. Using smoothness, we obtain
\begin{align}
&\E[f(\x_{k+1}) - f(\x_{k})|\x_k]\\&\leq \E[\nabla f(\x_{k})^\top (\x_{k+1} - \x_k)|\x_k]+ \frac{L_1}{2} \E[\| \x_{k+1} - \x_k \|^2]\\
&\leq - \frac{\sigma_1}{\sqrt{\mu_d}} \| \nabla f(\x_k) \| + \frac{L_1}{2} \sigma_1^2 \\
&\leq - \frac{\sigma_1}{\sqrt{\mu_d}} \epsilon + \frac{L_1}{2} \sigma_1^2,
\end{align}
where in the first inequality we used the fact that we can choose between $\s_1$ and $-\s_1$, and update with the perturbation which yields the best (i.e. the negative) step. Plugging-in our choice for $\sigma_1$ (which optimizes the quadratic upper bound above) we get the result.
\end{proof}

\begin{tcolorbox}
\RSSaddle*
\end{tcolorbox}

\begin{proof} Since $f(\x)$ is $L_2$-Lipschitz Hessian, we have (see e.g.~\cite{levy2016power})
\begin{align}
&f(\x_{k+1}) - f(\x_{k})\\
&\leq (\x_{k+1} - \x_k)^\top \nabla f(\x_k) + \frac12 (\x_{k+1} - \x_k)^\top \nabla^2 f(\x_k) (\x_{k+1} - \x_k) + \frac{L_2}{6} \| \x_{k+1} - \x_k \|^3.
\end{align}

We use Lemma~\ref{lemma:probability} on the quadratic $\tilde f(\x) : =\frac12 (\x - \x_k)^\top \nabla^2 f(\x_k) (\x - \x_k)$ to guarantee a decrease of $\gamma \sigma_2^2/2$ with probability $p_{\text{decr}} = \left(\frac{\gamma}{4L_1}\right)^{d/2}$. Therefore, with probability $p_{\text{decr}}$,
\begin{align}
&f(\x_{k+1}) - f(\x_{k})\\ &\leq (\x_{k+1} - \x_k)^\top \nabla f(\x_k) + \frac12 (\x_{k+1} - \x_k)^\top \nabla^2 f(\x_k) (\x_{k+1} - \x_k) + \frac{L_2}{6} \| \x_{k+1} - \x_k \|^3 \\
&\leq -\gamma \sigma_2^2/2 + \frac{L_2}{6} \sigma_2^3, \label{eq:saddle_decrease_f}
\end{align}
where we can ensure that $(\x_{k+1} - \x_k)^\top \nabla f(\x_k) \leq 0$ by testing for both $\s_2$ and $-\s_2$ in Algorithm~\ref{algo:random_search_PI} --- which does not affect $\frac12 \s_2^\top \nabla^2 f(\x_k) \s_2$.

\begin{figure}[ht!]
    \centering
    \includegraphics[width=0.28\textwidth]{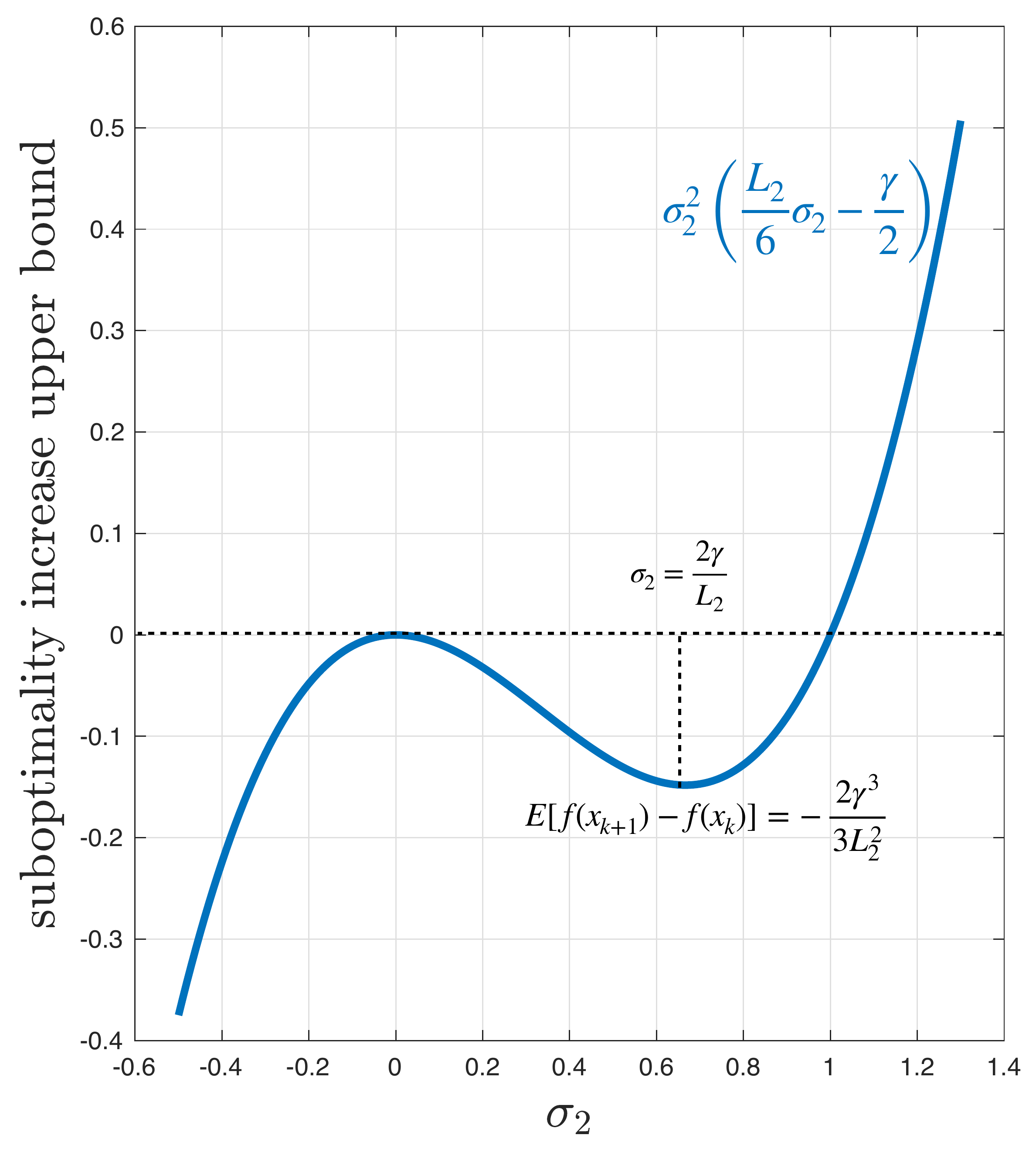}
    \caption{Selection of the value of $\sigma_2$ which yields the best decrease.}
    \label{fig:bound_lemma}
\end{figure}

Next, we seek to minimize Equation~\ref{eq:saddle_decrease_f} with respect to $\sigma_2$. To this, we take the derivative and set it to zero: $\sigma_2=0$ is a local maximizer, while $\sigma_2 = 2\gamma/L_2$ is the unique local minimizer for $\sigma_2\ge0$. Hence, since Equation~\ref{eq:saddle_decrease_f} goes to infinity for $\sigma_2\to \infty$, this minimizer is global (see Figure~\ref{fig:bound_lemma}). For this value of $\sigma_2$, we have
\begin{equation}
    f(\x_{k+1}) - f(\x_{k})\le -\frac{2}{3}\frac{\gamma^3}{L_2^2}.
\end{equation}
Therefore, for $\gamma = \epsilon^{2/3}$, we have $f(\x_{k+1}) - f(\x_{k})\le\Omega(\epsilon^2)$ for $\s_2\in E^*_{\text{decr}}$, defined in Equation~\ref{eq:e_star}. We proceed by computing the expected decrease using the law of total expectation
\begin{align}
    &\E[f(\x_{k+1}) - f(\x_{k})]\\ &=  \E[f(\x_{k+1}) - f(\x_{k})|E^*_{\text{decr}}]\cdot \Pr [E^*_{\text{decr}}] + \E[f(\x_{k+1}) - f(\x_{k})|(E^*_{\text{decr}})^c]\cdot \Pr [(E^*_{\text{decr}})^c]\\
    &\le \E[f(\x_{k+1}) - f(\x_{k})|E^*_{\text{decr}}]\cdot \Pr [E^*_{\text{decr}}]\\
    &= - p_{\text{decr}} \cdot \Omega(\epsilon^2).
\end{align}
where in the first inequality we used the fact that, by the algorithm definition, $f(\x_{k+1}) - f(\x_{k}) = 0$ (rejected step).
\end{proof}


\section{Analysis Random Search PI~(Algorithm~\ref{algo:random_search_PI})}

\begin{tcolorbox}
\RSPISaddle*
\end{tcolorbox}
\begin{proof}Since $f(\x)$ is $L_2$-Lipschitz Hessian, under Assumption~\ref{ass:error_eigenvector} we have
\begin{align}
&\E[f(\x_{k+1}) - f(\x_{k})|\x_k]\\ &\leq \E[(\x_{k+1} - \x_k)^\top \nabla f(\x_k)|\x_k] + \frac12 \E\left[(\x_{k+1} - \x_k)^\top \nabla^2 f(\x_k) (\x_{k+1} - \x_k)|\x_k\right] \\&+ \frac{L_2}{6} \| \x_{k+1} - \x_k \|^3 \\
&\leq -\gamma \sigma_2^2/2 + \gamma \sigma_2^2/4 + \frac{L_2}{6} \sigma_2^3  \\
&= \sigma_2^2 \left( -\frac{1}{4} \gamma + \frac{L_2}{6} \sigma_2 \right),
\label{eq:saddle_decrease_f2}
\end{align}
where we can ensure that $(\x_{k+1} - \x_k)^\top \nabla f(\x_k) \leq 0$ by testing for both $\s_2$ and $-\s_2$ in Algorithm~\ref{algo:random_search_PI} (which does not affect $\frac12 \s_2^\top \nabla^2 f(\x_k) \s_2 + \frac{L_2}{6} \| \s_2 \|^3$).

We therefore require $\sigma_2 \leq \frac{3}{2 L_2} \gamma$ for the RHS in Eq.~\eqref{eq:saddle_decrease_f} to be negative. Choosing, as for the random search case, $\sigma_2 = \frac{\gamma}{2 L_2}$,
\begin{align}
\E[f(\x_{k+1}) - f(\x_{k})|\x_k] \leq -\frac{1}{24} \frac{\gamma^3}{L_2^2}.
\end{align}

For $\gamma = \epsilon^{2/3}$, we obtain $\E[f(\x_{k+1}) - f(\x_{k})|\x_k] \leq -\Omega(\epsilon^2)$.
\end{proof}

\section{Analysis of DFPI~(Algorithm~\ref{alg:inexact_PI})}

\subsection{Proof of Lemma~\ref{prop:DFPI_error}}
We show that line 7 in Algorithm~\ref{alg:inexact_PI} can be written as a noisy power iteration step. That is,
$$\s_2^{(t+1)} = \s_2^{(t)} -\eta\frac{\g_+-\g_-}{2r} \overset{\text{to show}}{=} (\Im-\eta\nabla^2 f(\x))\s_2^{(t)} + \vxi_\DFPI^{(t)},$$
where $\vxi_\DFPI^{(t)}$ is an approximation error. We show that $\vxi_\DFPI^{(t)}$ can be made as small as needed if finite difference hyperparameters $r,c$ are chosen small enough. Therefore, Alg.~\ref{alg:inexact_PI} can be seen as a noisy power method; hence one can motivate the rate in Lemma~\ref{lemma:complexity_power} using the results in~\cite{kuczynski1992estimating,hardt2014noisy,balcan2016improved}, with an argument similar to~\cite{carmon2018accelerated}~(remark after the Assumption 1 of this reference).
\begin{tcolorbox}
\ERRORDFPI*
\end{tcolorbox}
\begin{proof}
We note that $\g_+$ and $\g_-$ are the finite-difference approximations of $\nabla f(\x+r\s_2^{(t)})$ and $\nabla f(\x-r\s_2^{(t)})$, respectively:
\begin{align}
    &\g_+ = \sum\limits_{i=1}^d \frac{f(\x + r\s^{(t)}_2 + c \e_i) - f(\x + r\s^{(t)}_2 - c  \e_i)}{2c} \e_i, \\
    &\g_- = \sum\limits_{i=1}^d \frac{f(\x - r\s^{(t)}_2 + c  \e_i) - f(\x - r\s^{(t)}_2 - c   \e_i)}{2c} \e_i.
    \label{eq:fd_grads}
\end{align}

where $r,c>0$. Let us denote by $\vxi_{1,+}^{(t)}$ and $\vxi_{1,-}^{(t)}$ the approximation errors in the estimation of $\nabla f(\x+r\s_2^{(t)})$ and $\nabla f(\x-r\s_2^{(t)})$, respectively (properties of this error discussed at the end of the proof). We have:
\vspace{-2mm}
\begin{align}
    \frac{\g_+-\g_-}{2r}&=\frac{\nabla f(\x+r\s_2^{(t)}) +\vxi_{1,+}^{(t)} - \nabla f(\x-r\s_2^{(t)})-\vxi_{1,-}^{(t)}}{2r}\\
    &=\nabla^2 f(x)\s_2^{(t)} + \vxi^{(t)}_2 + \frac{\vxi_{1,+}^{(t)}-\vxi_{1,-}^{(t)}}{2r}\\
    &=\nabla^2 f(x)\s_2^{(t)} + \vxi^{(t)}_\DFPI, 
    \label{eq:noisy_expr}
\end{align}
where $\vxi_2^{(t)}$ is the error on the Hessian-vector product. To conclude the proof, we bound the errors $\vxi_2^{(t)}$ and $\vxi_{1,\pm}^{(t)}.$
\vspace{-2mm}
\paragraph{Bound on $\vxi_2^{(t)}$.} This error vanishes as $r\to0$ under Assumption~\ref{ass:smoothness}~(see main paper):
\begin{align}
    \|\vxi_2^{(t)}\|&=\left\|\frac{\nabla f(\x+r\s_2^{(t)}) - \nabla f(\x-r\s_2^{(t)})}{2r}  - \nabla^2 f(\x)\s_2^{(t)}\right\|\\
    & = \left\|\frac{\int_0^1\nabla^2 f(\x- r\s_2^{(t)} +2u r\s_2^{(t)}) 2r\s_2^{(t)} du}{2r}  - \nabla^2 f(\x)\s_2^{(t)}\right\|\\
    & \le \int_{0}^1\left\|\nabla^2 f(\x+(2u-1) r\s_2^{(t)})-\nabla^2 f(\x)\right\| du \\
    & \le rL_2 \int_{0}^1 |2u-1| du\\ &= \frac{r L_2}{2},
    \label{eq:error_small}
\end{align}
where the second equality follows directly from the fundamental theorem of calculus~(see e.g. the introductory chapter in~\cite{nesterov2018lectures}, proof of Lemma 1.2.2). The first inequality comes from Cauchy–Schwarz and the definition of operator norm, after noting that $\| \s_2^{(t)} \|=1$. The second inequality from Hessian Lipschitzness. Note that for quadratics $L_2=0$ so $\xi^(t)_2$ is identically zero.
\vspace{-2mm}
\paragraph{Bound on $\vxi_{1,\pm}^{(t)}$.} These error also vanish as $c\to0$, and the proof is similar to the one above. This was already shown e.g. in Lemma 3 (Appendix D) from~\cite{ji2019improved}. We give a proof for completeness, again based on the fundamental theorem of calculus. 
\begin{align}
    \g_+ &= \frac{1}{2c}\sum\limits_{i=1}^d \left(f(\x + r\s^{(t)}_2 + c  \e_i) - f(\x + r\s^{(t)}_2 - c   \e_i)\right) \e_i\\
    &= \sum\limits_{i=1}^d\e_i\e_i^\top  \int_{0}^1 \nabla f(\x + r\s^{(t)}_2 + (2u-1)c\e_i) du. \label{eq:starting_point_quadratic}
\end{align}

Therefore, using the subadditivity of the Euclidean norm and gradient Lipschitzness,
\begin{align}
\|\vxi_{1,+}^{(t)}\|^2 &= \left\|\sum\limits_{i=1}^d\e_i \e_i^\top  \int_{0}^1 \left(\nabla f(\x + r\s^{(t)}_2 + (2u-1)c\e_i) - \nabla f(\x + r\s^{(t)}_2)\right) du\right\|^2\\
&\le\sum_{i=1}^d\left(\int_{0}^1 \left\|\nabla f(\x + r\s^{(t)}_2 + (2u-1)c\e_i) - \nabla f(\x + r\s^{(t)}_2)\right\|du\right)^2\\
&\le \sum_{i=1}^d L^2c^2 \left(\int_0^1 |2u-1|\right)^2\\
&=\frac{dL^2 c^2}{4},
\end{align}
where the first inequality holds true because the vectors in the sum are mutually orthogonal and $\|\e_i\e_i^\top\|^2=1$.Note that here an additional dependency on the dimension comes in --- which is due to the triangle inequality and the nature of the estimator~(sum of $d$ terms). The same bound can be derived for $\vxi_{1,-}^{(t)}$. This concludes the proof.
\vspace{-2mm}
\paragraph{The quadratic case.} As mentioned above, in the quadratic case the Hessian is constant; hence $L_2=0$ and  therefore $\|\vxi_{2}^{(t)}\|=0$. However, from the bound above it seems that the bound on $\|\vxi_{1,\pm}^{(t)}\|$ does not vanish, since $L_1\ne 0$. This is an artefact of the proof technique. Indeed, for the quadratic case we have $\g_+=f(\x + r\s^{(t)}_2)$ and $\g_-=f(\x - r\s^{(t)}_2)$. This can be seen by inspecting the integral in Equation~\ref{eq:starting_point_quadratic}: assuming $f(\x) = C + (\x-\x^*)^\top \Hm (\x-\x^*)$ we have
\begin{align}
    \int_{0}^1 \nabla f(\x + r\s^{(t)}_2 + (2u-1)c\e_i) du    &= \int_0^1 \Hm(\x + r\s^{(t)}_2 + (2u-1)c\e_i-\x^*) du \\
    &= \Hm(\x + r\s^{(t)}_2-\x^*) + \Hm\int_0^1 (2u-1)c\e_i du \\
    &= \Hm(\x + r\s^{(t)}_2-\x^*)\\
    &= \nabla f(\x + r\s^{(t)}_2).
\end{align}
This concludes the proof.
\end{proof}

\subsection{Lemma~\ref{lemma:Hardt} and results on convergence of (noisy) power methods}

Finding the smallest eigenvalue (assumed to be negative) of the Hessian $\nabla^2 f(\x_t)$) is equivalent to the one of finding the largest eigenvalue of $\Am = \Im - \eta \nabla^2 f(\x_t)$, where $\eta$ is a small positive number such that $\eta\le1/\|\nabla^2 f(\x_t)\|$ (the spectral norm of $\nabla^2 f(\x_t)$). For this choice of $\eta$, $\Im - \eta \nabla^2 f(\x_t)$ is positive semidefinite, hence one can use an (inexact) power method to retrieve the maximum eigenvalue. We first present the standard error analysis of the power iteration~(which we adapt from~\cite{golub2013matrix}), assuming we have access to the true Hessian. Then, we discuss the setting where we can only compute approximate Hessian-vector products~(analysis adapted from~\cite{hardt2014noisy}). Finally, we present the bound for the Derivative-Free Power Iteration (DFPI) algorithm~(Alg.~\ref{alg:inexact_PI}).

\subsubsection{Warm-up: error analysis for the exact power method}
Let $\Am\in\R^{d\times d}$ be a positive definite matrix with eigenvalues $a_1>a_2\ge\dots a_d>0$, and corresponding eigenvectors $\v_1, \v_2,\dots,\v_d$. Eigenvalues are counted together with their algebraic multiplicity. We seek an approximation for the dominant eigendirection $\v_1$. The power method on the positive semidefinite matrix $\Am$ can be found as Algorithm~\ref{alg:power_iteration}. 

\begin{algorithm}[ht] 
\begin{algorithmic}[1]
\STATE {\textbf{INPUT :} A matrix $\Am\in\R^{d\times d}$ with eigenvalues $a_1>a_2\ge\dots a_d>0$.}
\STATE Randomly initialize $\v^{(0)} \sim \mathcal{S}^{d-1}$\\
\FOR{$t = 0 \dots T-1$}
    \STATE $\v^+ = \Am \v^{(t)}$\\
    \STATE $\v^{(t+1)} =\v^+/\|\v^+\|$ 
\ENDFOR
\STATE {\textbf{OUTPUT :} $\v^{(T)}$ approximating $\v_1$, leading eigenvector of $\Am$. }
\end{algorithmic}
\caption{{\sc Power Method (exact, access to Hessian-vector products required)}}
\label{alg:power_iteration}
\end{algorithm}

We present the fundamental yet simple result, confirming that the power iteration step decreases the distance to the dominant eigendirection. We recall that $\angle(\v,\u) := \arccos{\frac{\langle \v,\u \rangle}{\|\u\|\cdot\|\v\|}}$.

\begin{tcolorbox}
\begin{lemma}
Consider a step of Alg.~\ref{alg:power_iteration}), $\tan(\angle(\v^{(t+1)},\v_1)) \le \frac{a_2}{a_1}\tan(\angle(\v^{(t)},\v_1))$.
\label{lemma:power}
\end{lemma}
\end{tcolorbox}
\begin{proof}
First, we write $\v^{(t)}$ in the eigenbasis $\{\v_i\}_{i=1}^d$: $\v^{(t)} = \sum_{i=1}^d \alpha^{(t)}_i \v_i$. Crucially, note that
\begin{equation}
    \tan(\angle(\v^{(t)},\v_1)) = \frac{\sqrt{\sum_{j=2}^d(\alpha_j^{(t)})^2}}{\alpha^{(t)}_1}.
    \label{eq:tangent}
\end{equation}
Since $\v^+= \sum_{i=1}^d a_i\alpha^{(t)}_i\v_i$, we have that 
\vspace{-3mm}
\begin{equation}
    \tan(\angle(\v^{(t+1)},\v_1)) = \tan(\angle(\v^{+},\v_1)) = \frac{\sqrt{\sum_{j=2}^d a_j^2(\alpha_j^{(t)})^2}}{a_1\alpha^{(t)}_1}\le\frac{a_2}{a_1}\tan(\angle(\v^{(t)},\v_1)).
\end{equation}
\end{proof}
\vspace{-2mm}
As noted by~\cite{hardt2014noisy}, the dependence on the eigenvalue separation arises already in the classical perturbation argument of Davis-Kahan~\cite{davis1970rotation}. If $a_1$ has multiplicity greater than 1, then of course the ratio will be $a_{k}/a_1$, where $a_k$ is the first eigenvalue strictly smaller than $a_1$. More on this point can be found in Remark~\ref{rmk:eigengap}.

From the lemma above, we can easily deduce the error on the eigenvalue computation

\begin{tcolorbox}
\begin{theorem}
\label{thm:power_iteration}
Algorithm~\ref{alg:power_iteration} outputs a vector $\v^{(T)}$ such that $|(\v^{(T)})^\top \Am \v^{(T)} - a_1|\le\epsilon a_1$ if
\begin{equation}
    T\ge\frac{a_1}{2(a_1-a_2)}\log \left(\frac{\tan^2(\angle(\v^{(0)},\v_1))}{\epsilon}\right).
\end{equation}
Moreover, as also mentioned in Lemma~2.5 in~\cite{hardt2014noisy} and Lemma 2.2 in~\cite{balcan2016improved}, if $\v^(0)$ is randomly initialized on the unit sphere, the main result in~\cite{rudelson2009smallest} implies that with probability $1-\delta-e^{\Omega(d)}$ we have
$\tan^2(\angle(\v^{(0)},\v_1))\le  d/\delta^2 $. Hence, with probability $1-\delta-e^{\Omega(d)}$, we have
\begin{equation}
    T\ge\frac{a_1}{2(a_1-a_2)}\log \left(\frac{d}{\epsilon\delta^2}\right).
\end{equation}
\end{theorem}
\end{tcolorbox}
\vspace{-4mm}
\begin{proof}
Note that since $\v^{(T)}$ is normalized,
\begin{equation}
    \tan^2(\angle(\v^{(T)},\v_1))^2=\frac{\sum_{j=2}^d(\alpha^{(T)}_j)^2}{(\alpha^{(T)}_1)^2}=\frac{1-(\alpha^{(T)}_1)^2}{(\alpha^{(T)}_1)^2},
\end{equation}
therefore
\begin{equation}
    (\alpha^{(T)}_1)^2=\frac{1}{1+\tan^2(\angle(\v^{(T)},\v_1))},\quad\quad\sum_{j=2}^d(\alpha^{(T)}_j)^2=\frac{\tan^2(\angle(\v^{(T)},\v_1))^2}{1+\tan^2(\angle(\v^{(T)},\v_1))}.
\end{equation}

We have the following bound:
\begin{align}
    |(\v^{(T)})^\top \Am \v^{(T)} - a_1|&=\left|a_1(\alpha^{(T)}_1)^2 + \sum_{i=2}^d a_i(\alpha^{(T)}_i)^2 - a_1\right|\\
    &= a_1-a_1(\alpha^{(T)}_1)^2 - \sum_{i=2}^d a_i(\alpha^{(T)}_i)^2\\ 
    &\le a_1-a_1\frac{1}{1+\tan^2(\angle(\v^{(T)},\v_1))} - a_d \frac{\tan^2(\angle(\v^{(T)},\v_1))}{1+\tan^2(\angle(\v^{(T)},\v_1))}\\
    &= \frac{\tan^2(\angle(\v^{(T)},\v_1))}{1+\tan^2(\angle(\v^{(T)},\v_1))}(a_1-a_d)\\
    &\le a_1 \tan^2(\angle(\v^{(T)},\v_1))
    \label{eq:power_to_tangent}
\end{align}
where the second equality is given by the fact that $a_1$ is the biggest eigenvalue of $\Am$. All in all, we need $\tan^2(\angle(\v^{(T)},\v_1))$ to be smaller than $\epsilon/(a_1-a_d)$.
Thanks to Lemma~\ref{lemma:power}, we have that
\begin{equation}
    \tan^2(\angle(\v^{(T)},\v_1))\le\left(\frac{a_2}{a_1}\right)^{2T}\tan^2(\angle(\v^{(0)},\v_1)).
\end{equation}
Therefore, we require $\left(\frac{a_2}{a_1}\right)^{2T} a_1\tan^2(\angle(\v^{(0)},\v_1))\le\epsilon a_1$, which can be written as,
\begin{equation}
    \left(\frac{a_1}{a_2}\right)^{2T}\ge\frac{\tan^2(\angle(\v^{(0)},\v_1))}{\epsilon}.
\end{equation}
\vspace{-2mm}
We conclude by taking the $\log$ on both sides:
\begin{equation}
    T\ge\frac{1}{2\log(a_1/a_2)}\log\left(\frac{\tan^2(\angle(\v^{(0)},\v_1))}{\epsilon}\right)
\end{equation}
Since, for all $x\in\R$, $\log(x)\ge 1-\frac{1}{x}$ and $a_1>a_2$, the above expression is verified if
\vspace{-2mm}
\begin{equation}
    T\ge\frac{a_1}{2(a_1-a_2)}\log\left(\frac{\tan^2(\angle(\v^{(0)},\v_1))}{\epsilon}\right)
\end{equation}
\vspace{-2mm}
\end{proof}

\begin{remark}
Note that Theorem~\ref{thm:power_iteration} is exactly equivalent to Theorem 8.2.1 in~\cite{golub2013matrix}. Here we followed a proof more similar to the one in~\cite{hardt2014noisy}.
\end{remark}

\begin{remark}[Eigen-gap dependency]
\label{rmk:eigengap}
The bound in Theorem~\ref{thm:power_iteration} depends on the eigen-gap $a_1-a_2$: as $a_1$ and $a_2$ get closer, the result suggests that we need a very large number of iterations to find a good approximation of $\v_1$. This is true because the power method is confounded by $\v_2$, and takes a long time to ``decide'' which one between $\v_1$ and $\v_2$ is dominant. However, this of course does not imply that the complexity in finding $\v^{(T)}$ such that $|(\v^{(T)})^\top \Am \v^{(T)} - a_1|\le\epsilon$ increases --- this is an artefact of our simple analysis~(inspired by~\cite{hardt2014noisy,golub2013matrix,balcan2016improved}), which crucially goes through Lemma~\ref{lemma:power} to derive the bound. Indeed, as the next theorem shows, \textbf{it is possible to directly remove this dependency}.
\end{remark}

\begin{tcolorbox}
\begin{theorem}[Consequence of Thm.~3.1 and Thm.~4.1 in~\cite{kuczynski1992estimating}] Let $\v^{(0)}$ be initialized randomly on the surface of the unit sphere. The power method returns a vector $\v^{(T)}$ such that $|(\v^{(T)})^\top \Am \v^{(T)} - a_1|\le\epsilon a_1$ in $T=\mathcal{O}(\log(d)/\epsilon)$ iterations, in expectation. For the result to hold with probability $1-\delta$, one instead needs at least $T=\mathcal{O}(\log(d/\delta^2)/\epsilon)$ iterations.
\label{thm:kuc}
\end{theorem}
\end{tcolorbox}
This results in also cited in~\cite{xu2018first}, where the bound above is used to conclude that, if $\lambda_\text{min}(\nabla^2 f(\x))\le-\gamma$ and $\|\nabla^2 f(\x)\|\le L_1$, the power method on $(\Im-\eta\nabla^2 f(\x))$ finds a direction $\v^{(T)}$ such that, with probability $1-\delta$, $(\v^{(T)})^\top\nabla^2 f(\x)\v^{(T)}\le-\frac{\gamma}{2}$ in $\mathcal{O}\left(\frac{L_1}{\gamma}\log(d/\delta^2)\right)$ iterations.

This proves directly a version of Lemma~\ref{lemma:Hardt} for the case of vanishing error.
\begin{tcolorbox}
\begin{lemma}[Noiseless version of Lemma~\ref{lemma:Hardt}]
Let the parameters of DFPI be such that the error $\vxi_\DFPI$ is vanishing~(possible within the limits of numerical stability by Lemma~\ref{prop:DFPI_error}). Let $\eta\le 1/L_1$. Let $\gamma=\epsilon^{2/3}$; for a fixed RSPI iteration, $T_\DFPI =\mathcal{O}\left(\epsilon^{-2/3}L_1\log\left(\frac{d}{\delta^2}\right)\right)$ DFPI iterations are enough to ensure validity of Assumption~\ref{ass:error_eigenvector} (without the expectation sign) at $\x_k$ with probability $1-\delta-e^{\Omega(d)}$.
\end{lemma}
\end{tcolorbox}
\begin{proof}Direct consequence of the reasoning above, supported by Lemma~\ref{prop:DFPI_error}.
\end{proof}

\subsubsection{Error analysis for the noisy power method}

We now consider the case where $\Am \v$ cannot be computed exactly~(Algorithm~\ref{alg:power_iteration_noisy}): we denote by $\vxi^{(t)}$ the error in computing the Hessian-vector product $\Am\v^{(t)}$.

\begin{algorithm}[ht] 
\begin{algorithmic}[1]
\STATE {\textbf{INPUT :} A matrix $\Am$ with eigenvalues $a_1>a_2\ge\dots a_d$.}
\STATE Randomly initialize $\v^{(0)} \sim \mathcal{S}^{d-1}$\\
\FOR{$t = 0 \dots T-1$}
    \STATE $\v^+ = \text{approx}(\Am\v^{(t)})=\Am \v^{(t)} + \vxi^{(t)}$\\
    \STATE $\v^{(t+1)} =\v^+/\|\v^+\|$ 
\ENDFOR
\STATE {\textbf{OUTPUT :} $\v^{(T)}$, approximating $\v_1$, leading eigenvector of $\Am$ }
\end{algorithmic}
\caption{{\sc Power Method (noisy, approximate Hessian-vector products permitted)}}
\label{alg:power_iteration_noisy}
\end{algorithm}

We are now ready to state the main result we are going to use on the noisy power method, presented in the main text in a less precise way, as Lemma~\ref{lemma:complexity_power}. This result was first derived in~\cite{hardt2014noisy}, and can be seen as an extension to Theorem~\ref{thm:power_iteration}. In plain english: \textbf{for small enough noise}, the bound in Theorem~\ref{thm:power_iteration} still holds with arbitrarily high probability.

\begin{tcolorbox}
\begin{restatable}[Direct consequence of Corollary 1.1 in~\cite{hardt2014noisy}]{theorem}{InexactPowerIteration}\label{thm:noisy_power_iteration}
In the context of Algorithm~\ref{alg:power_iteration_noisy}, fix the desired accuracy $\epsilon\le1/2$ and a failure probability $\delta$. Assume that for all iterations $t$ the noise is small enough: (1) $5\|\vxi^{(t)}\|\le \epsilon (a_1-a_2)$ and (2) $5|\v_1^\top\vxi^{(t)}|\le \delta(a_1-a_2)/\sqrt{d}$. With probability $1-\delta-e^{-\Omega(d)}$, Algorithm~\ref{alg:power_iteration_noisy} returns $\v^{(T)}$ such that $|(\v^{(T)})^\top \Am \v^{(T)} - a_1|\le\epsilon a_1$ if
\begin{equation}
    T\ge\mathcal{O}\left(\frac{a_1}{a_1-a_2}\log\left(\frac{d}{\epsilon\delta^2}\right)\right).
\end{equation}
\end{restatable}
\end{tcolorbox}
\begin{proof}
In the proof of Theorem~\ref{thm:power_iteration}, we showed that
\begin{equation}
    |(\v^{(T)})^\top \Am \v^{(T)} - a_1|\le \tan^2(\angle(\v^{(T)},\v_1))\cdot(a_1-a_d).
\end{equation}
This is enough to complete the result given Corollary 1.1 in~\cite{hardt2014noisy}. 
\end{proof}

The proof of Lemma~\ref{lemma:Hardt} then follows from a generalization of Theorem~\ref{thm:kuc} to the noisy case~(under the requirement of small enough noise). This is possible since the bounds in Theorem~\ref{thm:noisy_power_iteration} and Theorem~\ref{thm:power_iteration} are equivalent --- meaning that the geometry of convergence is not drastically affected by noise.

\subsection{How to speed up DFPI with SPSA: an experimental motivation}
\begin{figure}[ht!]
    \centering
    \includegraphics[width=0.49\textwidth]{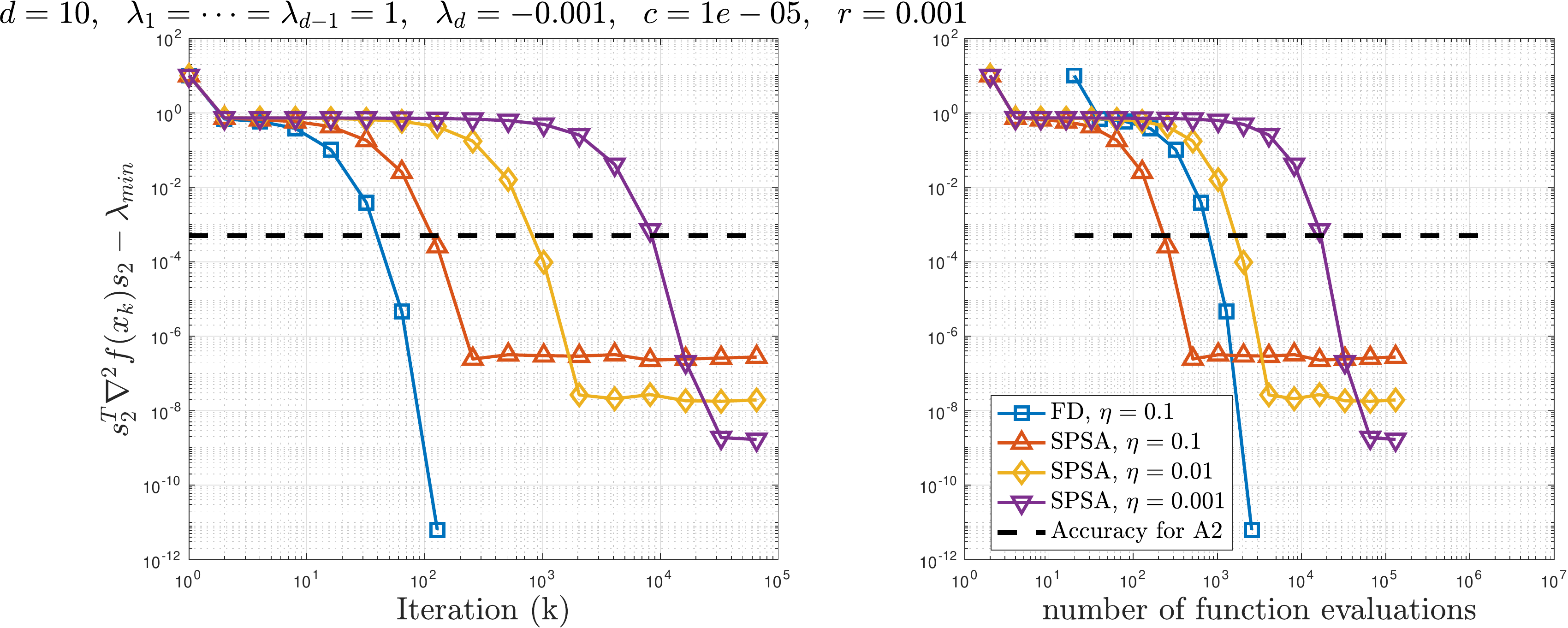}
    \includegraphics[width=0.49\textwidth]{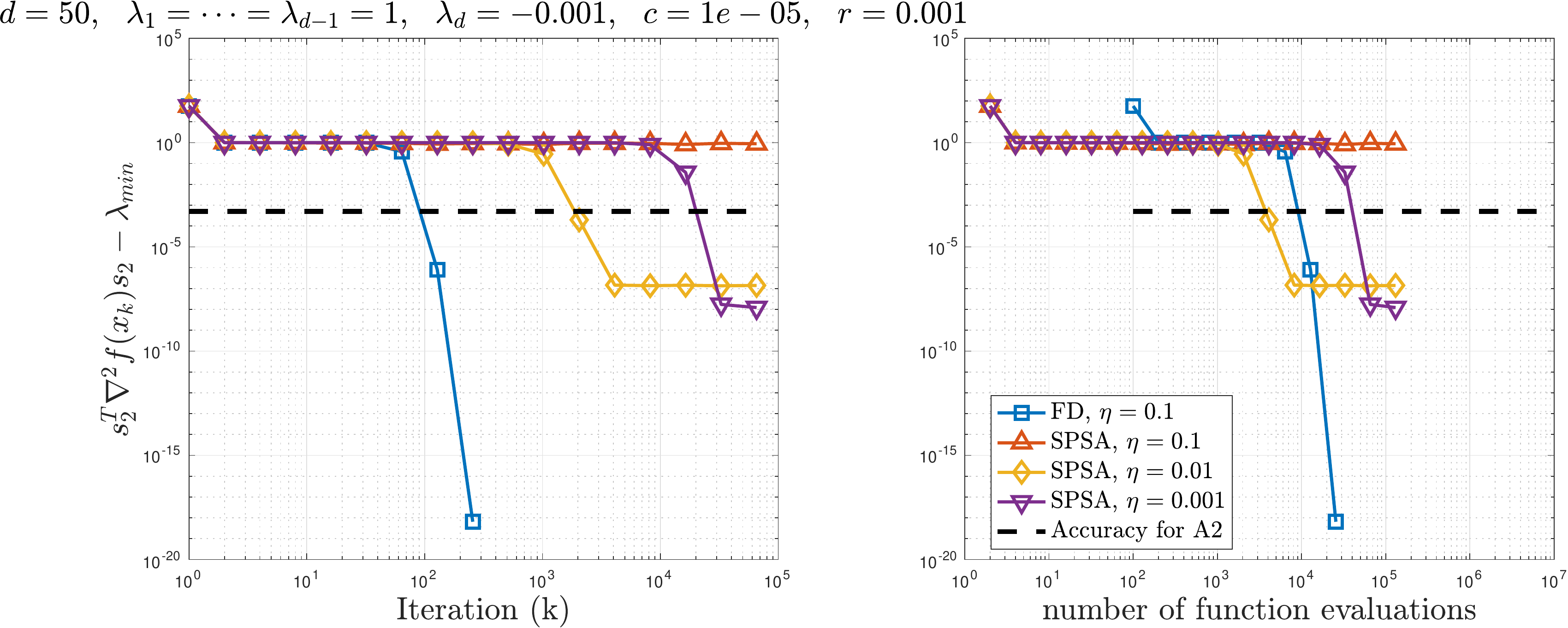}
    \includegraphics[width=0.49\textwidth]{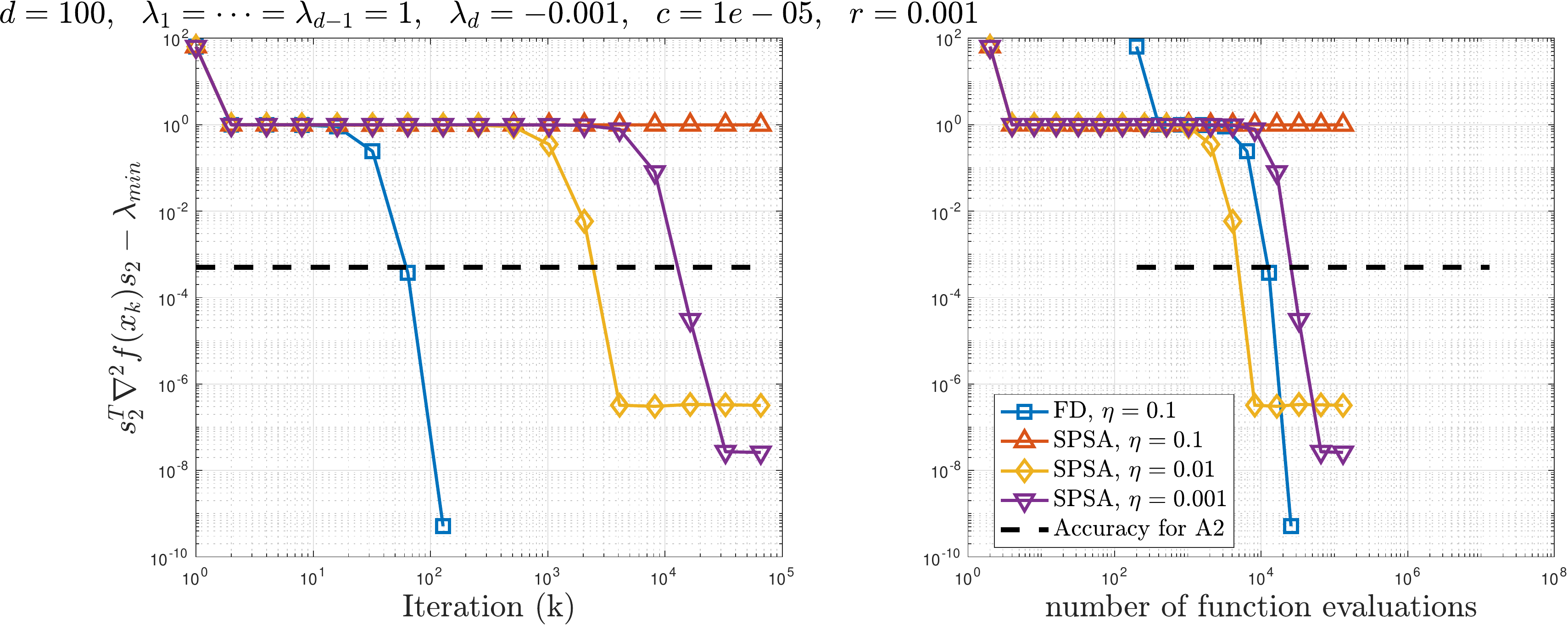}
    \includegraphics[width=0.49\textwidth]{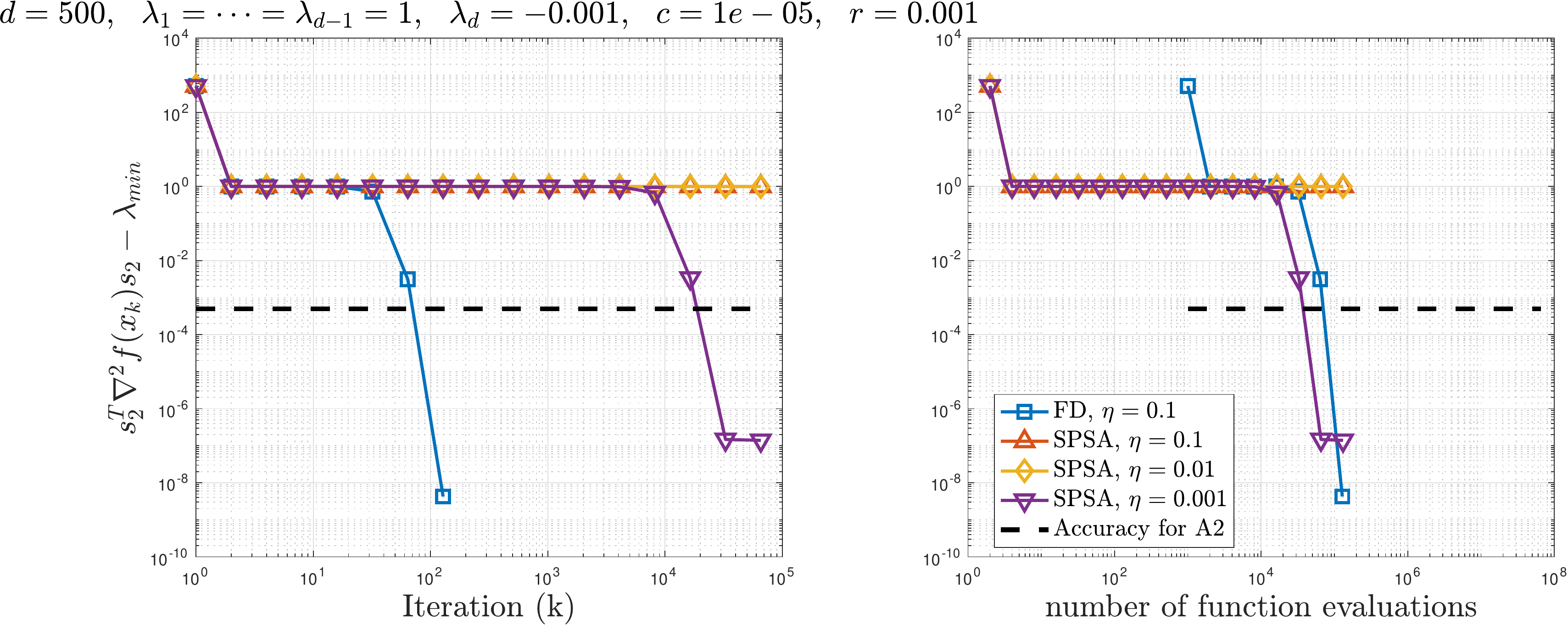}
    \caption{Experiment 1: $f(\x) = \x^\top diag(\lambda_1,\lambda_2,\cdots,\lambda_d)\x$, $\lambda_d = -0.001$. Settings described in the paragraph above.}
    \label{fig:SPSA1}
\end{figure}

\begin{figure}[ht!]
    \centering
    \includegraphics[width=0.49\textwidth]{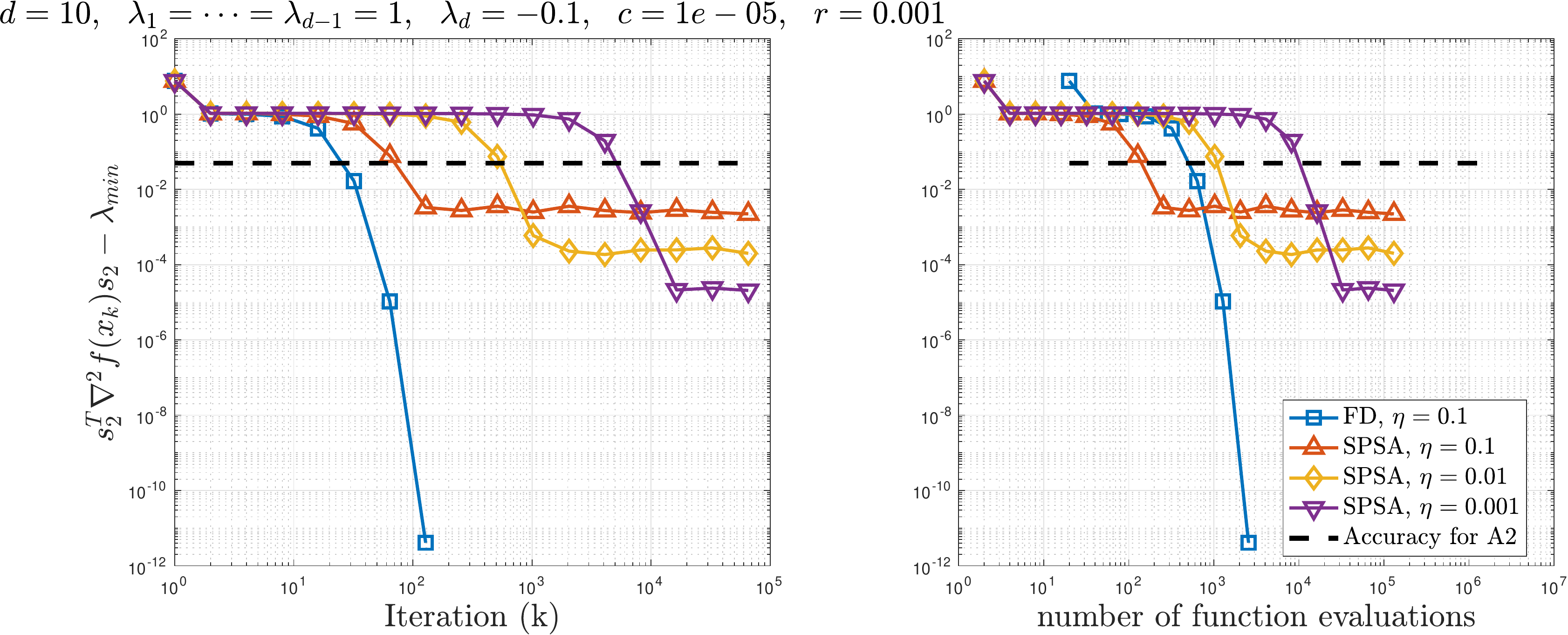}
    \includegraphics[width=0.49\textwidth]{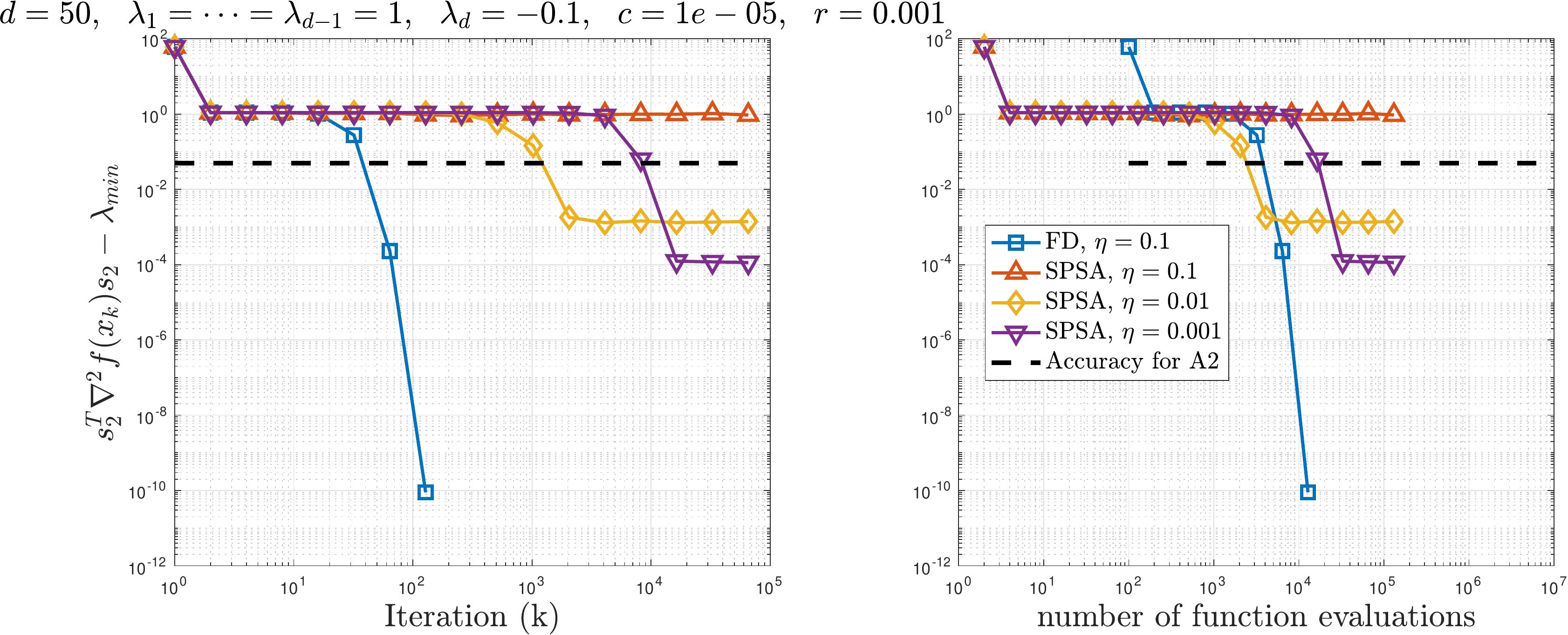}
    \includegraphics[width=0.49\textwidth]{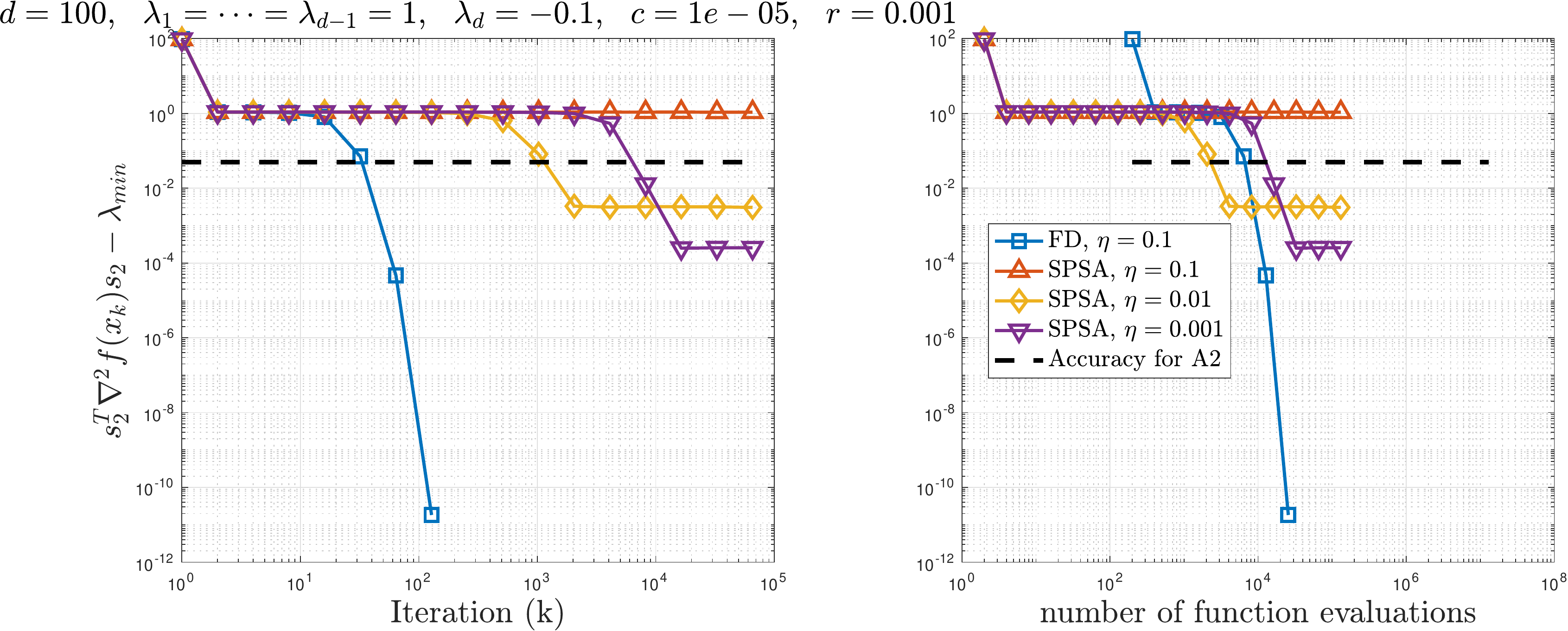}
    \includegraphics[width=0.49\textwidth]{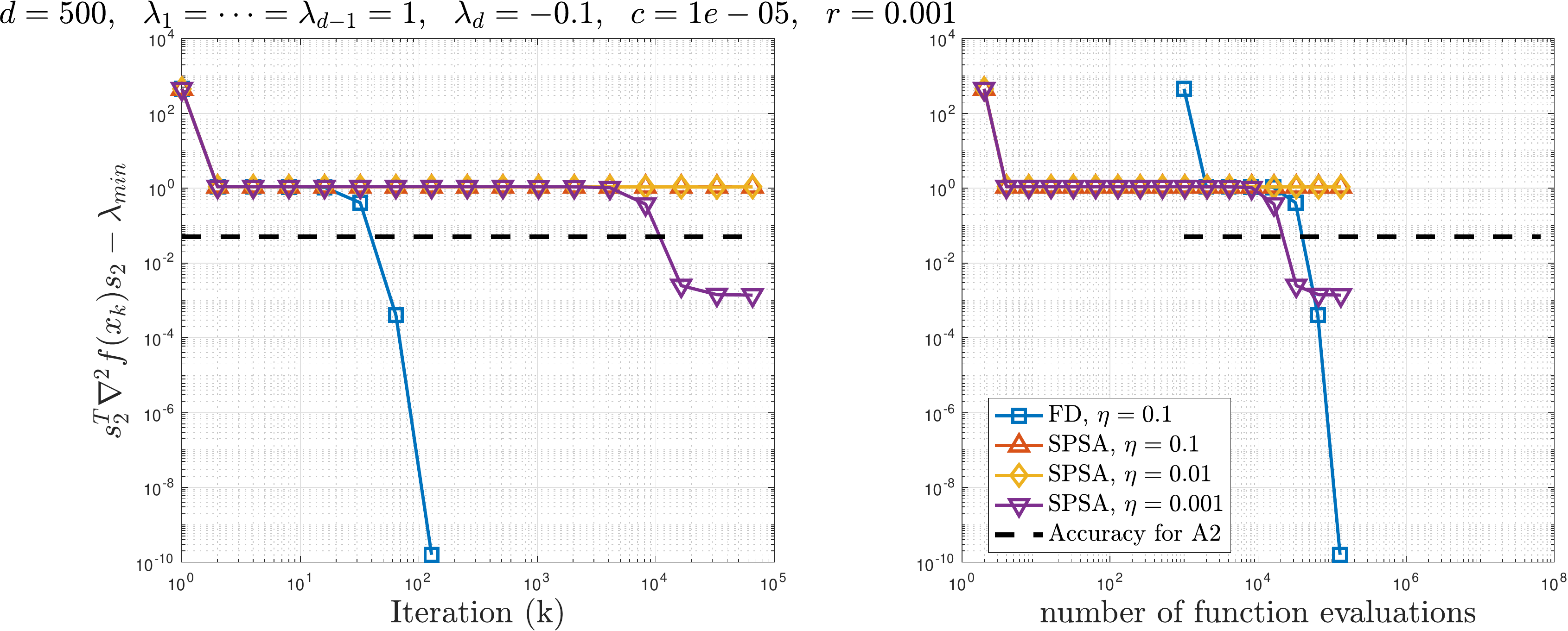}
    \caption{Experiment 2: $f(\x) = \x^\top diag(\lambda_1,\lambda_2,\cdots,\lambda_d)\x$, $\lambda_d = -0.1$. Settings described in the paragraph above.}
    \label{fig:SPSA2}
\end{figure}

\begin{figure}[ht!]
    \centering
    \includegraphics[width=0.49\textwidth]{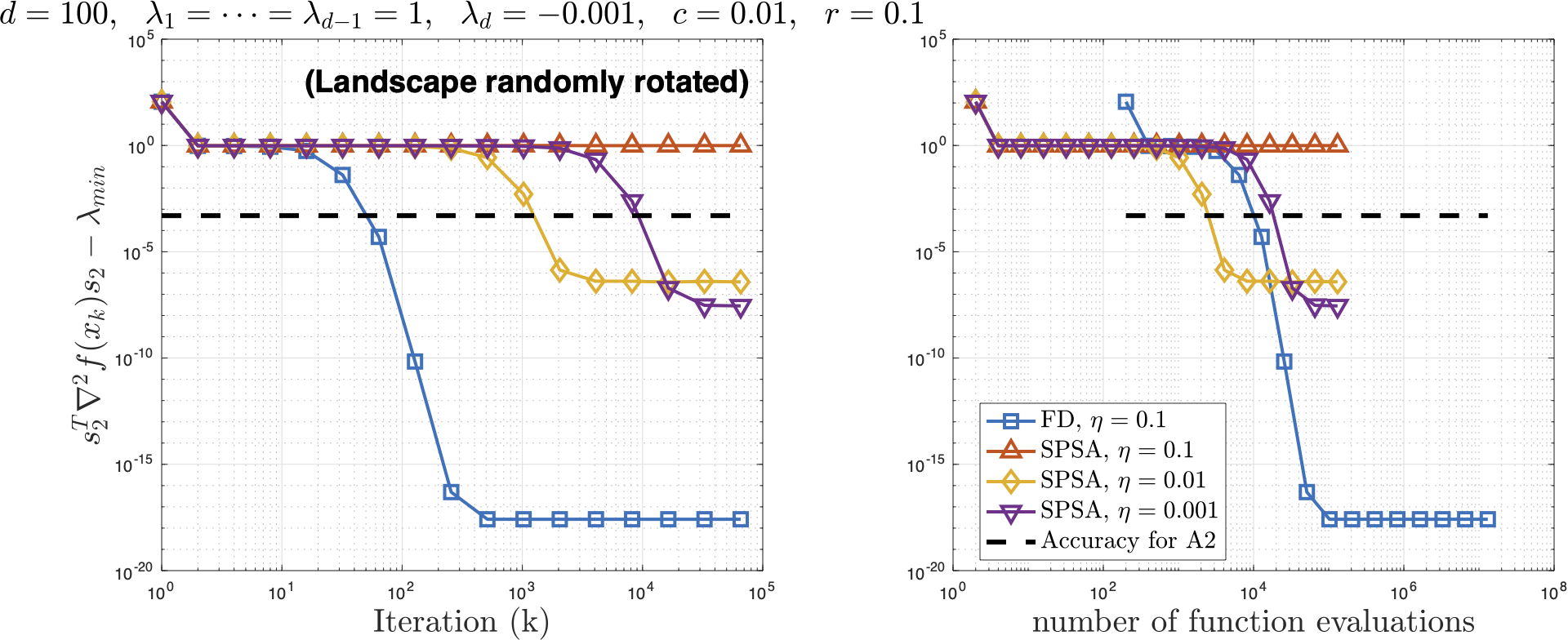}
    \includegraphics[width=0.49\textwidth]{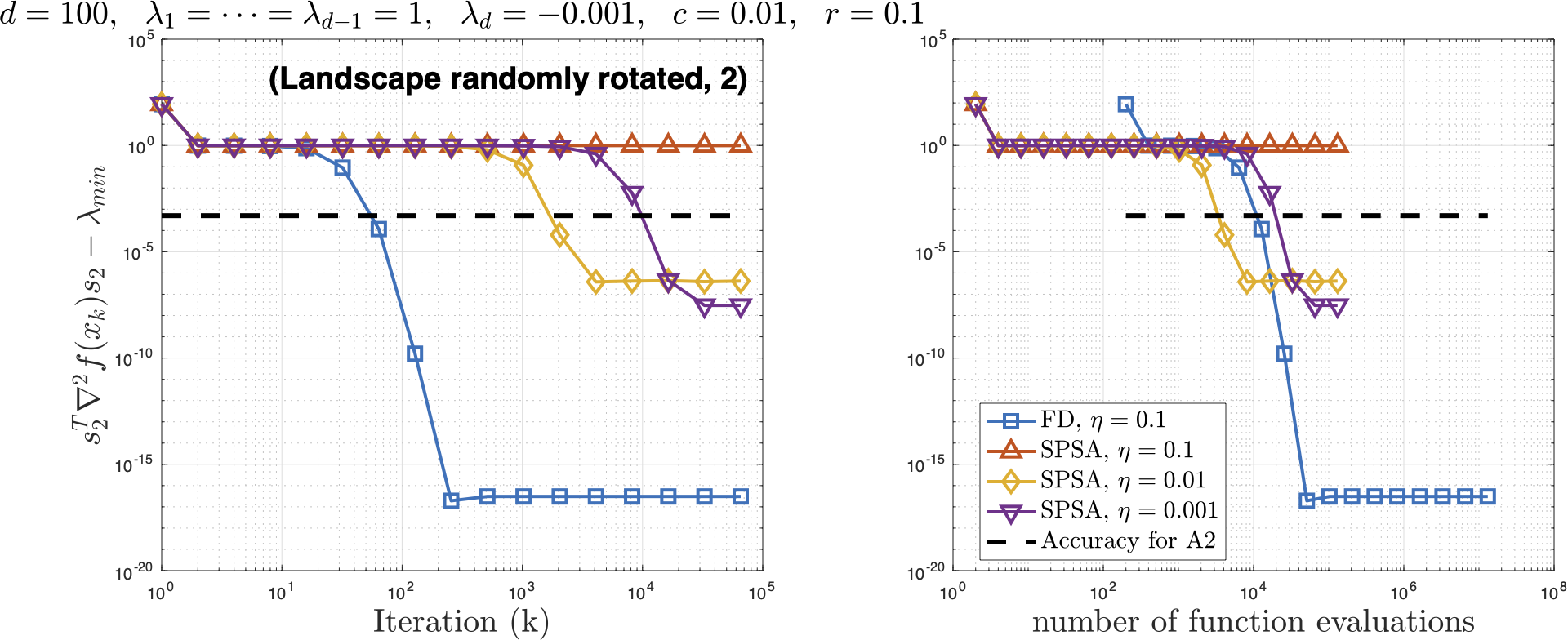}
    \caption{Experiment 2: $f(\x) = \x^\top \boldsymbol{U}^\top diag(\lambda_1,\lambda_2,\cdots,\lambda_d) \boldsymbol{U}\x$, $\lambda_d = -0.1$, where $\boldsymbol{U}$ is a random orthogonal matrix. Dynamics for two two different random orthogonal matrices are shown, where we additionally also decreased $c$ and $r$. The evolution is similar to the one in Figure~\ref{fig:SPSA1}, showing that SPSA is robust to both landscape rotations and hyperparameter choice.}
    \label{fig:SPSA3}
\end{figure}

We study some interesting properties of the SPSA gradient estimator, introduced by~\cite{spall1992multivariate}, in the context of DFPI~(Algorithm~\ref{alg:inexact_PI}, main paper). In particolar, we consider using SPSA instead of finite-difference(FD), which is the base for our theory~(Thm.~\ref{thm:conv_random_PI})
\vspace{-2mm}
\begin{align*}
    &\g^{\text{SPSA}}_{\pm} = \sum\limits_{i=1}^d \frac{f(\x \pm r\s^{(t)}_2 + c \vDelta) - f(\x \pm r\s^{(t)}_2 - c \vDelta)}{2c \Delta_i} \e_i.\\ &\quad \quad\quad\quad\quad\quad \quad\quad\quad\to \text{\textbf{$\boldsymbol{4}$ function evaluations} to get estimates of $\nabla f(\x\pm\s_2^{(t)})$.}\\
    & \g^{\text{FD}}_+  \ \ \ \ = \sum_{i=1}^d \frac{f(\x \pm r\s^{(t)}_2 + c \e_i) - f(\x \pm r\s^{(t)}_2 - c \e_i)}{2c} \e_i\\ &\quad \quad\quad\quad\quad \quad\quad\quad\quad \to \text{\textbf{$\boldsymbol{4d}$ function evaluations} to get estimates of $\nabla f(\x\pm\s_2^{(t)})$.}
\end{align*}
\paragraph{The SPSA estimator is asymptotically unbiased, but variance might be independent of the hyperparamerter $\boldsymbol{c}$.} Consider  $f(x_1,x_2) = x_1^2-x_2^2$, we want to approximate its gradient using SPSA. Since perturbation is $(\Delta_1,\Delta_2)$, we have $f(\x+c\vDelta) - f(\x-c\vDelta) = 4c\Delta_1x_1-\Delta_2x_2$. Therefore $\g^{\text{SPSA}} = \sum_{i=1}^2 \frac{f(\x+c\vDelta) - f(\x-c\vDelta)}{2c\Delta_i}\e_i =\sum_{i=1}^2 \frac{2\Delta_1x_1-2\Delta_2x_2}{\Delta_i}\e_i$. Since $\Delta_i$ are Bernoulli, then $\E[\g^{\text{SPSA}}] = \nabla f$. However, the estimator variance is finite and independent of $c$.
\vspace{-2mm}
\paragraph{Experimental comparison.} From the result in the paragraph above, one might conclude that SPSA cannot provide a satisfactory approximation of Hessian-vector products, and therefore cannot be used as a valid alternative to FD in the context of an approximate power method such as DFPI. However, in Figure~\ref{fig:SPSA1} \ \&~\ref{fig:SPSA2} we show that, for small enough $\eta$, the update $\s^{(t+1)}_2 = \s^{(t)}_2 - \eta \frac{\g^{\text{SPSA}}_+ - \g^{\text{SPSA}}_-}{2r}, \s^{(t+1)}_2 = \s^{(t+1)}_2/\|\s^{(t+1)}_2\|$ can effectively build a vector $\s_2$ aligned with negative curvature, even as the problem dimension increases. In these experiments, we consider applying DFPI to estimate the negative curvature direction $\e_d$ of $f(\x) = \x^\top diag(\lambda_1,\lambda_2,\cdots,\lambda_d)\x$, with $\lambda_d<0$~(non-axis aligned case discussed later). As we saw in Prop.~\ref{prop:conv_random_PI_it}, the finite difference estimator yields an exact power method on this function. Instead, SPSA yields an inexact power method where the error is independent of $r, c$~(see last paragraph and Fig.~\ref{fig:SPSA3}). As expected, SPSA does not actually converge to the leading eigenvector. However, it can always be tuned to yield an approximation which satisfies Assumption~\ref{ass:error_eigenvector}, in a total number of function evaluations which is actually smaller than FD. Further research is needed to better understand this phenomenon. However, this motivates the use of SPSA as a cheap alternative to FD in DFPI. In the experiments section of the main paper, we indeed show that this approximation is enough to yield a satisfactory improvement over vanilla method which do not consider computing negative curvature. As can be evinced from the last paragraph and from the proof of Prop.~\ref{prop:conv_random_PI_it}, the results in this case are independent of the values of $r$ and $c$; however, they could in principle depend on the landscape rotation. We show in Figure~\ref{fig:SPSA3} that this is not the case using two random rotations.

\clearpage

\section{Experimental Results}\label{app:exp}

All of our experiments are conducted on the Google Colaboratory \cite{bisong2019colab} environment without any hardware accelerators. 

\subsection{Function with growing dimension}

\begin{figure}[!ht]
    \centering
    \includegraphics[width=\textwidth]{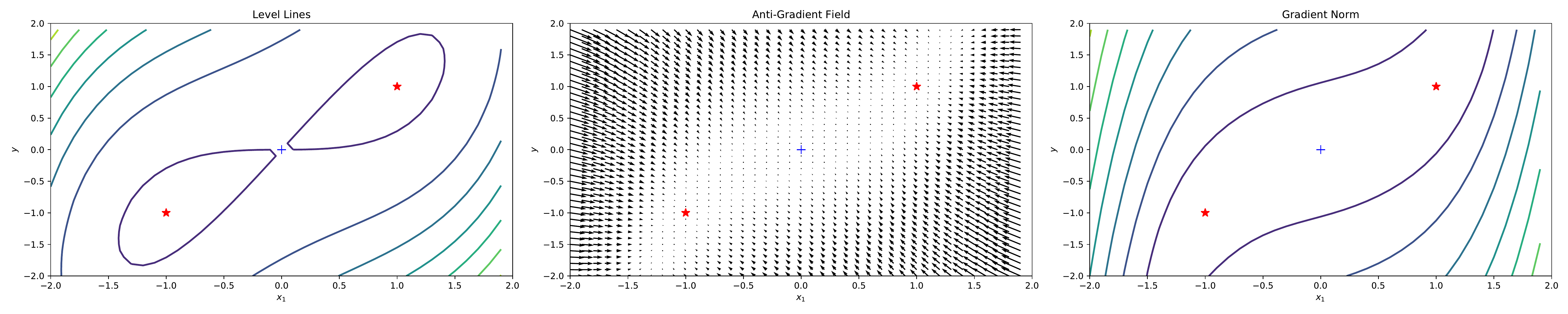}
    \caption{The landscape of the objective $f(x_1, \cdots, x_d, y) = \frac{1}{4}\sum_{i=1}^d x_i^4 - y\sum_{i=1}^d x_i + \frac{d}{2}y^2$ for $d=1$. A blue cross denotes a strict saddle point, whereas a red star corresponds to a global minimizer.}
    \label{fig:sb-landscape-vis}
\end{figure}

\begin{figure}[!ht]
    \centering
    \includegraphics[width=\textwidth]{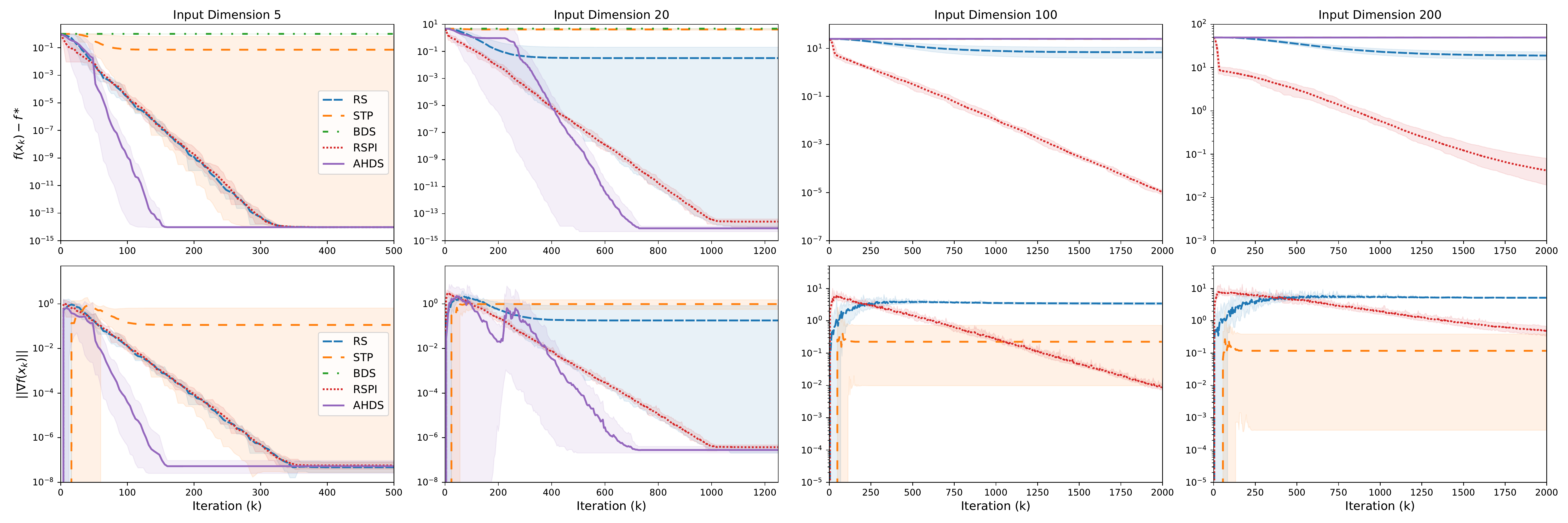}
    \caption{Empirical performance while minimizing $f(x_1, \cdots, x_d, y) = \frac{1}{4}\sum_{i=1}^d x_i^4 - y\sum_{i=1}^d x_i + \frac{d}{2}y^2$ against the number of iterations. Confidence intervals show min-max intervals over ten runs. All algorithms are initialized at the strict saddle point across all runs. }
\end{figure}

\begin{figure}[!ht]
    \centering
    \includegraphics[width=\textwidth]{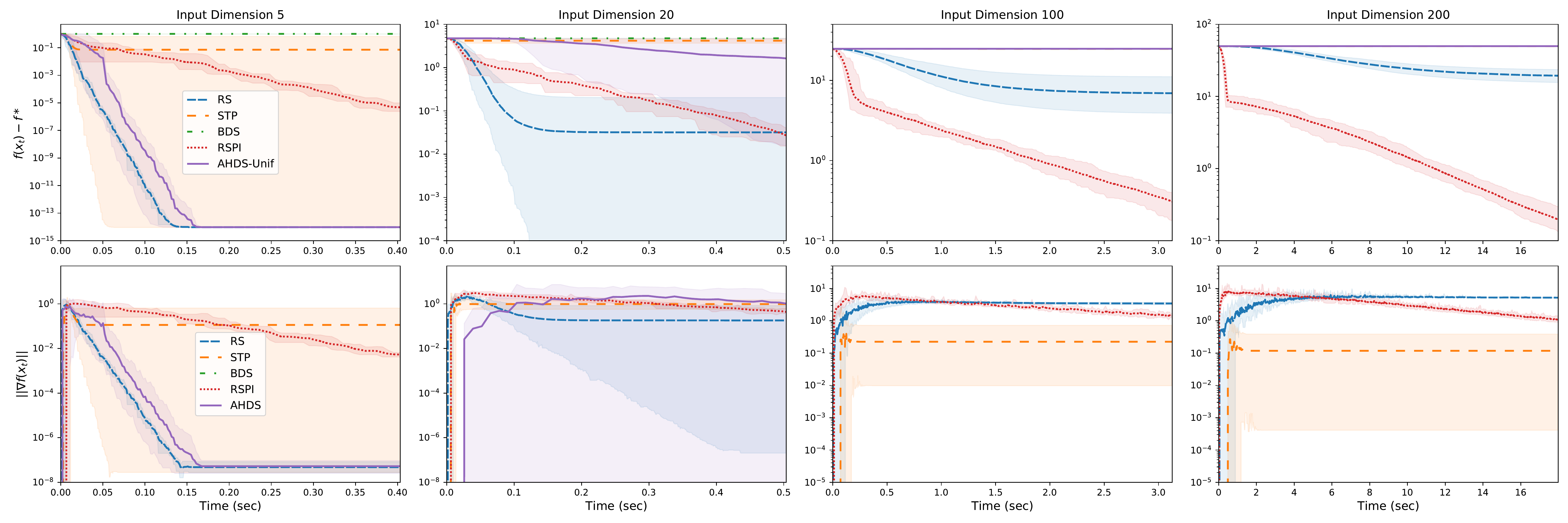}
    \caption{Empirical performance while minimizing the objective defined in the main paper against wall-clock time. Confidence intervals show min-max intervals over ten runs. All algorithms are initialized at the strict saddle point across all runs.}
\end{figure}

\clearpage
\subsection{Rastrigin function}

\begin{figure}[!ht]
    \centering
    \includegraphics[width=\textwidth]{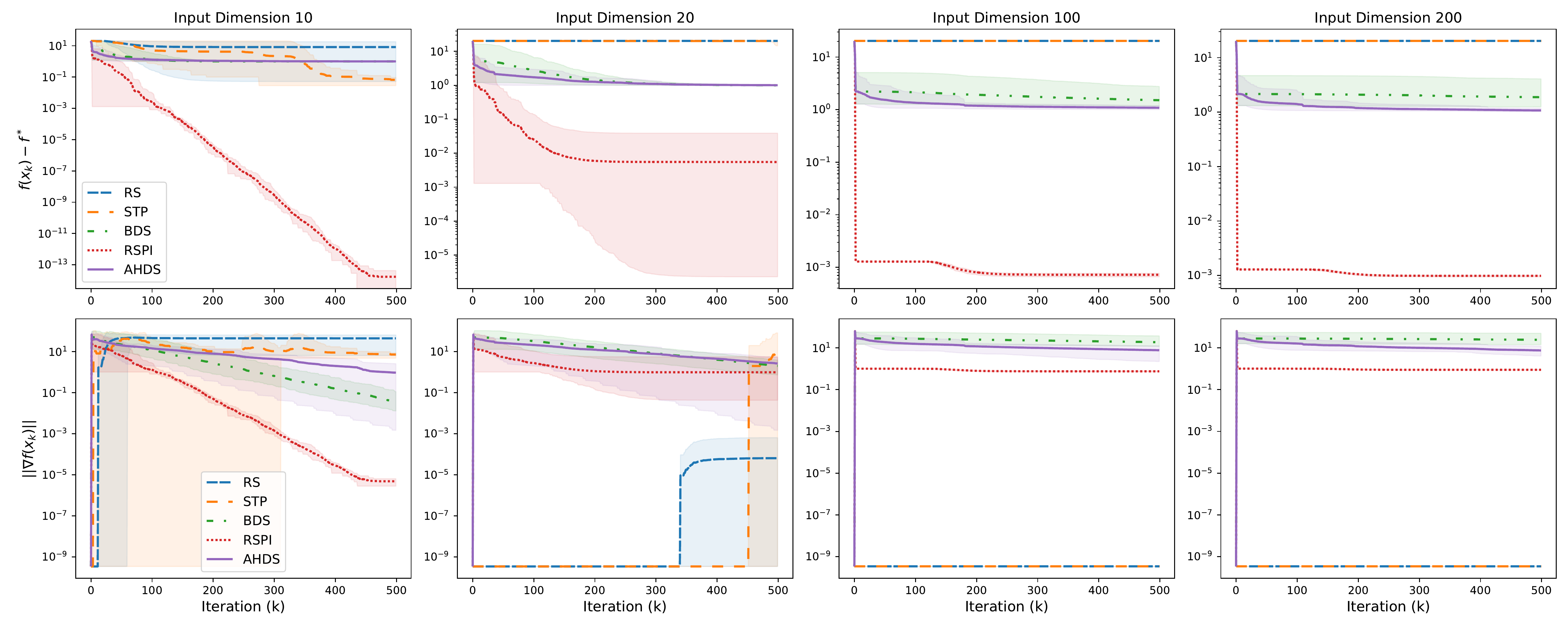}
    \caption{Empirical performance while minimizing the Rastrigin function against the number of iterations. Confidence intervals show min-max intervals over ten runs. All algorithms are initialized at a strict saddle point across all runs.}
    \label{fig:rastrigin-appendix-iter}
\end{figure}

\begin{figure}[!ht]
    \centering
    \includegraphics[width=\textwidth]{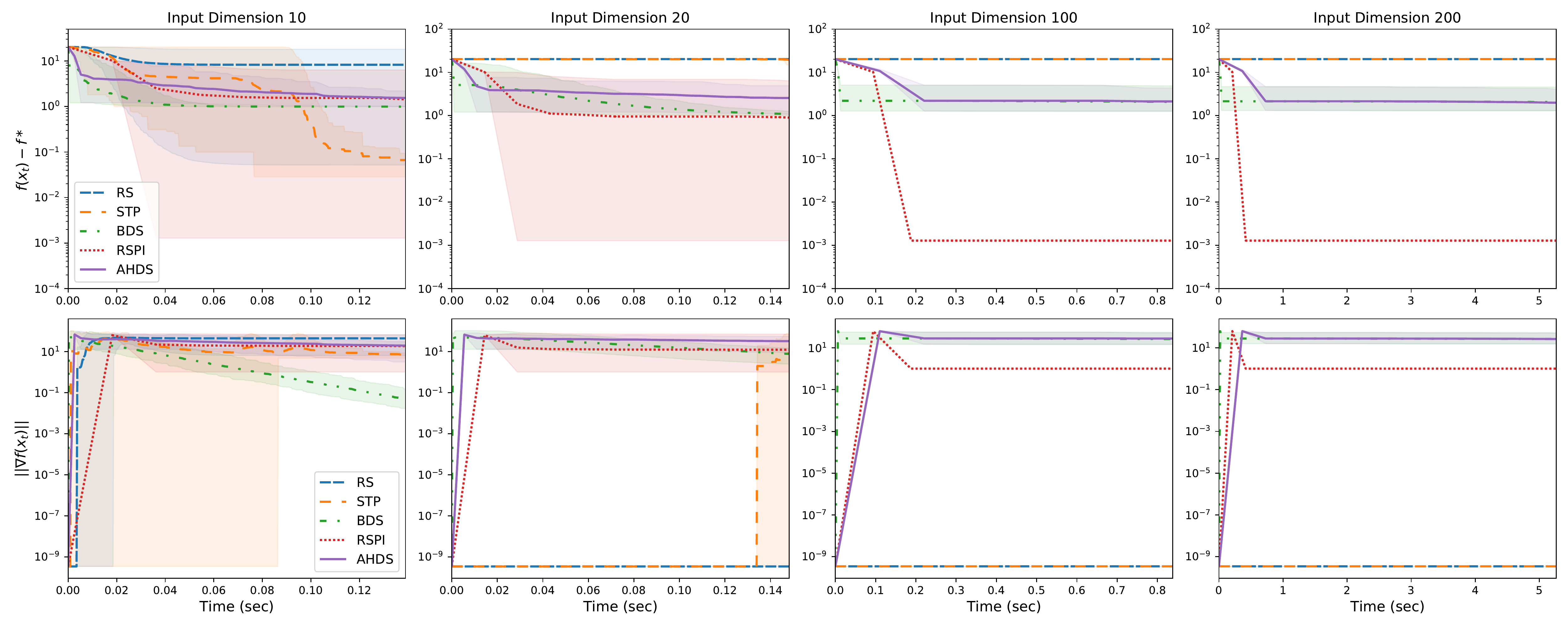}
    \caption{Empirical performance while minimizing the Rastrigin function against wall-clock time. Confidence intervals show min-max intervals over ten runs. All algorithms are initialized at a strict saddle point across all runs.}
\end{figure}

\paragraph{Initialization process.} The critical points of the Rastrigin function satisfy 
\begin{equation}\label{eq:rastrigin-stationary}
    x_i + 10\pi\sin(2\pi x_i) = 0
\end{equation}
for all $i=1,...,d$. The point $\x = \mathbf{0}$ is the unique global minimizer. Stationary points include local minimizers, local maximizers and saddle points. One solution is given by $x_i \approx 0.503$ (truncated to three decimal points). We consider the following initialization 
\begin{equation}
    x_i = 
    \begin{cases*}
        0.503 & if  $i \in \mathcal{I}$  \\
        0 & \textmd{otherwise},
    \end{cases*} 
\end{equation}
where $\mathcal{I}$ is a set of coordinates with cardinality strictly smaller than $d$. In this setup, each non-zero coordinate will be a direction of negative curvature. If we set $\mathcal{I} = \{1,...,d\}$ the point $\x$ is a local maximizer. 

In our experiments we choose $\mathcal{I}$ to have a single coordinate (picked randomly in each experiment repetition). Based on Lemma~\ref{lemma:probability}, we expect that having a single direction of negative curvature will challenge the core mechanism of each algorithm while trying to escape the saddle, especially as the input dimension increases. The results in Figure~\ref{fig:rastrigin-appendix-iter} support our theoretical argument that as the problem dimension grows the probability of sampling a direction that is aligned with the direction of negative curvature decreases exponentially. As a result, RS and STP fail to escape the saddle point for $d=100,200$. 

\begin{figure}
    \centering
    \includegraphics[width=0.6\textwidth]{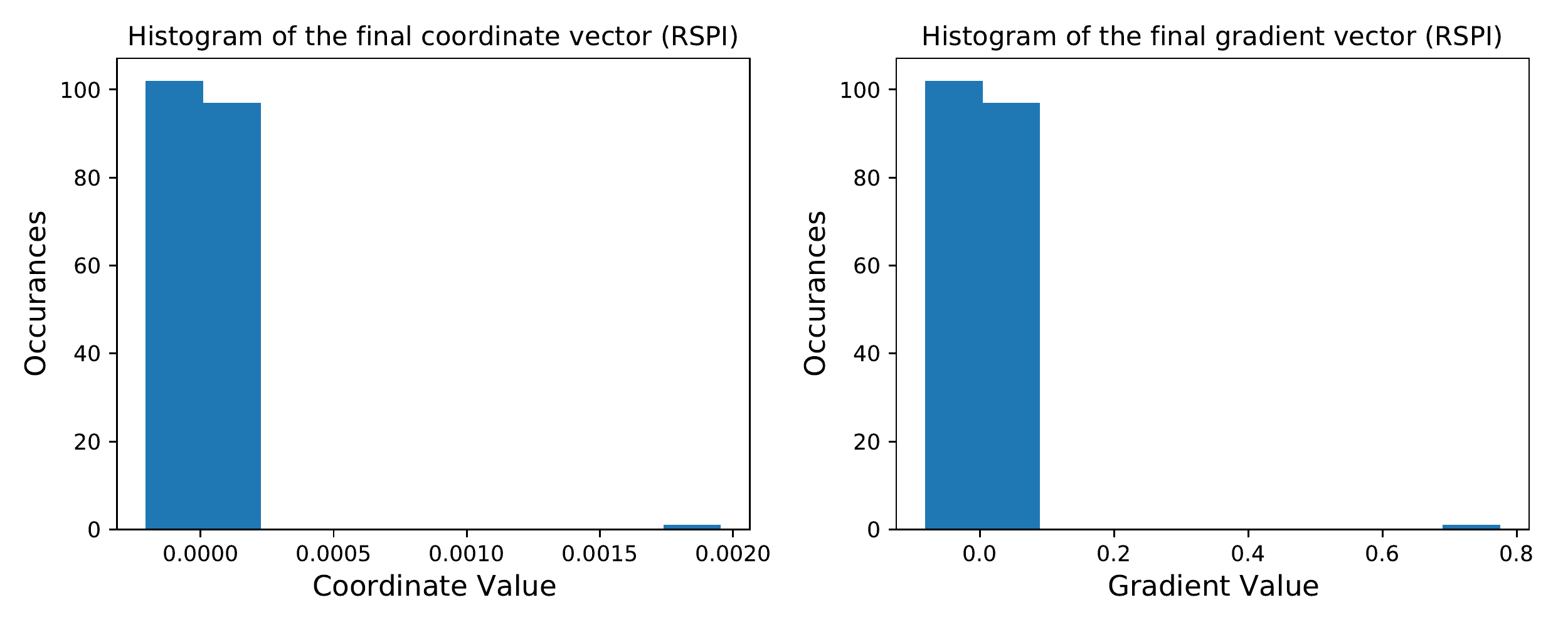}
    \caption{Histogram of the point coordinates (left) and gradient values (right) at the final iterate of RSPI for $d=200$.}
    \label{fig:rastrigin-histogram}
\end{figure}

The rest of the algorithms converge quickly (for $d=100,200$ there is no significant progress after 25 iterations). We speculate that this behaviour is related to the initialization choice. Figure \ref{fig:rastrigin-histogram} shows the distribution of the point coordinates and gradient values at the final iterate of RSPI for $d=200$. In both plots, we observe a cluster of values around zero and a stand-alone component. The later corresponds to the same coordinate that was initialized to non-zero in order to give rise to a saddle point. We observe that the coordinate moves closer to zero (the final coordinate value is less than $0.0020$, whereas the initial value was $0.503$) where the global minimizer occurs. This improvement is achieved through the successful usage of DFPI. That is, RSPI successfully approximates the direction of negative curvature in order to escape the saddle point and move closer to the minimum. Afterwards, no significant progress is achieved via random sampling and that is why the performance curve flattens out after a few iterations. The reason is that in order to achieve further progress via random sampling, it is required to sample a direction that aligns with the single direction of non-zero gradient (see Figure \ref{fig:rastrigin-histogram} (right)) and we expect that probability to decrease exponentially as the dimension increases. That is why further progress can be achieved for $d=10,20$ but not for $d=100,200$. 

\subsection{Leading eigenvector problem}

\begin{figure}[!ht]
    \centering
    \includegraphics[width=0.49\textwidth]{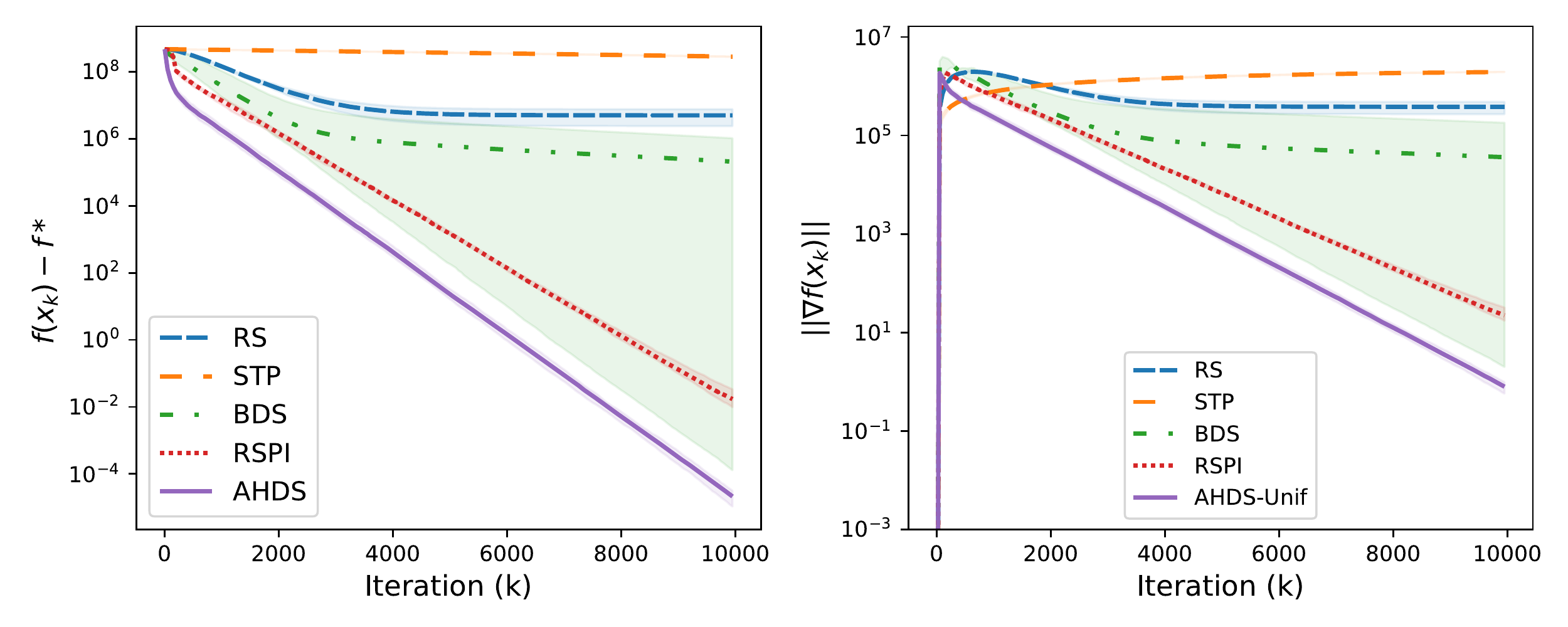}
    \includegraphics[width=0.49\textwidth]{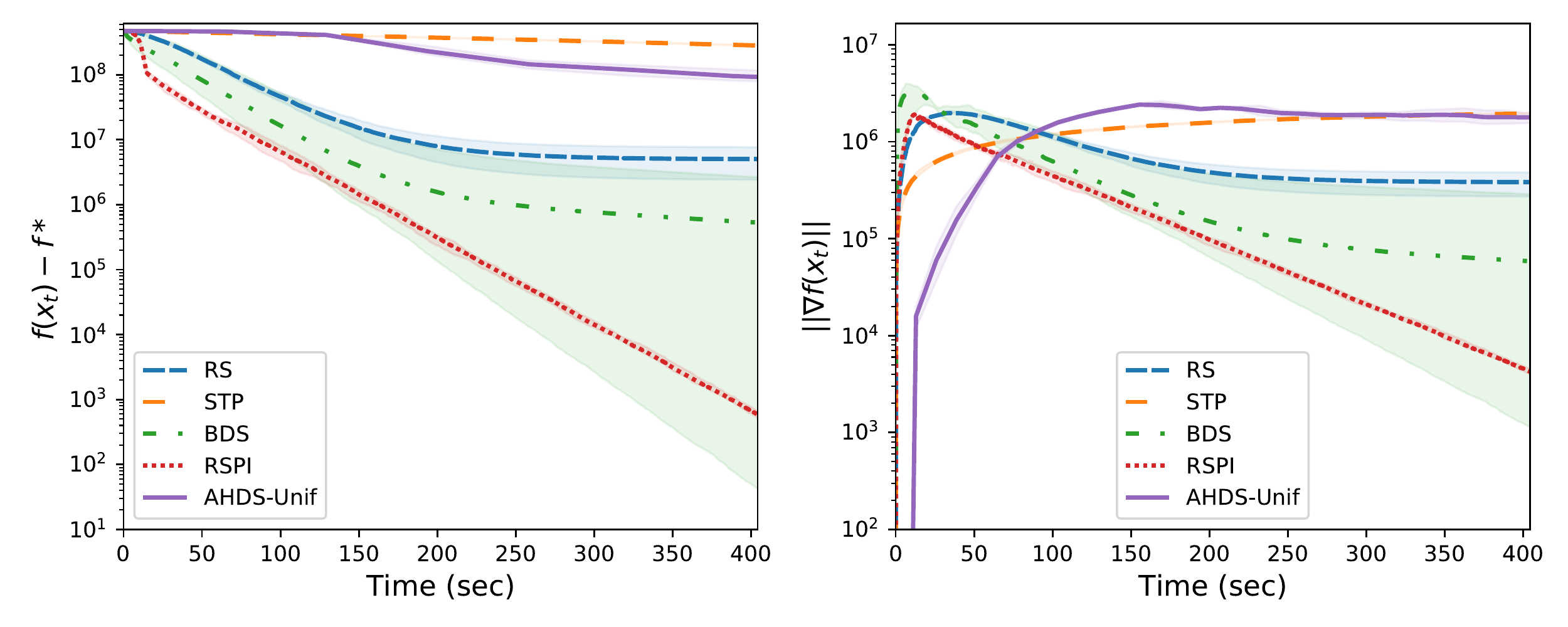}
    \caption{Empirical performance in finding the leading eigenvector of a $350$-dimensional random matrix against iteration and wall-clock time. Confidence intervals show min-max intervals over ten runs. All algorithms are initialized at a strict saddle point across all runs.}
\end{figure}
\newpage
\section{Algorithm Descriptions}

\begin{algorithm}[!ht] 
\begin{algorithmic}[1]
\STATE {\bf{INPUTS : } $\eta_0 \in \mathbb{R}^+$, $\phi : \mathbb{R}^+ \rightarrow \mathbb{R}^+$}
\STATE {Initialize $\x_0$}  
\FOR{$k = 0 \dots K$}
    \STATE $\s_k \sim \S^{d-1}$
    \STATE $\x_{k+1} = \argmin \big\{f(\x_k), f(\x_k + \eta_k \s_k), f(\x_k - \eta_k \s_k)\big\}$
    \STATE $\eta_{k+1} = \phi(\eta_0)$ 
\ENDFOR
\end{algorithmic}
\caption{Stochastic Three Points (STP)}
\label{algo:stp}
\end{algorithm}

\begin{algorithm}[!ht] 
\begin{algorithmic}[1]
\STATE {\bf{INPUTS : } $\eta_{\textmd{max}} > \eta_0 > 0$, $\gamma > 1 > \theta > 0$, $\rho : \mathbb{R}^+ \rightarrow \mathbb{R}^+ $}
\STATE {Set $k=0$ and initialize $\x_0$.}\\
\STATE{Generate a polling set $\mathcal{D}_k$.}\\
\STATE{If there exists $\s_k \in \mathcal{D}_k$ such that 
$$ f(\x_k + \eta_k \s_k) <  f(\x_k) - \rho(\eta_k)$$
then declare the iteration \textbf{successful}, set $\x_{k+1} = \x_k +  \eta_k \s_k$, $\eta_{k+1} = \min\{\gamma \eta_k, \eta_{\textmd{max}}\}$, $k = k+1$ and go to step 4.} \\
\STATE{Otherwise, declare the iteration \textbf{unsuccessful}, set $\x_{k+1} = \x_k$, $\eta_{k+1} = \theta \eta_{k}$, $k = k+1$ and go to step 4.}
\end{algorithmic}
\caption{Basic Direct Search (BDS)}
\label{algo:bds}
\end{algorithm}

\begin{algorithm}[!ht] 
\begin{algorithmic}[1]
\STATE {\bf{INPUTS : } $\eta_{\textmd{max}} > \eta_0 > 0$, $\gamma > 1 > \theta > 0$, $\rho : \mathbb{R}^+ \rightarrow \mathbb{R}^+ $} \\
\STATE {Set $k=0$ and initialize $\x_0$.}\\
\STATE {Generate a polling set $\mathcal{D}_k$. If there exists $\s \in \mathcal{D}_k$ such that
\begin{equation}\label{eq:ahds_suff_dec}
    f(\x_k + \eta_k \s) <  f(\x_k) - \rho(\eta_k)
\end{equation}
then declare iteration $k$ successful with $\s_k = \s$ and go to 8. Otherwise go to 5.} \\
\STATE{If there exists $\s \in \mathcal{D}_k$ such that Eq. (\ref{eq:ahds_suff_dec}) is satisfied with~$-\s$, then declare the iteration successful with $\s_k = -\s$ and go to step 8. Otherwise, go to step 6.} \\
\STATE{Choose $\mathcal{B}_k$ as a subset of $\mathcal{D}_k$ with $d$ linearly independent directions, which we index by $\u_1, ..., \u_d$. If there exists $\s \in \{\u_i + \u_j,  1 \leq i < j \leq d\} $ such that Eq. (\ref{eq:ahds_suff_dec}) holds, then declare the iteration successful with $\s_k = \s$ and go to step 8. Otherwise, go to step 7.} \\
\STATE{Define the Hessian approximation at iteration $k$ as
$$(H_k)_{{i,j}} = \frac{f(\x_k + \eta_k \u_i) - f(\x_k) + f(\x_k - \eta_k \u_i)}{\eta_k^2} \quad \textmd{if } i=j,$$
and 
$$(H_k)_{{i,j}} = \frac{f(\x_k + \eta_k \u_i + \eta_k \u_j) - f(\x_k + \eta_k \u_i) - f(\x_k + \eta_k \u_j) + f(\x_k)}{\eta_k^2} \quad \textmd{if } i<j, $$
for all $i,j \in \{1,...,d\}^2$. Compute a unitary eigenvector $\v_k$ associated with the minimum eigenvalue of $H_k$. If $\v_k$ or $- \v_k$ satisfy the decrease condition in Eq. (\ref{eq:ahds_suff_dec}), then declare the iteration successful with $\s_k$ equal to $\v_k$ or $-\v_k$. Otherwise, declare the iteration unsuccessful and go to step 8.}\\
\STATE{If the iteration was successful, set $\x_{k+1} = \x_k + \eta_k \s_k$ and $\eta_{k+1} = \min\{\gamma \eta_k, \eta_{\textmd{max}}\}$. Otherwise, set $\x_{k+1} = \x_{k}$ and $\eta_{k+1} = \theta \eta_k$.} \\
\STATE{Increment $k$ and go to step 4.} 
\end{algorithmic}
\caption{Approximate Hessian Direct Search (AHDS)}
\label{algo:ahds}
\end{algorithm}

\clearpage
\section{Hyperparameter selection}\label{sec:app-hyperparameters}
For all tasks, the hyperparameters of each method are selected based on a coarse grid search procedure that is refined heuristically by trial and error. The hyperparameters of RS and RSPI are initialized and updated in the same manner, hence the only difference between the two is that RSPI explicitly extracts negative curvature whereas the two-step RS samples a direction at random. In our experiments, we keep $\sigma_2$ constant and only update $\sigma_1$ every $T_{\sigma_1} \in \mathbb{Z}^+$ iterations using the update rule $\sigma_1 \leftarrow \rho \sigma_1$ where $\rho \in (0,1)$. The parameters $\rho$ and $T_{\sigma_1}$ are also selected based on a coarse grid search. We run DFPI for $20$ iterations for all the results shown in the paper and we clarify in the following tables whether Finite Differences (DFPI-FD) or SPSA (DFPI-SPSA) is used to approximate the gradient evaluations within DFPI. 

We illustrate the effect that some crucial parameters have on the performance of the two-step Random Search and the Random Search PI algorithms. In the following figures, confidence intervals show min-max intervals across five runs. All algorithms are initialized at the strict saddle point of the objective 
\begin{equation}
    f(x_1, \cdots, x_d, y) = \frac{1}{4}\sum_{i=1}^d x_i^4 - y\sum_{i=1}^d x_i + \frac{d}{2}y^2.
    \label{eq:objective-saddle}
\end{equation}

\begin{figure*}[!h]
    \centering
    \includegraphics[width=\textwidth]{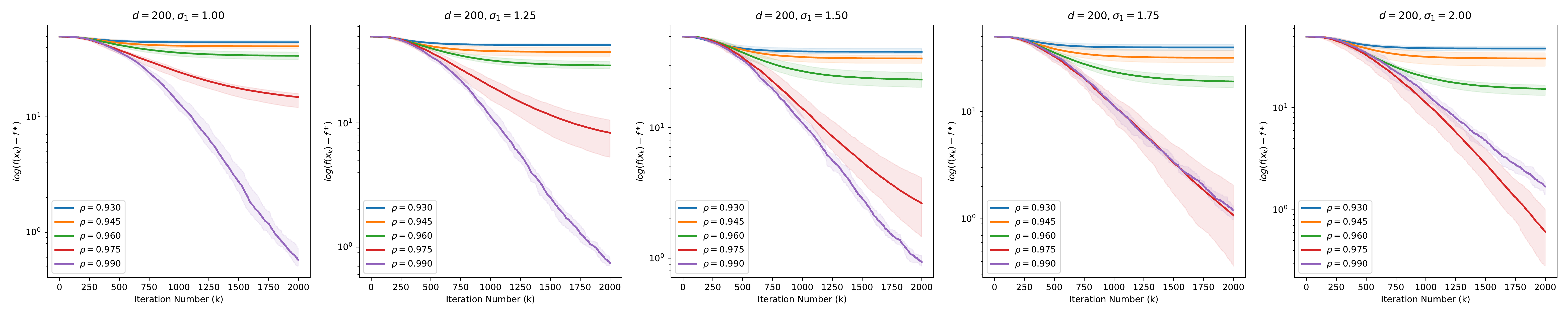}
    \label{fig:sb-params-d200-rspi}
    \caption{Empirical behaviour of the vanilla RS algorithm while minimizing the objective defined in Eq. \ref{eq:objective-saddle} for $d=200$ across different settings of the pair of parameters $(\sigma_1, \rho)$. The parameter $T_{\sigma_1}$ is fixed to $10$.}
\end{figure*} 

\begin{figure*}[!h]
    \centering
    \includegraphics[width=\textwidth]{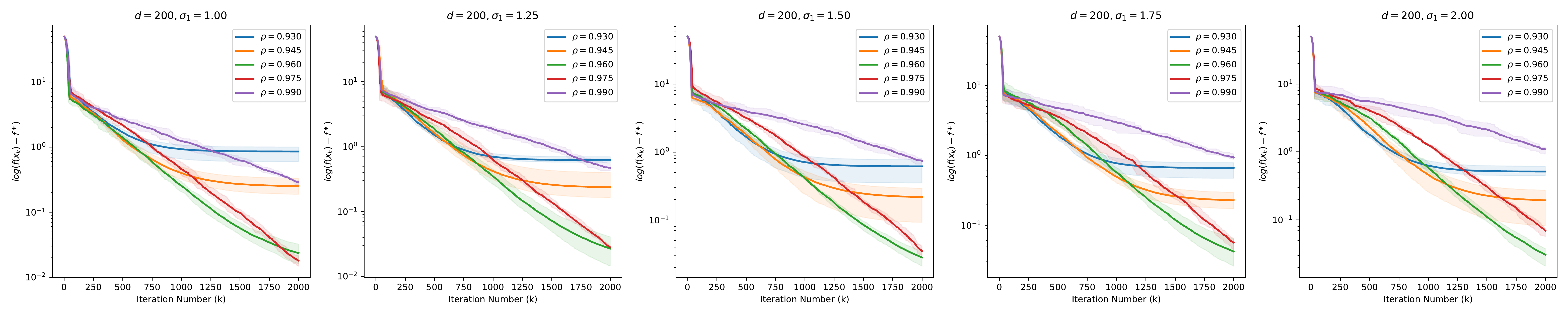}
    \label{fig:sb-params-d200-rspi-2}
    \caption{Empirical behaviour of RSPI while minimizing the objective defined in Eq. \ref{eq:objective-saddle} for $d=200$ across different settings of the pair of parameters $(\sigma_1, \rho)$. The parameter $T_{\sigma_1}$ is fixed to $10$.}
\end{figure*} 

\begingroup
\begin{table}[!h]
  \label{table:leading-eig}
  \caption{Hyperparameters for the leading eigenvector task.}
  \centering
  \renewcommand{\arraystretch}{1.25}
  \begin{tabular}{cc}
    \toprule
    \textbf{Method} & \textbf{Parameters} \\
    \midrule
    \multicolumn{2}{c}{$\bf{d=350}$} \\
    \midrule
    RS & $\sigma_1 = 9.25$, $\sigma_2 = 4.5$, $\rho=0.97$, $T_{\sigma_1} = 25, T_{\text{DFPI}} = 20$      \\
    \hline
    RSPI & $\sigma_1 = 9.25$, $\sigma_2 = 4.5$, $\rho=0.97$, $T_{\sigma_1} = 25, T_{\text{DFPI}} = 20$, DFPI-SPSA      \\
    \hline
    BDS & $\eta_0 = 5.8$, $\eta_{max} = 35.0$, $\gamma=1.25$, $\theta=0.5$, $\rho(x) = 0 $ \\
    \hline
    AHDS & $\eta_0 = 5.8$, $\eta_{max} = 35.0$, $\gamma=1.25$, $\theta=0.5$, $\rho(x) = 0$  \\
    \bottomrule
  \end{tabular}
\end{table}
\endgroup

\begingroup
\begin{table}[!ht]
  \caption{Hyperparameters for the objective in Eq. (\ref{eq:objective-saddle}).}
  \centering
  \renewcommand{\arraystretch}{1.25}
  \begin{tabular}{cc}
    \toprule
    \textbf{Method} & \textbf{Parameters} \\
    \midrule
    \multicolumn{2}{c}{$\bf{d=5}$} \\
    \midrule 
    RS & $\sigma_1 = 1.8$, $\sigma_2 = 0.65$, $\rho=0.6$, $T_{\sigma_1} = 10, T_{\text{DFPI}} = 20$ \\
    \hline
    RSPI & $\sigma_1 = 1.8$, $\sigma_2 = 0.65$, $\rho=0.6$, $T_{\sigma_1} = 10, T_{\text{DFPI}} = 20$, DFPI-SPSA     \\
    \hline
    STP & $\eta_0 = 2.5$,  $\phi(\eta_k) = 0.5\eta_k$ if $k \equiv \mod 10$ (every $10$ iterations)     \\
    \hline
    BDS & $\eta_0 = 0.8$, $\eta_{max} = 10.0$, $\gamma=1.25$, $\theta=0.5$, $\rho(x) = 0 $ \\
    \hline
    AHDS & $\eta_0 = 0.8$, $\eta_{max} = 10.0$, $\gamma=1.25$, $\theta=0.5$, $\rho(x) = 0$  \\
    \midrule
    \multicolumn{2}{c}{$\bf{d=20}$} \\
    \midrule
    RS & $\sigma_1 = 1.75$, $\sigma_2 = 0.65$, $\rho=0.78$, $T_{\sigma_1} = 15, T_{\text{DFPI}} = 20$      \\
    \hline
    RSPI & $\sigma_1 = 1.75$, $\sigma_2 = 0.65$, $\rho=0.78$, $T_{\sigma_1} = 15, T_{\text{DFPI}} = 20$, DFPI-SPSA       \\
    \hline
    STP & $\eta_0 = 2.5$,  $\phi(\eta_k) = 0.5\eta_k$ if $k \equiv \mod 10$ (every $10$ iterations)     \\
    \hline
    BDS & $\eta_0 = 0.8$, $\eta_{max} = 10.0$, $\gamma=1.25$, $\theta=0.5$, $\rho(x) = 0 $ \\
    \hline
    AHDS & $\eta_0 = 0.8$, $\eta_{max} = 10.0$, $\gamma=1.25$, $\theta=0.5$, $\rho(x) = 0$  \\
    \midrule
    \multicolumn{2}{c}{$\bf{d=100}$} \\
    \midrule
    RS & $\sigma_1 = 1.0$, $\sigma_2 = 0.65$, $\rho=0.95$, $T_{\sigma_1} = 15, T_{\text{DFPI}} = 20$      \\
    \hline
    RSPI & $\sigma_1 = 1.0$, $\sigma_2 = 0.65$, $\rho=0.95$, $T_{\sigma_1} = 15, T_{\text{DFPI}} = 20$, DFPI-SPSA      \\
    \hline
    STP & $\eta_0 = 2.5$,  $\phi(\eta_k) = 0.5\eta_k$ if $k \equiv \mod 10$ (every $10$ iterations)     \\
    \hline
    BDS & $\eta_0 = 5.0$, $\eta_{max} = 20.0$, $\gamma=1.25$, $\theta=0.5$, $\rho(x) = 0 $ \\
    \hline
    AHDS & $\eta_0 = 5.0$, $\eta_{max} = 20.0$, $\gamma=1.25$, $\theta=0.5$, $\rho(x) = 0$  \\
    \midrule
    \multicolumn{2}{c}{$\bf{d=200}$} \\
    \midrule
    RS & $\sigma_1 = 1.75$, $\sigma_2 = 0.65$, $\rho=0.96$, $T_{\sigma_1} = 15, T_{\text{DFPI}} = 20$      \\
    \hline
    RSPI & $\sigma_1 = 1.75$, $\sigma_2 = 0.65$, $\rho=0.96$, $T_{\sigma_1} = 15, T_{\text{DFPI}} = 20$, DFPI-SPSA       \\
    \hline
    STP & $\eta_0 = 2.5$,  $\phi(\eta_k) = 0.5\eta_k$ if $k \equiv \mod 10$ (every $10$ iterations)     \\
    \hline
    BDS & $\eta_0 = 5.0$, $\eta_{max} = 20.0$, $\gamma=1.25$, $\theta=0.5$, $\rho(x) = 0 $ \\
    \hline
    AHDS & $\eta_0 = 5.0$, $\eta_{max} = 20.0$, $\gamma=1.25$, $\theta=0.5$, $\rho(x) = 0$  \\
    \bottomrule
  \end{tabular}
\end{table}
\endgroup

\begingroup
\begin{table}[!ht]
  \caption{Hyperparameters for the Rastrigin function.}
  \centering
  \renewcommand{\arraystretch}{1.25}
  \begin{tabular}{cc}
    \toprule
    \textbf{Method} & \textbf{Parameters} \\
    \midrule
    \multicolumn{2}{c}{$\bf{d=10}$} \\
    \midrule 
    RS & $\sigma_1 = 0.25$, $\sigma_2 = 0.25$, $\rho=0.83$, $T_{\sigma_1} = 5, T_{\text{DFPI}} = 20$\\
    \hline
    RSPI & $\sigma_1 = 0.25$, $\sigma_2 = 0.25$, $\rho=0.83$, $T_{\sigma_1} = 5, T_{\text{DFPI}} = 20$, DFPI-FD    \\
    \hline
    STP & $\eta_0 = 0.25$,  $\phi(\eta_k) = \eta_0/\sqrt{k+1}$    \\
    \hline
    BDS & $\eta_0 = 0.25$, $\eta_{max} = 10.0$, $\gamma=1.1$, $\theta=0.9$, $\rho(x) = 0 $ \\
    \hline
    AHDS & $\eta_0 = 0.25$, $\eta_{max} = 10.0$, $\gamma=1.1$, $\theta=0.9$, $\rho(x) = 0$  \\
    \midrule
    \multicolumn{2}{c}{$\bf{d=20}$} \\
    \midrule
    RS & $\sigma_1 = 0.255$, $\sigma_2 = 0.25$, $\rho=0.83$, $T_{\sigma_1} = 5, T_{\text{DFPI}} = 20$      \\
    \hline
    RSPI & $\sigma_1 = 0.255$, $\sigma_2 = 0.25$, $\rho=0.83$, $T_{\sigma_1} = 5, T_{\text{DFPI}} = 20$, DFPI-FD      \\
    \hline
    STP & $\eta_0 = 0.25$,  $\phi(\eta_k) = \eta_0/\sqrt{k+1}$    \\
    \hline
    BDS & $\eta_0 = 0.25$, $\eta_{max} = 10.0$, $\gamma=1.1$, $\theta=0.9$, $\rho(x) = 0 $ \\
    \hline
    AHDS & $\eta_0 = 0.25$, $\eta_{max} = 10.0$, $\gamma=1.1$, $\theta=0.9$, $\rho(x) = 0$  \\
    \midrule
    \multicolumn{2}{c}{$\bf{d=100}$} \\
    \midrule
    RS & $\sigma_1 = 0.15$, $\sigma_2 = 0.25$, $\rho=0.83$, $T_{\sigma_1} = 5, T_{\text{DFPI}} = 20$      \\
    \hline
    RSPI & $\sigma_1 = 0.15$, $\sigma_2 = 0.25$, $\rho=0.83$, $T_{\sigma_1} = 5, T_{\text{DFPI}} = 20$, DFPI-FD      \\
    \hline
    STP & $\eta_0 = 0.25$,  $\phi(\eta_k) = \eta_0/\sqrt{k+1}$    \\
    \hline
    BDS & $\eta_0 = 0.25$, $\eta_{max} = 10.0$, $\gamma=1.1$, $\theta=0.9$, $\rho(x) = 0 $ \\
    \hline
    AHDS & $\eta_0 = 0.25$, $\eta_{max} = 10.0$, $\gamma=1.1$, $\theta=0.9$, $\rho(x) = 0$  \\
    \midrule
    \multicolumn{2}{c}{$\bf{d=200}$} \\
    \midrule
    RS & $\sigma_1 = 0.15$, $\sigma_2 = 0.25$, $\rho=0.83$, $T_{\sigma_1} = 5, T_{\text{DFPI}} = 20$      \\
    \hline
    RSPI & $\sigma_1 = 0.15$, $\sigma_2 = 0.25$, $\rho=0.83$, $T_{\sigma_1} = 5, T_{\text{DFPI}} = 20$, DFPI-FD      \\
    \hline
    STP & $\eta_0 = 0.25$,  $\phi(\eta_k) = \eta_0/\sqrt{k+1}$    \\
    \hline
    BDS & $\eta_0 = 0.25$, $\eta_{max} = 10.0$, $\gamma=1.1$, $\theta=0.9$, $\rho(x) = 0 $ \\
    \hline
    AHDS & $\eta_0 = 0.25$, $\eta_{max} = 10.0$, $\gamma=1.1$, $\theta=0.9$, $\rho(x) = 0$  \\
    \bottomrule
  \end{tabular}
\end{table}
\endgroup

%% file: main.bbl
\begin{thebibliography}{10}

\bibitem{abramson2005second}
Mark~A Abramson.
\newblock Second-order behavior of pattern search.
\newblock {\em SIAM Journal on Optimization}, 16(2):515--530, 2005.

\bibitem{abramson2014subclass}
Mark~A Abramson, Lennart Frimannslund, and Trond Steihaug.
\newblock A subclass of generating set search with convergence to second-order
  stationary points.
\newblock {\em Optimization Methods and Software}, 29(5):900--918, 2014.

\bibitem{anagnostidis2021direct}
Sotirios-Konstantinos Anagnostidis, Aurelien Lucchi, and Youssef Diouane.
\newblock Direct-search for a class of stochastic min-max problems.
\newblock In {\em International Conference on Artificial Intelligence and
  Statistics}, pages 3772--3780. PMLR, 2021.

\bibitem{arjevani2020second}
Yossi Arjevani, Yair Carmon, John~C Duchi, Dylan~J Foster, Ayush Sekhari, and
  Karthik Sridharan.
\newblock Second-order information in non-convex stochastic optimization: Power
  and limitations.
\newblock In {\em Conference on Learning Theory}, pages 242--299. PMLR, 2020.

\bibitem{balcan2016improved}
Maria-Florina Balcan, Simon~Shaolei Du, Yining Wang, and Adams~Wei Yu.
\newblock An improved gap-dependency analysis of the noisy power method.
\newblock In {\em Conference on Learning Theory}, pages 284--309. PMLR, 2016.

\bibitem{bergou2019stochastic}
El~Houcine Bergou, Eduard Gorbunov, and Peter Richtarik.
\newblock Stochastic three points method for unconstrained smooth minimization.
\newblock {\em SIAM Journal on Optimization}, 2020.

\bibitem{bergstra2012random}
James Bergstra and Yoshua Bengio.
\newblock Random search for hyper-parameter optimization.
\newblock {\em Journal of machine learning research}, 13(2), 2012.

\bibitem{bisong2019colab}
Ekaba Bisong.
\newblock {\em Google Colaboratory}, pages 59--64.
\newblock Apress, Berkeley, CA, 2019.

\bibitem{bubeck2012regret}
S{\'e}bastien Bubeck and Nicolo Cesa-Bianchi.
\newblock Regret analysis of stochastic and nonstochastic multi-armed bandit
  problems.
\newblock {\em arXiv preprint arXiv:1204.5721}, 2012.

\bibitem{carmon2018accelerated}
Yair Carmon, John~C Duchi, Oliver Hinder, and Aaron Sidford.
\newblock Accelerated methods for nonconvex optimization.
\newblock {\em SIAM Journal on Optimization}, 28(2):1751--1772, 2018.

\bibitem{chen2017zoo}
Pin-Yu Chen, Huan Zhang, Yash Sharma, Jinfeng Yi, and Cho-Jui Hsieh.
\newblock Zoo: Zeroth order optimization based black-box attacks to deep neural
  networks without training substitute models.
\newblock In {\em Proceedings of the 10th ACM Workshop on Artificial
  Intelligence and Security}, pages 15--26, 2017.

\bibitem{chen2020accelerated}
Yuwen Chen, Antonio Orvieto, and Aurelien Lucchi.
\newblock An accelerated dfo algorithm for finite-sum convex functions.
\newblock In {\em Proceedings of the 30th International Conference on Neural
  Information Processing Systems}, 2020.

\bibitem{cheung2001generalizations}
Wing-Sum Cheung.
\newblock Generalizations of {H{\"o}lder}'s inequality.
\newblock {\em International Journal of Mathematics and Mathematical Sciences},
  26, 2001.

\bibitem{choromanska2015loss}
Anna Choromanska, Mikael Henaff, Michael Mathieu, G{\'e}rard~Ben Arous, and
  Yann LeCun.
\newblock The loss surfaces of multilayer networks.
\newblock In {\em Artificial intelligence and statistics}, pages 192--204,
  2015.

\bibitem{conn2009introduction}
Andrew~R Conn, Katya Scheinberg, and Luis~N Vicente.
\newblock {\em Introduction to derivative-free optimization}.
\newblock SIAM, 2009.

\bibitem{daneshmand2018escaping}
Hadi Daneshmand, Jonas Kohler, Aurelien Lucchi, and Thomas Hofmann.
\newblock Escaping saddles with stochastic gradients.
\newblock {\em arXiv preprint arXiv:1803.05999}, 2018.

\bibitem{davis1970rotation}
Chandler Davis and William~Morton Kahan.
\newblock The rotation of eigenvectors by a perturbation. iii.
\newblock {\em SIAM Journal on Numerical Analysis}, 7(1):1--46, 1970.

\bibitem{de2018accelerated}
Christopher De~Sa, Bryan He, Ioannis Mitliagkas, Christopher R{\'e}, and Peng
  Xu.
\newblock Accelerated stochastic power iteration.
\newblock {\em Proceedings of machine learning research}, 84:58, 2018.

\bibitem{du2017gradient}
Simon~S Du, Chi Jin, Jason~D Lee, Michael~I Jordan, Aarti Singh, and Barnabas
  Poczos.
\newblock Gradient descent can take exponential time to escape saddle points.
\newblock In {\em Advances in neural information processing systems}, pages
  1067--1077, 2017.

\bibitem{duchi2015optimal}
John~C Duchi, Michael~I Jordan, Martin~J Wainwright, and Andre Wibisono.
\newblock Optimal rates for zero-order convex optimization: The power of two
  function evaluations.
\newblock {\em IEEE Transactions on Information Theory}, 61(5):2788--2806,
  2015.

\bibitem{ge2017no}
Rong Ge, Chi Jin, and Yi~Zheng.
\newblock No spurious local minima in nonconvex low rank problems: A unified
  geometric analysis.
\newblock In {\em International Conference on Machine Learning}, pages
  1233--1242. PMLR, 2017.

\bibitem{ge2016matrix}
Rong Ge, Jason~D Lee, and Tengyu Ma.
\newblock Matrix completion has no spurious local minimum.
\newblock In {\em Proceedings of the 30th International Conference on Neural
  Information Processing Systems}, pages 2981--2989, 2016.

\bibitem{ghadimi2013stochastic}
Saeed Ghadimi and Guanghui Lan.
\newblock Stochastic first-and zeroth-order methods for nonconvex stochastic
  programming.
\newblock {\em SIAM Journal on Optimization}, 23(4):2341--2368, 2013.

\bibitem{golovin2019gradientless}
Daniel Golovin, John Karro, Greg Kochanski, Chansoo Lee, Xingyou Song, et~al.
\newblock Gradientless descent: High-dimensional zeroth-order optimization.
\newblock {\em arXiv preprint arXiv:1911.06317}, 2019.

\bibitem{golub2013matrix}
Gene~H Golub and Charles~F Van~Loan.
\newblock Matrix computations, 4th.
\newblock {\em Johns Hopkins}, 2013.

\bibitem{gratton2016second}
Serge Gratton, CW~Royer, and Luis~Nunes Vicente.
\newblock A second-order globally convergent direct-search method and its
  worst-case complexity.
\newblock {\em Optimization}, 65(6):1105--1128, 2016.

\bibitem{gu2015subspace}
Ming Gu.
\newblock Subspace iteration randomization and singular value problems.
\newblock {\em SIAM Journal on Scientific Computing}, 37(3):A1139--A1173, 2015.

\bibitem{hansen2009rastrigin}
Nikolaus Hansen, Steffen Finck, Raymond Ros, and Anne Auger.
\newblock {Real-Parameter Black-Box Optimization Benchmarking 2009: Noiseless
  Functions Definitions}.
\newblock Research Report RR-6829, {INRIA}, 2009.

\bibitem{hardt2014noisy}
Moritz Hardt and Eric Price.
\newblock The noisy power method: A meta algorithm with applications.
\newblock In {\em Advances in Neural Information Processing Systems}, pages
  2861--2869, 2014.

\bibitem{hopcroft2012computer}
John Hopcroft and Ravi Kannan.
\newblock {\em Lecture notes: Computer science theory for the information age}.
\newblock Citeseer, 2012.

\bibitem{ji2019improved}
Kaiyi Ji, Zhe Wang, Yi~Zhou, and Yingbin Liang.
\newblock Improved zeroth-order variance reduced algorithms and analysis for
  nonconvex optimization.
\newblock {\em arXiv preprint arXiv:1910.12166}, 2019.

\bibitem{jin2017escape}
Chi Jin, Rong Ge, Praneeth Netrapalli, Sham~M Kakade, and Michael~I Jordan.
\newblock How to escape saddle points efficiently.
\newblock In {\em Proceedings of the 34th International Conference on Machine
  Learning-Volume 70}, pages 1724--1732. JMLR. org, 2017.

\bibitem{jin2019nonconvex}
Chi Jin, Praneeth Netrapalli, Rong Ge, Sham~M Kakade, and Michael~I Jordan.
\newblock On nonconvex optimization for machine learning: Gradients,
  stochasticity, and saddle points.
\newblock {\em arXiv preprint arXiv:1902.04811}, 2019.

\bibitem{Jin2019leadingeigval}
Chi Jin, Praneeth Netrapalli, Rong Ge, Sham~M. Kakade, and Michael~I. Jordan.
\newblock On nonconvex optimization for machine learning: Gradients,
  stochasticity, and saddle points.
\newblock {\em arXiv: Learning}, 2019.

\bibitem{jin2018accelerated}
Chi Jin, Praneeth Netrapalli, and Michael~I Jordan.
\newblock Accelerated gradient descent escapes saddle points faster than
  gradient descent.
\newblock In {\em Conference On Learning Theory}, pages 1042--1085. PMLR, 2018.

\bibitem{kawaguchi2016deep}
Kenji Kawaguchi.
\newblock Deep learning without poor local minima.
\newblock In {\em Proceedings of the 30th International Conference on Neural
  Information Processing Systems}, pages 586--594, 2016.

\bibitem{kuczynski1992estimating}
Jacek Kuczy{\'n}ski and Henryk Wo{\'z}niakowski.
\newblock Estimating the largest eigenvalue by the power and {Lanczos}
  algorithms with a random start.
\newblock {\em SIAM journal on matrix analysis and applications},
  13(4):1094--1122, 1992.

\bibitem{larson2019derivative}
Jeffrey Larson, Matt Menickelly, and Stefan~M Wild.
\newblock Derivative-free optimization methods.
\newblock {\em arXiv preprint arXiv:1904.11585}, 2019.

\bibitem{levy2016power}
Kfir~Y Levy.
\newblock The power of normalization: Faster evasion of saddle points.
\newblock {\em arXiv preprint arXiv:1611.04831}, 2016.

\bibitem{lewis2000direct}
Robert~Michael Lewis, Virginia Torczon, and Michael~W Trosset.
\newblock Direct search methods: then and now.
\newblock {\em Journal of computational and Applied Mathematics},
  124(1-2):191--207, 2000.

\bibitem{liu2017noisy}
Mingrui Liu and Tianbao Yang.
\newblock On noisy negative curvature descent: Competing with gradient descent
  for faster non-convex optimization.
\newblock {\em arXiv preprint arXiv:1709.08571}, 2017.

\bibitem{maheswaranathan2019guided}
Niru Maheswaranathan, Luke Metz, George Tucker, Dami Choi, and Jascha
  Sohl-Dickstein.
\newblock Guided evolutionary strategies: Augmenting random search with
  surrogate gradients.
\newblock In {\em International Conference on Machine Learning}, pages
  4264--4273. PMLR, 2019.

\bibitem{mania2018simple}
Horia Mania, Aurelia Guy, and Benjamin Recht.
\newblock Simple random search provides a competitive approach to reinforcement
  learning.
\newblock {\em arXiv preprint arXiv:1803.07055}, 2018.

\bibitem{matyas1965random}
J~Matyas.
\newblock Random optimization.
\newblock {\em Automation and Remote control}, 26(2):246--253, 1965.

\bibitem{musco2015randomized}
Cameron Musco and Christopher Musco.
\newblock Randomized block {Krylov} methods for stronger and faster approximate
  singular value decomposition.
\newblock In {\em Proceedings of the 28th International Conference on Neural
  Information Processing Systems-Volume 1}, pages 1396--1404, 2015.

\bibitem{nelder1965simplex}
John~A Nelder and Roger Mead.
\newblock A simplex method for function minimization.
\newblock {\em The computer journal}, 7(4):308--313, 1965.

\bibitem{nesterov2018lectures}
Yurii Nesterov et~al.
\newblock {\em Lectures on convex optimization}, volume 137.
\newblock Springer, 2018.

\bibitem{nesterov2017random}
Yurii Nesterov and Vladimir Spokoiny.
\newblock Random gradient-free minimization of convex functions.
\newblock {\em Foundations of Computational Mathematics}, 17(2):527--566, 2017.

\bibitem{oja1985stochastic}
Erkki Oja and Juha Karhunen.
\newblock On stochastic approximation of the eigenvectors and eigenvalues of
  the expectation of a random matrix.
\newblock {\em Journal of mathematical analysis and applications},
  106(1):69--84, 1985.

\bibitem{rudelson2009smallest}
Mark Rudelson and Roman Vershynin.
\newblock Smallest singular value of a random rectangular matrix.
\newblock {\em Communications on Pure and Applied Mathematics: A Journal Issued
  by the Courant Institute of Mathematical Sciences}, 62(12):1707--1739, 2009.

\bibitem{salimans2017evolution}
Tim Salimans, Jonathan Ho, Xi~Chen, Szymon Sidor, and Ilya Sutskever.
\newblock Evolution strategies as a scalable alternative to reinforcement
  learning.
\newblock {\em arXiv preprint arXiv:1703.03864}, 2017.

\bibitem{shamir2015stochastic}
Ohad Shamir.
\newblock A stochastic {PCA} and {SVD} algorithm with an exponential
  convergence rate.
\newblock In {\em International Conference on Machine Learning}, pages
  144--152, 2015.

\bibitem{spall1992multivariate}
James~C Spall et~al.
\newblock Multivariate stochastic approximation using a simultaneous
  perturbation gradient approximation.
\newblock {\em IEEE transactions on automatic control}, 37(3):332--341, 1992.

\bibitem{vicente2013worst}
Lu{\'\i}s~Nunes Vicente.
\newblock Worst case complexity of direct search.
\newblock {\em EURO Journal on Computational Optimization}, 1(1-2):143--153,
  2013.

\bibitem{vlatakis2019efficiently}
Emmanouil-Vasileios Vlatakis-Gkaragkounis, Lampros Flokas, and Georgios
  Piliouras.
\newblock Efficiently avoiding saddle points with zero order methods: No
  gradients required.
\newblock In {\em Advances in Neural Information Processing Systems 32: Annual
  Conference on Neural Information Processing Systems 2019}, 2019.

\bibitem{xu2018first}
Yi~Xu, Rong Jin, and Tianbao Yang.
\newblock First-order stochastic algorithms for escaping from saddle points in
  almost linear time.
\newblock {\em Advances in Neural Information Processing Systems},
  31:5530--5540, 2018.

\bibitem{yang2020devil}
Yingxiang Yang, Negar Kiyavash, Le~Song, and Niao He.
\newblock The devil is in the detail: A framework for macroscopic prediction
  via microscopic models.
\newblock {\em Advances in Neural Information Processing Systems}, 33, 2020.

\bibitem{ye2018hessian}
Haishan Ye, Zhichao Huang, Cong Fang, Chris~Junchi Li, and Tong Zhang.
\newblock Hessian-aware zeroth-order optimization for black-box adversarial
  attack.
\newblock {\em arXiv preprint arXiv:1812.11377}, 2018.

\bibitem{zabinsky2013stochastic}
Zelda~B Zabinsky.
\newblock {\em Stochastic adaptive search for global optimization}, volume~72.
\newblock Springer Science \& Business Media, 2013.

\end{thebibliography}
